\documentclass[a4paper,oneside]{memoir}
\usepackage[utf8]{inputenc}
\usepackage{amssymb,amsmath,latexsym}
\usepackage[english]{babel}

\usepackage{tikz-cd}

\begin{document}

\title{Points Of Small Heights in Certain Nonabelian Extensions}
\date{}
\author{Soumyadip Sahu}
\begin{titlingpage}

    \begin{center}
        \vspace*{1cm}
        
        \huge
        \textbf{Points Of Small Heights in Certain Nonabelian Extensions\footnote{Updated version}}
        
        \vspace{1.5cm}
        \large
        \textbf{Soumyadip Sahu}\\
        MMA201605\\
        \textbf{Adviser :} Prof. Sinnou David
        
        \vfill
        
        A thesis presented for the degree of\\
        M.Sc. in Mathematics
        
        \vspace{0.8cm}

        \Large
        Chennai Mathematical Institute\\
        2018

    \end{center}

\thispagestyle{empty}
\end{titlingpage}
\setcounter{chapter}{-1}
\newpage
\begin{abstract}
Let $E$ be an elliptic curve without complex multiplication defined over a number field $K$ which has at least one real embedding. The field $F$ generated by all torsion points of $E$ over $K$ is an infinite, non-abelian Galois extension of the ground field which has unbounded, wild ramification above all primes. Following the treatment in \cite {hab} we prove that the absolute logarithmic Weil height of an element of $F$ is either zero or bounded from below by a positive constant depending only on $E$ and $K$. We also show that the N\'eron-Tate height has a similar gap on $E(F)$. \footnote{As of now, the proof of a technical lemma is incomplete (see lemma 2.3.7). Hence strictly speaking, we have proved the results mentioned here assuming the statement of lemma 2.3.7. (appendix B).}\\
In appendix-A we have included some new results about Galois properties of division points of formal groups which are generalizations of results proved in chapter 2.                               
\end{abstract}  
\newpage  

\tableofcontents  

\newpage
\addcontentsline{toc}{chapter}{Notations and conventions}
\chapter*{Notations and conventions}
We shall use $\mathbb{R}, \mathbb{C}, \mathbb{Q}, \mathbb{Z}, \mathbb{N}$ in usual sense. In our convention $0 \notin \mathbb{N}$.\\
For a finite set $X$ we shall use $|X|$ to denote the cardinality of $X$.\\
In different places of the thesis we have introduced notation which are supposed to be valid for that part (chapter, section, subsection etc.) only. If we use them in other parts we indicate it.\\~\\
A number field is a finite extension of $\mathbb{Q}$. An algebraic number is an element of a number field. \\
While referring to these terms we shall fix an algebraic closure of $\mathbb{Q}$ and all the number fields can be considered as subfield of that.\\
\chapter*{Acknowledgements}
I am indebted to Prof. S. David and this can be considered as a joint work with him.\\
I am thankful to Prof. P. Rath for introducing me to Prof. S. David and giving me a chance to work in this project.\\
I am thankful to my pet C-infinity who has been a great support during these tiresome and troubled days.

\chapter{Introduction}
For an algebraic number one can introduce the notion of `height' which in some sense measures its arithmetic complexity (see section 1, chapter 1 of this thesis). Kronecker's theorem states that an algebraic number has height zero if and only if it is zero or a root of unity.\\
By a well known result of Northcott, we know that there are only finitely many algebraic numbers of bounded degree and bounded height. So a nonzero element of a number field that is not a root of unity has height bounded below by a positive real number which depends only on the number field under consideration.\\
A field of algebraic numbers is said to have the Bogolomov property if zero is an isolated point among its height values. By previous remark all the number fields have this property.  \\~\\
Consider the following example :\\
Let $n$ be a positive integer.\\
Define, $a_{n} = 2^{ \frac {1} {n} }$.\\
Now one can easily show, the height of $a_{n}$ (denoted $h(a_{n})$) is given by $\frac { \log 2} {n}$.\\
Hence, $h(a_{n}) \to 0$ as $ n \to \infty$.\\
This proves that the field of all algebraic numbers does not have the Bogolomov property.\\~\\
Some fields that are known to have the Bogolomov property are :\\
i) a number field,\\
ii) maximal abelian extension of a number field (Amoroso and Zannier, see \cite{az}),\\
iii) $\mathbb{Q}^{\text{tr}}$,the maximal totally real extension of $\mathbb{Q}$ (Schinzel, see  \cite{sch}).\\~\\
Let $E$ be an elliptic curve defined over a number field $K$.\\
Let $K(E_{\text{tor}})$ be the field generated by torsion points of $E$.\\
If $E$ has complex multiplication then $K(E_{\text{tor}}) \subseteq K^{\text{ab}}$, the maximal abelian extension of $K$, provided $K$ contains field of complex multiplication. Thus it has Bogolomov property. So we can assume that $E$ does not have complex multiplication.\\
In \cite{hab} Habegger proves that if $K = \mathbb{Q}$ then $K(E_{\text{tor}})$ has Bogolomov property.\\
He further showed that in this case similar discreteness result holds for N\'eron - Tate height (denoted $\hat{h}$) on the elliptic curve $E$. More precisely, there is an $\epsilon > 0$ such that if $ A \in E(K(E_{\text{tor}}))$ is nontorsion point, then $\hat{h}(A) > \epsilon $.\\~\\
The goal of this thesis is to prove the following theorem:\\~\\
\textbf{Theorem 0.1 :} Let $E$ be an elliptic curve  without complex multiplication defined over a number field $K$ such that $K$ has at least one real embedding.\\
Then $K(E_{\text{tor}})$ has Bogolomov property.\\~\\
We also prove an analogue of the result in elliptic case:\\~\\
\textbf{Theorem 0.2 :} $K$ and $E$ as in statement of theorem 0.1.\\
There exists an $\epsilon > 0$  depending only on $K$ and $E$ such that if $ A \in E(K(E_{\text{tor}}))$ is a nontorsion, then $\hat{h}(A) > \epsilon$.\\~\\  
\textbf{Remark 0.3 :} i) The starting point of Habegger's paper is a result due to Elkies which states that every elliptic curve defined over $\mathbb{Q}$ has infinitely many primes of supersingular reduction.\\
In a later paper Elkies has also shown that if $E$ is an elliptic curve over a number field $K$ with at least one real embedding, then it has infinitely many places of supersingular reduction (see chapter 1, section 3).\\
We take this as main input for our arguments presented in the thesis.\\
The proof can be generalized to an elliptic curve over any number field if that curve has infinitely many places of supersingular reduction.\\
Since nothing is known in general for elliptic curves over totally imaginary field (though it is known for some specific curves), we restrict ourselves to the case stated in our theorem.\\
ii) If $E$ has complex multiplication, analogue of theorem 0.2 is consequence of height bound on elliptic curves for abelian extension. (Silverman, \cite{sil 3}).\\

\chapter{Preliminaries}
\section{Heights}
Let $K$ be a number field. A place $v$ of $K$ is an absolute value $|\cdot|_{v} : K \to [ 0 , \infty )$ whose restriction $w$ to $\mathbb{Q}$ is either the standard complex absolute value $ w = \infty$ or $ w = p$, the $p-$adic absolute value for a positive prime number $p$. In the former case we write $ v | \infty$ and call $v$ infinite or Archimedean. In the latter case we write $ v | p$ or in general $ v \nmid \infty$ and call $v$ finite or non-Archimedean. A place is finite if and only if it satisfies the ultrametric triangle inequality. The completion of $K$ with respect to $v$ is denoted by $K_{v}$. We use the same symbol $|\cdot|_{v}$ for the absolute value on $K_v$. The set of finite places can be naturally identified with the set of nonzero prime ideals of the ring of integers of $K$. The infinite places are in bijection with field embeddings $ K \to \mathbb{C}$ upto complex conjugation.\\
Let $v$ be a place of $K$ which restricts to a place $w$ of $\mathbb{Q}$. We define local degree of $v$ as $d_{v} = [ K_{v} : \mathbb{Q}_{w} ]$. \\
For each positive prime $p$ we normalize $|\cdot|_{p}$ by setting $|p|_p = \frac{1}{p}$ . We normalize the finite places of $K$ accordingly.\\~\\
\textbf{Definition 1.1.1 :} Let $K$ be a number field and let $\alpha \in K$. The logarithmic Weil height, or in short height, of $\alpha$ with respect to $K$ is defined to be 
$$h_{K}(\alpha) = \frac{1}{[K : \mathbb{Q}]} \sum_{v} d_{v}\log^{+} (|\alpha|_{v})$$
where $v$ runs over all normalized places of $K$ and $\log^{+}(t) = \log (\max \{ 1, t \})$ for a real number $t$.\\~\\
\textbf{Remark 1.1.2 :} a) If $\alpha = 0$ then $h_{K}(\alpha) = 0$.\\
b) If $\alpha \neq 0$ then for all but finitely many places $|\alpha|_v = 1$. Thus the sum on the right hand side is actually a finite sum.\\
c) Let $\alpha \in \overline{\mathbb{Q}}$ and let $L, K$ be two number fields both containing $\alpha$. Then $h_{L}(\alpha) = h_{K}(\alpha)$. For a proof of this the readers are referred to \cite{heights}, chapter 1.\\
As a consequence of this statement we can drop the subscript $K$ from the definition and define height uniformly over $\overline{\mathbb{Q}}$.\\~\\
\textbf{Proposition 1.1.3 :} a) Let $\alpha, \beta$ be two be two algebraic numbers. Then $h(\alpha\beta ) \leq h(\alpha) + h(\beta)$.\\
b) Let $\alpha_{1},\cdots,\alpha_{n}$ be $n$ algebraic numbers . Then $h(\alpha_{1} +\cdots+ \alpha_{n}) \leq \log(n) + h(\alpha_{1}) +\cdots+ h(\alpha_{n})$.\\
c) $h(\alpha^{k}) = |k|h(\alpha)$, where $k$ is an integer and $|\cdot|$ denotes the usual absolute value.\\
d) Let $\alpha \neq 0$. Then $h(\alpha) = 0$ if and only if $\alpha$ is a root of unity. Moreover, if $\zeta$ is a root of unity then $h(\zeta\alpha) = h(\alpha)$.\\
e) If $\alpha$ is a conjugate of $\beta$ over $\mathbb{Q}$ then $h(\alpha) = h(\beta)$.\\~\\
Proof of these basic properties of height can be found in many textbooks, for example see the book of Bombeiri and Gubler (\cite{heights}, chapter 1).\\~\\

\section{Local fields}
Let $K$ be a valued field and let $w : K \to \mathbb{Z} \cup \{ \infty \}$ denote the surjective valuation. Let $O_{K}$ denote the ring of integers and $k_{K}$ denote the residue field of $K$. In all our computations $k_{K}$ will be finite and $K$ will have characteristic 0.\\
$K/F$ be a Galois extension of valued fields and $w$ be the surjective valuation associated with $K$. We shall assume that valuation on $F$ is non trivial.\\
Let $i \geq -1$ be an integer. Define $$ G_{i}(K/F) = \{ \sigma \in \, \text{Gal} \, (K/F) \mid w (\sigma(a) - a ) \geq i + 1, \, \forall a \in O_K \}$$.\\
Note that $G_{-1}(K/F) = \text{Gal} \, (K/F)$ . $G_{0}(K/F)$ is called the inertia group associated to the extension $K/F$.\\~\\
Let $p$ be a prime, and let $\mathbb{Q}_{p}$ be the field of $p$-adic numbers with absolute value $|\cdot|_{p}$. We shall work with a fixed algebraic closure $\overline{\mathbb{Q}_{p}}$ of $\mathbb{Q}_{p}$ and extend $|\cdot|_p$ to $\overline{\mathbb{Q}_{p}}$.\\~\\
\textbf{Notation :} We shall use these notations till the end of the thesis.\\ 
i) If $f \in \mathbb{N}$, we call $\mathbb{Q}_{p^f}$ to be the unique unramified extension of $\mathbb{Q}_{p}$ of degree $f$ inside $\overline{\mathbb{Q}_{p}}$. The ring of integers of $\mathbb{Q}_{p^f}$ is denoted by $\mathbb{Z}_{p^f}$. \\
ii) Let $\mathbb{Q}_{p}^{\text{nr}}$ be the maximal unramified extension of $\mathbb{Q}_{p}$ inside $\overline{\mathbb{Q}_{p}}$. Further, let $\phi_{p} \in \text{Gal}\,(\mathbb{Q}_{p}^{\text{nr}} | \mathbb{Q}_{p} )$ be the unique lift of the Frobenius automorphism. We write $\phi_{p^f} = \phi^{f}$.\\~\\
Let $\overline {\mathbb{Q}}$ denote the algebraic closure of $\mathbb{Q}$ in $\overline{\mathbb{Q}_{p}}$. If $K$ is a finite extension of $\mathbb{Q}$ inside $\overline{\mathbb{Q}}$ then $|\cdot|_{p}$ restricts to a finite place $v$ of $K$. Then the completion $K_v$ is nothing but the topological closure of $K$ in $\overline{\mathbb{Q}_{p}}$. So, if $K$ is Galois extension of $\mathbb{Q}$ then $\text{Gal}\,(K_{v}|\mathbb{Q}_{p})$ can be thought as a subgroup of $\text{Gal}\,(K|\mathbb{Q})$ by restriction.\\
Now we introduce a notation:\\~\\
\textbf{Notation :} Let $L$ and $K$ be two finite extensions of $\mathbb{Q}_{p}$ such that $ L \subseteq K$. Then we use the notation $ e(K|L)$ to denote the ramification index of $K$ with respect to $L$.\\
We shall use this notation till the end of the thesis.\\~\\
Now we have a lemma :\\~\\
\textbf{Lemma 1.2.1 :} Let $ F \subseteq \overline{\mathbb{Q}_{p}}$ be a finite extension of $\mathbb{Q}_{p}$. Let $ K , L \subseteq \overline{\mathbb{Q}_{p}}$  be finite Galois extensions of $F$ with $K / F$ totally ramified and $L / F$ unramified.\\
i) We have $ K \cap L = F$, and the map $$ \text{Gal} ( KL / F) \to \text{Gal} ( K / F) \times \text{Gal} ( L / F )$$ given by $\sigma \to (\sigma|_{K} , \sigma  |_{L})$ is an isomorphism of groups.\\
ii) The extension $KL / K$ is unramified of degree $ [L : F]$ and the extension $KL / L$ is totally ramified of degree $[K : F]$.\\
iii) Say $ i \geq -1$. If $\sigma \, \in  \text{Gal} ( KL / L) \cap G_{i} ( KL / F )$, then $\sigma |_{K} \in G_{i}( K / F)$. Moreover, the induced map $ \text{Gal} (KL / L) \cap G_{i}( KL / F) \to G_{i} (K / F)$ is an isomorphism of groups.\\~\\  
\textbf{Proof :} \\
\emph{Proof of part (i):} Note that $K\cap L$ is a totally ramified as well as unramified extension of $F$. So this must be a trivial extension ie $K\cap L = F$. Then the statement follows from standard Galois theory (\cite{lang}, chapter VI, theorem 1.14).\\
\emph{Proof of part (ii):} Say $[K : F] = e$ and $[ L : F] =f$. Note that $e(K|F) = e$ and $e(L|F) = f$.\\
Now $e(KL|F) = e(KL|L)e(L|F)= e(KL|L) \leq [KL:L] = [K:F] = e$ where the last equality is by part (i).\\
Again $e(KL|F) = e(KL|K)e(K|L) \geq e$.\\
Combining these $e(KL|F) = e$, $e(KL|L) = [KL:L]$ and $e(KL|K) = 1$. The assertion follows from the last two equality.\\
\emph{Proof of part (iii):} Let $O_K$, $O_L$, $O_{KL}$ be ring of integers of $K, L, KL$ respectively.\\
Use $\mathfrak{p}$ to denote the maximal ideal of $O_{L}$
Let $\pi$ be a generator of the maximal ideal of $O_K$.\\
Before proving the assertion we would like to show the equality \[ O_{KL} = \sum_{l = 0}^{e - 1} \pi^{l}O_{L} \tag{1.2.1} \]
Since $KL/L$ is totally ramified (part (ii)) and $\pi O_{KL} \cap O_{L}= \mathfrak{p}$ we conclude that $O_{KL} = O_{L} + \pi O_{L}$. Iterating this $e$ times one obtains $O_{KL} = O_{L} + \pi O_{L} +...+\pi^{e- 1}O_{L} + \mathfrak{p}O_{KL}$ where in the last step we are using $\pi^{e}O_{KL} = \mathfrak{p}O_{KL}$ which is a easy consequence of part (i) and (ii).\\
But $KL/L$ is a separable extension. So $O_{KL}$ is a finitely generated $O_{L}$ module. Now the desired equality follows from Nakayama lemma.\\
Now let $w : KL \to e^{-1}\mathbb{Z} \cup \{\infty\}$ denote the unique extension of the surjective valuation $F \to \mathbb{Z} \cup \{\infty\}$. Since $e(K|F) = e$ restriction of $w$ to $K$ is also the unique extension of the surjective normalized valuation on $F$.\\
Assume $\sigma \in \text{Gal}(KL|L) \cap G_{i}(KL|F)$. Then $ew(\sigma(a) - a) \geq i + 1$ for all $a \in O_{KL}$. Since $O_{K} \subseteq O_{KL}$ the first part of the assertion follows from the observation above we conclude the first part of the assertion.\\
Injectivity is already a consequence of the result in part $(i)$.\\
Let $\widetilde{\sigma} \in G_{i}(K|F)$. Use part (i) to get a $\sigma \in \text{Gal} (KL|L)$ such that $\sigma|_{K} = \widetilde{\sigma}$. If we can show $\sigma \in G_{i}(KL|F)$ then we shall be done.\\
Assume $ a \in O_{KL}$. Use $(1.2.1)$ to write $a = \sum_{l = 0}^{e -1} \pi^{l}a_l$ with $a_l \in O_L$ for all $ 0 \leq l \leq e - 1$.\\ 
$\sigma$ fixes each of $a_l$. By the assumption on $\widetilde{\sigma}$ we have $ew(\sigma(\pi^l) - \pi^l) \geq i + 1$ for each $l$. Thus by ultrametric triangle inequality $ew(\sigma(a) - a) \geq i + 1$.\\
Holds for all $a \in O_{KL}$. This completes the proof.\hfill $\square$\\~\\ 
The rest of the section is concerned with results in ramification theory. For details one should consult book of Serre (\cite{Serre}, Chapter-4).\\
Let $K$ be a finite extension of $\mathbb{Q}_{p}$. Let $L$ be a finite Galois extension of $K$.\\
Let $k_L$ be the field of residues associated with $L$ . Let $k_L^{*}$ denote the multiplicative group of nonzero elements in the field.\\
For $i \geq -1 $, let $G_{i}(L / K)$ be the $i$-th ramification group as defined before.\\
From standard computations in ramification theory one has $G_{i}$ is normal in $\text{Gal}(L/K)$ for each $i \geq -1$.\\
Now, there are injective homomorphisms \[\theta_{0} : G_{0}/G_{1} \to k_{L}^{*}\] and \[\theta_{i} : G_{i} / G_{i+1} \to k_{L}\] for each $ i \, \geq 1$ where $k_L$ denotes the additive group of the field.\\
With this set up we have the following lemma :\\~\\
\textbf{Lemma 1.2.2 :} Let $s \in G_{0}$ and $\tau \in G_{i}$ for some $ i \geq 1$. Since $G_{i}$ is normal in $G_{0}$ we have  $ s\tau s^{-1} \in G_{i}$ .   Then $$\theta_{i}(s\tau s^{-1}) = \theta_{0}(s)^{i} \theta_{i}(\tau)$$ where we are precomposing by obvious projection maps and the product in the right hand side is taken in the field $k_L$.\\~\\
\textbf{Proof :} \cite{Serre}, chapter 4 , proposition 9.\hfill $\square$ \\~\\

\section{Elliptic curves and supersingular reduction}
Let $K$ be a perfect field with characteristic of $K = p > 0$. Let $E / K$ be an elliptic curve.\\
Then one has the following theorem :\\~\\
\textbf{Theorem 1.3.1 :} Let $$ \phi_{r} : E \to E ^{p^r} \; \text{and} \; \widehat{\phi_{r}} : E ^{p^r} \to E$$ 
be the $p^r$ Frobenius map and its dual for some integer $ r \geq 1$. 
Then the following are equivalent :\\
i)  $E[p^{r}] = 0$ for one (all) $r \geq 1$.\\
ii) $\widehat{\phi_{r}}$ is purely inseparable for one (all) $r \geq 1$.\\
iii) The map $[p] : E \to E $ is purely inseparable and $j(E) \in \mathbb{F}_{p^2}$ where $\mathbb{F}_{p^2}$  is the field with $p^{2}$ elements.\\
iv) The endomorphism ring of $E$ over $\overline{K}$ is an order in a quaternion  algebra , where $\overline{K}$ denotes the algebraic closure of $K$.\\~\\
\textbf{Proof :} See \cite{sil1} theorem 3.1, chapter V.\hfill $\square$ \\~\\
\textbf{Definition 1.3.2 :} Let $ E \,\text{and} \,K$ be as above. If $E$ satisfies one of the  properties in the statement of the theorem above, then $E$ is said to be supersingular.\\~\\
\textbf{Remark 1.3.3 :} Let $K$ be a finite field with $ |K| = q $ where $ q = p^{r}$ where $r$ is an integer $ r \geq 1$. Let $ a_{q} = |E(K)|- q -1 $. Then $E$ is supersingular if and only if $ p \,|\, a_q $.\\~\\
Now let $K$ be a finite extension of $\mathbb{Q}_{p}$. Let \[E : y^{2} + a_{1}xy + a_{3}y = x^{3} + a_{2}x^{2} + a_{4}x + a_{6}\] be an elliptic curve defined over $K$ presented in minimal Weierstrass form.\\
 $E$ is said to have a supersingular reduction at the maximal ideal $\mathfrak{p}$  of $O_{K}$ the ring of integers associated to $K$, $E$ has a good reduction at $\mathfrak{p}$ and the reduced elliptic curve is supersingular over the residue field.\\
If $K$ is a number field, ie a finite extension of $\mathbb{Q}$ and let $E$ be an elliptic curve over $K$. Let $O_K$ be the ring of integers of $K$ and let $\mathfrak{p}$ be a non zero prime ideal of $O_{K}$ .\\
Let $K_{\mathfrak{p}}$ be the completion of $K$ with respect to $\mathfrak{p}$-adic valuation. Clearly it is a finite extension of $\mathbb{Q}_{p}$ for some prime integer $p$. Now one can think $E$ as an elliptic curve over $K_{\mathfrak{p}}$ . If this elliptic curve has supersingular reduction in the sense of the previous paragraph then $E$ is said to have supersingular reduction at the prime ideal $\mathfrak{p}$ of $O_{K}$.\\
Now one has the following theorem :\\~\\
\textbf{Theorem 1.3.3 :} Let $K$ be a number field with at least one real embedding and let $E$ be an elliptic curve over $K$. Then $E$ has infinitely many distinct prime ideals of supersingular reduction.\\~\\
\textbf{Proof :} See \cite{Elkies}.\hfill $\square$ \\~\\

\section{Formal groups}
\textbf{Notation :} R - A commutative ring with unity\\~\\
In this section we recall the basics of formal groups and formal groups over discrete valuation rings. For details one can see \cite{sil1}, chapter 4.\\~\\
\textbf{Definition 1.4.1 :} A (one-parameter, commutative) formal group law $\mathfrak{F}$ defined over $R$ is a power series $F(X,Y) \in R[[X, Y]] $ satisfying :\\
a) $F(X,Y) = X + Y + (\text{terms of degree} \geq 2)$,\\
b) $F(X,F(Y,Z)) = F(F(X,Y),Z)$,\\
c) $F(X,Y) = F(Y,X)$,\\
d) $F(X,0) = X$ and $ F(0,Y) = Y$. \\~\\
\textbf{Remark 1.4.2 :} i) Let $F(X,Y) \in R[[X,Y]]$ be such that it satisfies the conditions in a) and b) above, then there is a unique power series $i(T) \in R[[T]]$ such that the constant term of $i(T)$ is zero and $F(T,i(T)) = 0$.  \\
ii) It can be shown $(a)$ and $(b)$ implies $(d)$. Further $(a)$ and $(b)$ implies $(c)$ provided $R$ has no element which is both torsion and nilpotent (see \cite[Ch.4]{sil1}, \cite[Ch.1]{haz2}). Note that $(d)$ implies $F(X,Y) \in X + Y + XYR[[X,Y]]$.\\ 
iii) If the formal group $\mathfrak{F}$ is defined by the power series $F$, then we write it as a ordered pair $(\mathfrak{F},F)$.\\~\\
\textbf{Definition 1.4.3 :} Let $(\mathfrak{F}, F)$ and $(\mathfrak{G}, G)$ be formal groups defined over $R$. A homomorphism from $\mathfrak{F}$ to $\mathfrak{G}$ defined over $R$ is a power series (with no constant term) $f(T) \in R[[T]]$  satisfying  \[f(F(X,Y)) = G(f(X),f(Y)).\] 
$f$ is said to be an isomorphism if there is another power series $g$ in $R[[T]]$ such that $g$ is a homomorhism from $\mathfrak{G}$ to $\mathfrak{F}$ and $f(g(T)) = g(f(T)) = T$.\\~\\
\textbf{Example 1.4.4 :} Let $(\mathfrak{F}, F)$ be a formal group defined over some commutative ring $R$. For each $m \in \mathbb{Z}$ we construct homomorphisms  $[m]: \mathfrak{F} \to \mathfrak{F} $ as follows :\\
Define $[0](T) = 0 $.  For each non negative integer $m$  define $[m+1](T) = F ([m]T, T)$ . If $m$ is a non positive integer define $[m - 1](T) = F([m]T ,i(T))$ where $i(T)$ is the unique element in $R[[T]]$ satisfying $F(T ,i(T)) = 0$. \\
It is not very difficult to verify that $[m]$ is indeed an homomorphism $\mathfrak{F} \to \mathfrak{F}$.\\~\\
Now one has the following theorem :\\~\\
\textbf{Theorem 1.4.5 :} Let $\mathfrak{F}$ be a formal group over $R$ and let $ m \in \mathbb{Z}$. Then we have:\\
a) $[m](T) = mT + \text {higher order terms}$.\\
b) If $m$ is an unit in $R$ then $[m]$ is an isomorphism of the formal group $\mathfrak{F}$.\\~\\
\textbf{Proof :} See \cite{sil1} chapter 4 , proposition 2.3. \hfill $\square$ \\~\\
\subsection {Groups associated to formal groups}
\textbf {Notation :} $( R, \mathfrak{M} )$ -  A local ring whose unique maximal ideal is $\mathfrak{M}$ and which is complete with respect to $\mathfrak{M}$-adic topology,\\
$(\mathfrak{F}, F)$ - A formal group law defined over $R$,\\
$ k_{R}$ - The field $ R / \mathfrak{M} $,\\
$ i(T)$ - The unique element in $R[[T]]$ such that $F(T,i(T)) = 0$.\\~\\
First note that since $R$ is complete, $F(X, Y)$ converges when evaluated at a point of $\mathfrak{M} \times \mathfrak{M}$. Similarly $i(T)$ converges when evaluated at a point of $\mathfrak{M}$. This allows us to make the following definition :\\~\\
\textbf{Definition 1.4.6 :} The group associated to $\mathfrak{F}$, denoted $\mathfrak{F}(\mathfrak{M})$, is the set $\mathfrak{M}$ with the group operations $$ x \oplus_{\mathfrak{F}} y = F(x, y) \hspace{2cm} \text{(Addition)} \hspace{1cm}\text{for}\, x, y \in \mathfrak{M},$$  $$ \ominus_{\mathfrak{F}} x = i (x) \hspace{2cm} \text{(Inverse)} \hspace{1cm} \text{for} \, x \in \mathfrak{M}.$$
It is not hard to check that it really defines a group law.\\~\\
\textbf{Remark 1.4.7 :} i) Let $n \geq 1$. Note that $ F (\mathfrak{M}^{n} \times \mathfrak{M}^{n}) \subseteq \mathfrak{M}^{n}$ and $i(\mathfrak{M}^{n}) \subseteq \mathfrak{M}^{n}$. Thus the subset $\mathfrak{M}^{n}$ of $\mathfrak{M}$ is a subgroup of $\mathfrak{F}(\mathfrak{M})$ which is denoted as $\mathfrak{F}(\mathfrak{M}^{n})$.\\
ii) We shall use the notation $ x \ominus_{\mathfrak{F}} y $ to denote $ x \oplus_{\mathfrak{F}} (\ominus_{\mathfrak{F}} y)$.\\~\\
Now we have :\\~\\
\textbf{Proposition 1.4.8 :} a) For each $ n \geq 1$ , the map $$ \mathfrak{F}(\mathfrak{M}^{n}) / \mathfrak{F}(\mathfrak{M}^{n+1}) \to \mathfrak{M}^{n} / \mathfrak{M}^{n+1}$$ induced by the identity map on the underlying sets is actually an isomorphism of groups.\\
b) Let $p$ be the characteristic of $k_{R}$. Then every torsion element of $\mathfrak{F}$ has order a power of $p$ (if $p = 0$ then there is no nontrivial torsion point).\\~\\
\textbf{Proof :} \cite{sil1}, chapter 4, proposition 3.2. \hfill $\square$ \\~\\
For the rest of this subsection we restrict to the case where $R$ is a discrete valuation ring and $k_{R}$ has characteristic $p > 0$.\\
Let $K$ be the fraction field of $R$ and we assume that characteristic of $K$ is $0$. Let $v : K - \{0\} \to \mathbb{Z}$ denote the surjective valuation.\\
Assume \[[p](T) = pT + a_{2}T^{2} + a_{3}T^{3} + \cdots\]
where $a_{i} \in R$ for each $i \geq 2$.\\
We define \[\alpha = - \min_{i \geq 2}\, \frac{1} {i - 1} \,v \left(\frac {a_{i}}{p}\right) \tag{1.4.1}.\]
Then one has the following proposition :\\~\\
\textbf{Proposition 1.4.9 :} One can extend the $\mathbb{Z}$ module structure on $\mathfrak{F}(\mathfrak{M})$ to a $\mathbb{Z}_{p}$ module structure,  such that $\mathfrak{F}(\mathfrak{M}^{n})$ is submodule for each $n \geq 1$.\\
Further if $k_{R}$ is finite then $\mathfrak{F}(\mathfrak{M}^{n})$ is free of finite rank for $n > \alpha $.\\~\\
\textbf{Proof :} See \cite{oshikawa}, remark 2 and corollary 1 after theorem 1 in section 1. \hfill $\square$ \\~\\
\subsection{Formal group law of an elliptic curve}
Let $F$ be a field and let $E$ be an elliptic curve defined over $F$ which is presented in the Weierstrass form \[ y^{2} + a_{1}xy + a_{3} y = x^{3} + a_{2} x^{2} + a_{4}x + a_{6} . \tag{1.4.2}\]
Now if we make the change of variable $z = -\frac {x} {y}$ and $w = -\frac {1} {y}$ then $O$ on $E$ is now on the point $(z, w) = (0, 0)$, and $z$ is a local uniformizer at $O$ since it has a zero of order $1$ at $O$. \\
With the new coordinates the Weierstrass equation looks like \[ w = z^{3} + a_{1}zw + a_{2}z^{2}w + a_{3} w^{2} + a_{4}zw^{2} + a_{6} w^{3} (= f(z, w)).\]
Now one has the following proposition:\\~\\
\textbf{Proposition 1.4.10 :} Let $R$ be any ring containing $\mathbb{Z} [a_{1},\cdots,a_{6}]$. Then the following holds:\\
i) There is an unique formal power series $w(Z) = Z^{3} (1 + A_{1} Z + A_{2} Z +\cdots) \in \, \mathbb{Z}[a_{1},\cdots,a_{6}] [[Z]] \subseteq R[[Z]]$ such that it is a solution to the  equation $ w(Z) = f (Z, w(Z))$ in $R[[Z]]$.\\
ii) There are formal Laurent series $x(Z), y(Z) $ with coefficients in $\mathbb{Z}[a_{1},\cdots,a_{6}]$ defined by the formula \[x(Z) = \frac {Z} {w(Z)} \; \text{and} \; y(Z) = \frac {- 1} {w(Z)}\] such that they are formal solution of Weierstrass equation  $(1.4.2)$.\\
iii) Now assume that $R$ is a complete discrete valuation ring with maximal ideal $\mathfrak{M}$ and let $K$ denote the fraction field of $R$. Clearly $ F \subseteq K$.  Then $ (x(Z), y(Z))$ converges on $(\mathfrak{M} - \{0\}) \times ( \mathfrak{M} - \{0\})$ and it describes a point on $E(K)$.\\
Further the map \[\mathfrak {M} \to E(K) \tag{1.4.3}\] obtained this way ($0$ is mapped to the point at infinity) is injective.\\
iv) There is a unique power series $F(Z_{1}, Z_{2}) \in R\,[Z_{1}, Z_{2}]$ with coefficients in $\mathbb{Z}[a_{1},\cdots,a_{6}]$ such that it defines a formal group law $\mathfrak{F}$ and the map in $(1.4.3)$ induces a homomorphism from $\mathfrak{F}(\mathfrak{M})$ to $E(K)$.\\  
v) Let $k$ denote the field of residues associated with $R$ and we assume that $k$ is finite.
Consider the reduction modulo $\mathfrak{M}$ map $$\mathbb{P}^{2}(K) \to \mathbb{P}^{2}(k)$$ and let $\widetilde{O}$ denote the image of $O$ under this map.\\
Let $E_{1}(K) = \{P \in E(K) \; | \; \widetilde{P} = \widetilde{O}\}$ . Then the image of the map in $(1.4.3)$ is exactly $E_{1}(K)$.\\~\\
\textbf{Proof :} For proof of (i), (ii), (iii), (iv) see chapter 4 of \cite{sil1}.\\
For proof of (v) see chapter 7 in the same reference. \hfill $\square$

\chapter{Galois properties of p-torsion points}
\section{Introduction}
\textbf{Notation :} Let $E$ be an elliptic curve defined over a field $K$ and let $N \in \mathbb{N}$. Further let $\overline{K}$ be a fixed algebraic closure of $K$ and $L$ be a subfield of $\overline{K}$. Fix an Weierstrass model for $E/K$. Then $E[N]$ denotes the group of $N$-torsion points of $E$ over $\overline{K}$ and $L(N)$ denotes the field generated by $x$ and $y$ coordinates of points in $E[N]$ over $L$. This field is independent of choice of Weierstrass model of $E$ over $K$. \\~\\ 
In this chapter we work with the set up:\\
Let $f$ be a fixed positive integer, let $p$ be a prime $\geq 5$ and let $w=p^{2f}$. Let \[ E: y^{2}= x^{3}+ Ax + B \tag{2.1.1}\] be a fixed elliptic curve defined over $\mathbb{Z}_{w}$ (notation as in (1.2.1)) presented in minimal Weierstrass form such that its reduction modulo $p$,\[\widetilde{E}: y^{2}=x^{3}+ \widetilde{A}x + \widetilde{B} \tag{2.1.2}\] is a supersingular elliptic curve (in particular $\widetilde{E}$ is nonsingular). We shall also assume that the $j-$invariant of the reduced elliptic curve is not among $0$ or $1728$. \\
For notational convenience we put $ q = p^{2}$ in what follows.\\~\\
The main goal of this chapter is to prove the following theorem :\\~\\
\textbf{Theorem 2.1.1 :} Let $ n \in \mathbb{N}$.\\
i) The extension $\mathbb{Q}_{w}(p^{n}) / \mathbb{Q}_{w}$ is totally ramified and abelian of degree $ q^{n - 1}(q - 1)$.\\
Moreover, \[\text{Gal}(\mathbb{Q}_{w}(p^{n})|\mathbb{Q}_{w})\,\cong \,\mathbb{Z} /(q - 1)\mathbb{Z} \times (\mathbb{Z} / p^{n - 1}\mathbb{Z} )^{2} \tag{2.1.3}\] and \[\text{Gal}(\mathbb{Q}_{w}(p^{n})|\mathbb{Q}_{w}(p^{n - 1}))\,\cong \, (\mathbb{Z} / p \mathbb{Z})^{2} \tag {2.1.4} \] for $ n \, \geq 2$.\\
ii) Let $k$ and $i$ be integers with $ 1 \leq k \leq n $ and $ q^{k - 1} \leq i \leq q^{k} - 1 $. Then \[ G_{i} (\mathbb{Q}_{w} (p^n)|\mathbb{Q}_{w}) = \text{Gal} ( \mathbb{Q}_{w}(p^{n}) | \mathbb{Q}_{w}(p^{k}) )\tag{2.1.5} \]
iii) Let $ M \in \mathbb{N} $ be coprime to $p$. The image of the representation \[\text{Gal} (\mathbb{Q}_{w}(p^n)|\mathbb{Q}_{w}) \to \text{Aut} E[p^{n}] \tag{2.1.6} \] contains multiplication by $\pm M $ and it acts transitively on torsion points of order $p^{n}$, where $\text{Aut} E[p^{n}]$ is the group of automorphisms of $E[p^{n}]$ as an abelian group.\\~\\
\textbf{Remark 2.1.2 :} The results in chapter 2 are analogues of the results in section 3 of \cite{hab}. Habbeger proves  these results for the case $f = 1$ using Lubin-Tate theory and twisting of elliptic curves. Our technique puts more emphasis on theory of formal groups and ramification theory. This is the only place where the treatment in our work differs significantly from that of \cite{hab}.

\section{Setting it up}
\subsection{An important exact sequence}
Let $K$ be any finite extension of $\mathbb{Q}_{w}$ and let $k_{K}$ denote its residue field.\\
Now if one considers $E$ as an elliptic curve over $K$, then $E/K$ has a supersingular reduction because $E$ is supersingular over $\mathbb{Q}_{w}$ implies $\widetilde{E}(\overline{\mathbb{F}_{p}})$ has no non-trivial $p-$torsion.\\
We shall abuse notation and use $\widetilde{E}$ to denote the image of this elliptic curve under the reduction map.\\
Let $\widetilde{O}$ denote the origin of the reduced elliptic curve.\\
Put $E_{1}(K) = \{P \in E(K) \mid \widetilde{P} = \widetilde{O}\}.$\\
One can define a complex of abelian groups \[0 \to E_{1}(K) \to E(K) \to \widetilde{E}(k_K) \to 0 \tag{2.2.1}\] where the first map is inclusion and the second map is reduction.\\
By Proposition-2.1, Chapter-VII in (\cite{sil1}) one concludes that $(2.2.1)$ is exact.\\
Since $E$ is supersingular, all its $p^{n}-$torsions ($n \geq 1$) come from $E_{1}(K)$. \\ ~\\
Let $O_{K}$ be the ring of integers of $K$. We know that $O_{K}$ is a complete discrete valuation ring. Let $\mathfrak{p}_{K}$ denote the unique maximal ideal of $O_{K}$.\\
Let $(\mathfrak{F}(\mathfrak{p}_{K}), F)$ be the formal group law associated to $E$.\\
By proposition 1.4.10, we have an isomorphism \[ \mathfrak{F}(\mathfrak{p}_{K}) \to E_{1}(K) \tag{2.2.2} \]
which is defined by $z \to (x(z),y(z))$ ($0$  goes to the point at infinity, $O$).\\~\\
\textbf{Remark 2.2.1 :} i) Since the formal group law is defined over $\mathbb{Q}_{w}$, the isomorphism in $(2.2.2)$ commutes with the action of absolute Galois group of $\mathbb{Q}_{w}$.\\
ii) If $(x,y) \in E_{1}(K)$ and $z$ is the preimage of $(x,y)$ under $(2.2.2)$ then we call $z$ to be the local parameter of $(x,y)$.\\
iii) Let $n$ be a positive integer. It follows from the argument above that :\\
$(x,y) \in E(K)$ has order $p^{n}$ if and only if there is a $z \in \mathfrak{p}_{K} - \{0\}$ such that $x = x(z)$, $y = y(z)$ and $z$ has order $p^{n}$ in $\mathfrak{F}(\mathfrak{p}_{K})$.\\~\\
\textbf{Lemma 2.2.2 :} Let $z \in \mathfrak{p}_{K} - \{0\}$. Then $\mathbb{Q}_{w}(z) = \mathbb{Q}_{w} (x(z), y(z))$.\\~\\
\textbf{Proof :} Note that $z = \frac {- x(z)} {y(z)}$.\\
Hence, $\mathbb{Q}_{w}(z) \subseteq  \mathbb{Q}_{w} (x(z),y(z))$.\\
The other inclusion follows from the remark 2.2.1 (i).\\
This proves the lemma. \hfill $\square$  

\subsection{Multiplication by \emph{m} map}
Let $K$ be a field such that its characteristic is not 2 or 3 and let $E$ be an elliptic curve defined over $K$ presented in Weierstrass form $E: y^{2} = x^{3} + Ax + B$. Now one has,\\~\\
\textbf{Lemma 2.2.3 :} Let $m$ be a odd positive integer.\\
Then $[m]: E \to E$ is given by a rational function $(\frac{\phi_{m}(x)}{\psi_{m}(x)^{2}} , \frac{\omega_{m}(x)}{\psi_{m}(x)^{3}})$ where $\phi_{m}(x), \psi_{m}(x), y^{-1}\omega_{m}(x)$ are polynomials in $\mathbb{Z}[A, B, x]$.\\
Further degree of $\phi_{m}$ and $\psi_{m}$ are polynomials of degree $m^{2}$ and $\frac{m^{2} - 1}{2}$ (in $x$) respectively.\\~\\
\textbf{Proof :} A sketch of proof can be found in Exercise - 3.7 of (\cite{sil1}).\hfill $\square$ \\~\\
\textbf{Remark 2.2.4 :} Let $K, E, m$ be as above. Let $(x_{0},y_{0})$ be a $m$-th division point of $E$ defined over $\overline{K}$. Then clearly one has $\phi_{m}(x_{0}) = 0$. Thus $[K(x_{0}): K] \leq \frac{m^{2} - 1}{2}$. Now $[K (x_{0})(y_{0}) : K (x_{0})] \leq 2$. Thus $[K(x_{0}, y_{0}) : K] \leq m^{2} - 1$. Hence the finiteness assumption in previous section is not restrictive.\\   

\section{Computations}
Now we are in the setup considered in the introduction.
\subsection{Computation of $\alpha$}
Let $K$ be a finite extension of $\mathbb{Q}_{w}$.\\
Let $O_{K}$ denote the ring of integers of $K$ and let $\nu_{K}$ be the usual valuation on $O_{K}$. Since $E$ is defined over $\mathbb{Z}_{w}$ it is also defined over $O_{K}$ and we can apply the theory in previous subsection.\\
Consider the curve $\widetilde{E}$. It is defined over $\mathbb{F}_{w}$. Let $\widetilde{\phi_{w}} \in \text{Gal}(\overline{\mathbb{F}_{p}}/ \mathbb{F}_{w})$ be the map $ x \to x^{w}$. $\widetilde{\phi_{w}}$ defines an endomorphism of $\widetilde{E}$ which is purely inseparable and of degree $w$. Now $[p^{f}]$ map of $\widetilde{E}$ has degree $w$ and this map is purely inseparable since $\widetilde{E}$ is supersingular. So one can find an automorphism $u$ of $\widetilde{E}$ such that $[p^{f}] = u \circ \widetilde{\phi_{w}}$ (see proposition 13.5.4 in \cite{Huse}).\\
Since $j-$invariant of the reduced curve is not $0$ or $1728$, consulting the table of isomorphism in Silverman's book (section-10, chapter-3 of \cite{sil1}) we conclude that $[p^{f}]= \pm \widetilde{\phi_{w}}$ on $\widetilde{E}$. \\
From the previous observation one concludes that thought as a morphism of the formal group associated to $E$, $[p^{f}](X) \equiv \pm X^{w} \mod{p\mathbb{Z}_{w}}$ .\\
Now we have the following claim : \\~\\
\textbf{Lemma 2.3.1 :} Let $[p](X) = pX + d_{2}X^{2} + d_{3}X^{3} + \cdots$ be the multiplication by $p$ map. Then $\min \{i \geq 2 \mid p \nmid d_{i}\} = p^{2}(= q)$.\\~\\
\textbf{Proof :} We know that $[p^{f}](X) \equiv \pm X^{w} \mod p \mathbb{Z}_{w}$. \\
This shows that at least one of the coefficients of $[p](X)$ is not divisible by $p$.\\
Now if the minimum is $m$ then a simple calculation shows that the coefficient of $X^{m^{f}}$ in $[p^{f}](X)$ not divisible by $p$.\\ Thus one must have $m = p^{2}$. \hfill $\square$ \\~\\ 
From $(1.4.1)$ we have, \[\alpha_{K} = - \min_{i \geq 2}\, \frac{1} {i - 1} \,v_{K} \left(\frac {a_{i}}{p}\right).\]
Since all the coefficients are in $\mathbb{Z}_{w}$, an unramified extension of $\mathbb{Z}_{p}$, the minimum will be attained at $ i = m$ where $m$ is as in lemma 2.3.1.\\
Thus we have,\[\alpha_{K} =  \frac {v_{K}(p)} {q - 1} \tag{2.3.1}.\]

\subsection{Computation of Galois group}
We start off by noting a simple fact :\\~\\
\textbf{Lemma 2.3.2 :} Let $n$ be a positive integer. Then $\mathbb{Q}_{w}(p^{n})$ is a Galois extension of $\mathbb{Q}_{w}$.\\~\\
\textbf{Proof :} Follows from the fact that group law of $E$ is defined over $\mathbb{Q}_{w}$. \hfill $\square$\\~\\
Next we record a small observation: \\~\\
\textbf{Remark 2.3.3 :} i) If $(x,y)$ and $(x,y)$ are two distinct points on $E$ then by $(2.1.1)$ we have $y = -y'$.
Thus with respect to the group law on $E$ these two points are inverse of each other and hence have same order.\\
ii) Let $1 \leq i \leq k$ and assume $\sigma \in \text{Gal}(\overline{\mathbb{Q}}_{p} |\mathbb{Q}_{w})$. $E$ is defined over $\mathbb{Q}_{w}$. So, $(\sigma (x_i),\sigma (y_i))$ has same order as $(x_i, y_i)$. Hence all conjugates of $x_i$ over $\mathbb{Q}_{w}$  are $x$-coordinates of points which have order $p^{k + 1 - i}$.\\~\\
\textbf{Lemma 2.3.4 :}
Let $n$ be a positive integer.\\
Consider the statement :\\
$\textbf{P(n)}$ : i) Let $(x_{0}, y_{0})$ be a point on $E$ with exact order $p^{n}$ defined over $\overline{\mathbb{Q}_{p}}$.\\ Then $[\mathbb{Q}_{w}(x_{0}) : \mathbb{Q}_{w}] = \frac{q^{n-1}(q-1)}{2}$ and $[\mathbb{Q}_{w}(x_{0}, y_{0}) : \mathbb{Q}_{w}] = q^{n-1}(q-1)$.\\
ii) The extension $\mathbb{Q}_{w}(x_{0}, y_{0}) / \mathbb{Q}_{w}$ is totally ramified.\\
Then $P(n)$ is true for all $n \geq 1$.\\~\\
\textbf{Proof :} We shall show that $P(n)$ holds for each $n \geq 1$ by induction.\\
The main idea of the proof is to obtain an upper bound for $[\mathbb{Q}_{w}(x_0, y_0) : \mathbb{Q}_{w}]$ and a lower bound for $e(\mathbb{Q}_{w}(x_{0}, y_{0})\,|\,\mathbb{Q}_{w})$ and to compare these two bounds.\\~\\
Let $(x_{0},y_{0})$ be a point of exact order $p$.\\
By remark 2.2.3 we already know \[[\mathbb{Q}_{w}(x_{0}) : \mathbb{Q}_{w}] \leq \frac{q-1}{2}\,\text{and}\, [\mathbb{Q}_{w}(x_{0}, y_{0}):\mathbb{Q}_{w}] \leq (q-1). \tag{2.3.2}\]
Put $K = \mathbb{Q}_{w}(x_{0}, y_{0})$ and let $\mathfrak{p}_{K}$ be the maximal ideal of $O_{K}$, the ring of integers of $K$.\\
Since $(x_{0}, y_{0})$ is a $p$-torsion, we conclude that $(x_{0}, y_{0}) \in E_{1}(K)$. Now $z_{0} \in \mathfrak{p}_{K}$ be the local parameter of $(x_{0}, y_{0})$. By lemma $2.2.2$ we have $K = \mathbb{Q}_{w}(z_{0})$.\\
Here $\mathfrak{F}(\mathfrak{p}_{K})$ has a nontrivial torsion point namely $z_0$. Hence $[p]$ is not injective on $\mathfrak{p}_{K}$.\\
Thus by proposition 1.4.9 we have $\alpha_{K} \geq 1$.\\
From $(2.3.1)$ we have $\nu_{K}(p) \geq (q-1)$.\\
So \[e(K |\mathbb{Q}_{w}) \geq (q-1). \tag{2.3.3} \] 
Since $[K:\mathbb{Q}_{w}] \geq e(K|\mathbb{Q}_{w})$, $(2.3.2)$ and $(2.3.3)$ imply $[K:\mathbb{Q}_{w}] = (q - 1)$.\\
Now $[K :\mathbb{Q}_{w}(x_{0})] \leq 2 $. Hence $[\mathbb{Q}_{w}(x_0) : \mathbb{Q}_{w}] \geq \frac{q-1}{2}$. Using $(2.3.2)$ we have $[\mathbb{Q}_{w}(x_{0}) : \mathbb{Q}] = \frac{(q - 1)} {2}$.\\
This proves $P(1)(i)$.\\
$P(1)(ii)$ follows from $(2.3.3)$ and the fact that $ [K : \mathbb{Q}_{w}] = (q - 1)$.\\
This finishes the proof in base case.\\~\\
Now we assume that $P(1),\cdots,P(k)$ hold for some positive integer $k$. We shall show that $P(k+1)$ also holds.\\~\\
Let $(x_{0}, y_{0})$ be a point of order exactly $p^{k+1}$. Put $K = \mathbb{Q}_{w}(x_{0},y_{0})$ and let $z_{0}$ be the local parameter associated to $(x_{0}, y_{0})$. We know $K = \mathbb{Q}_{w}(z_0)$.\\
Define $z_{i} = [p^{i}](z_{0})$ for each $1 \leq i \leq k$. Assume that the image of the point $z_{i}$ on the elliptic curve is $(x_{i}, y_{i})$. The point $z_{i}$ has exact order of $p^{(k+1)-i}$.\\
Let $f_{i}(X)$ denote the minimal polynomial of $x_{i}$ over $\mathbb{Q}_{w}$.\\
By induction hypothesis \[\text{degree}(f_{i}(X)) = \frac{q^{k-i}(q-1)}{2}.\tag {2.3.4}\]
Consider the polynomial $\psi_{p^{k+1}}(X)$ as in lemma 2.2.3.\\ 
Clearly $\psi_{p^{k+1}}(x_i) = 0$ for all $ 0 \leq i \leq k$.\\
Hence $ f_{i}(X) \, | \, \psi_{p^{k+1}}$  for all $ 0 \leq i \leq k$, where $f_{i}(X)$ is the minimal polynomial of $x_{i}$ over $\mathbb{Q}_{w}$.\\
By remark 2.3.3 it follows that these polynomials have no common roots in $\overline{\mathbb{Q}_{p}}$.\\
So \[\prod_{i=0}^{k} f_{i}(X)\, | \,\psi_{p^{k+1}}(X) .\tag{2.3.5}\]  \\
From lemma 2.2.3, $(2.3.4)$ and $(2.3.5)$ one has \[\text{degree}(f_{0}(X)) \leq \frac{q^{k+1}-1} {2} - \sum_{i =1}^{k} \frac{q^{k - i}(q - 1)} {2} \\ \leq \frac{q^{k}(q - 1)} {2} .\tag{2.3.6} \]  
Now $(2.3.6)$ translates into $ [\mathbb{Q}_{w} (x_0) : \mathbb{Q}_{w}] \leq \frac {q^{k}(q - 1)} {2}$. Using the fact that \[[\mathbb{Q}_{w}(x_0, y_0) : \mathbb{Q}_{w} (x_0)] \leq 2\] we conclude $[K : \mathbb{Q}_{w}] \leq q^{k} (q - 1) $.\\
Put $\mathbb{Q}_{w}(z_{1}) = K_{1}$. Let $O_{K_1}$ denote the ring of integers and let $\mathfrak{p}_{K_{1}}$ denote the maximal ideal in $O_{K_1}$. Note that, $z_{1} \in \mathfrak{p}_{K_1}$.\\
Further by induction hypothesis $[K_{1} : \mathbb{Q}_{w}] = q^{k - 1}(q - 1)$ and $K_{1}$ is totally ramified over $\mathbb{Q}_{w}$. It follows that \[ z_{1} \in \mathfrak{p}_{K_{1}} -  \mathfrak{p}_{K_{1}}^{2}. \tag{2.3.7}\]\\  
Note that $ K_{1} \subseteq K$. Hence \[e(K | \mathbb{Q}_{w} ) \geq e ( K_{1} | \mathbb{Q}_{w}) \geq q^{k - 1}(q - 1) \geq (q - 1).\] So $\nu_{K}(p) \geq (q - 1)$.\\
Now using lemma 2.3.1 one deduces $[p](z_{0}) \in \mathfrak{p}_{K}^{p^{2}}$. But $[p](z_{0}) = z_{1} $.\\
By observation $(2.3.7)$ and the remark above we conclude that $e(K|K_{1}) \geq p^{2}$.\\
Using induction hypothesis $e(K_{1}|\mathbb{Q}_{w})= q^{k-1}(q-1)$. Thus \[e(K|\mathbb{Q}_{w}) \geq q^{k}(q-1).\]
This with previous bound, $[K : \mathbb{Q}_{w}] \leq q^{k}(q - 1) $,  proves $P(k+1)(i)$ and $P(k+1)(ii)$.\\~\\
Hence by principle of mathematical induction we are done. \hfill $\square$ \\~\\
\textbf{Remark 2.3.5 :} Let $z_{0}$ be a torsion point of $\mathfrak{F}(p\mathbb{Z}_{w})$ defined on $\overline{\mathbb{Q}}_{p}$. It's order will be a power of $p$. \\
From the proof of the lemma it follows that $\mathbb{Q}_{w}(z_{0})$ is a totally ramified extension and $z_{0}$ is a generator of the unique maximal ideal of its ring of integers.\\ ~\\ 
\textbf{Lemma 2.3.6 :} Let $(x_{0}, y_{0})$ be a nontrivial $p$-torsion point of $E$.\\
Then \[\mathbb{Q}_{w}(x_{0}, y_{0}) = \mathbb{Q}_{w}(p).\] 
Further, \[\text{Gal}\,(\mathbb{Q}_{w}(p)|\mathbb{Q}_{w})\cong\mathbb{Z}/(q - 1)\mathbb{Z}.\]
\textbf {Proof :} Let $z_0$ be the local parameter associated to the point $(x_{0}, y_{0})$. Put $\mathbb{Q}_{w}(x_0, y_0) = \mathbb{Q}_{w}(z_{0}) = K$.
By lemma 2.3.4 $K$ is a totally ramified extension of $\mathbb{Q}_{w}$ whose degree is $(q - 1)$.\\
Thus $z_{0}$ has $(q-1)$ many conjugates over $\mathbb{Q}_{w}$ all of which has order exactly $p$ in $\mathfrak{F}(\mathfrak{p}_{K})$.
But from standard results in theory of elliptic curves we know that there are exactly $(q - 1)$ many points of order $p$. So all of them must be conjugates of $z_0$ over $\mathbb{Q}_{w}$ (here we are using the isomorphism in $(2.2.2)$).\\
Since $\mathbb{Q}_{w}(p)$ is Galois over $\mathbb{Q}_{w}$ it follows that this must be the Galois closure of $K$ over $\mathbb{Q}_{w}$.\\
$K$ is a tamely ramified extension of $\mathbb{Q}_{w}$ ($\gcd (p, q - 1) = 1$). Applying a standard lemma in algebraic number theory (proposition-12 in chapter-2; \cite{lang2}) we know that there is a generator $\Pi$ of $\mathfrak{p}_{K}$ and a generator $\pi$ of $p\mathbb{Z}_{w}$ such that $\Pi^{q-1} = \pi$. Thus $K$ is Kummer extension of $\mathbb{Q}_{w}$.\\
Now $X^{q-1} - 1$ has $q - 1$ many distinct solutions in $\mathbb{F}_{w}$ (since $(q -1)\mid w-1$) where $\mathbb{F}_{w}$ is the field with $w$ elements.
All these solutions lift to a solution in $\mathbb{Q}_{w}$ because $\gcd (q - 1,p) = 1$. Hence all the $(q - 1)$ roots of unities are in $\mathbb{Q}_{w}$. \\
So $K$ is Galois over $\mathbb{Q}_{w}$ and the Galois group is $\mathbb{Z}/(q-1)\mathbb{Z} $ .\\
This finishes proof of the lemma.\hfill $\square$ \\~\\ 
\textbf{Notation :} Let $n \geq 1$. We shall use the notation $\mathfrak{F}[p^{n}]$ to denote the group of $p^{n}$ torsion points of $\mathfrak{F}$. This is isomorphic to $E[p^{n}]$ by the isomorphism in $(2.2.2)$.\\~\\
\textbf{Lemma 2.3.7 :} Let $(x_0, y_0)$ be a point of order $p^{n}$, where $n$ is a positive integer. Then $\mathbb{Q}_{w}(x_0, y_0) = \mathbb{Q}_{w}(p^{n})$.\\~\\
Lemma 2.3.7 is quite important for rest of the developments of this chapter. We shall assume it for time being and differ discussion on it till appendices (appendix B).\\~\\
\textbf{Lemma 2.3.8 :} Let $G = \text{Gal}\,(\mathbb{Q}_{w} (p^{n+1}) / \mathbb{Q}_{w} (p^{n}))$ where $n$ is a positive integer. Then $G \cong \mathfrak{F}[p]$.\\~\\
\textbf{Proof :} Let $P$ be a point such that it has order $p^{n+1}$ in $\mathfrak{F}$. Then $[p]_{\mathfrak{F}}(P)$ has order $p^{n}$ and thus it lies in the ground field.\\
Let $\sigma \in G$. Clearly $\sigma ([p]_{\mathfrak{F}}(P)) = [p]_{\mathfrak{F}}(P)$. The group law is defined over the ground field. So $\sigma (P) \ominus_{\mathfrak{F}} P \in E[p]$.\\
This is true for all $\sigma \in G$.\\ 
So we have a well defined map $\Delta_{P} : G \to \mathfrak{F}[p]$ given by \[\Delta_{P}(\sigma) := \sigma (P) \ominus_{\mathfrak{F}} P.\] 
This map is an injection since $P$ generates $\mathbb{Q}_{w} (p^{n+1})$ over $\mathbb{Q}_{w}$.\\
From lemma 2.3.4 and lemma 2.3.7 we have $[\mathbb{Q}_{w}(p^{n+1}) : \mathbb{Q}_{w}(p^{n})] = q$. So $|G| = q$. But $|E[p]| = |\mathfrak{F}[p]| = q $.\\
Comparing the cardinality of the sets, one has $\Delta_{P}$ is onto.\\
Now let $\sigma, \tau \in G$. Then 
\[
\begin{split}
\Delta_{P}(\sigma \tau) & = \sigma \tau (P)\,\ominus_{\mathfrak{F}}\,P \\ 
                        & = \sigma(\tau (P)\,\ominus_{\mathfrak{F}}\,P)\,\oplus_{\mathfrak{F}}\,(\sigma (P)\, \ominus_{\mathfrak{F}}\,P)\\ 
                        & = \sigma (\Delta_{P}(\tau))\,\oplus_{\mathfrak{F}} \, \Delta_{P}(\sigma).
\end{split}\]
Since $n \geq 1$, we conclude that $\sigma$ fixes $\mathfrak{F}[p]$. Thus
\[\Delta_{P}(\sigma \tau) = \Delta_{P}(\sigma) \, \oplus_{\mathfrak{F}} \, \Delta_{P}(\tau).\]
True for all $\sigma, \tau \in G$. Hence $\Delta_{P}$ is a homomorphism of groups.\\
This proves the lemma. \hfill $\square$\\~\\
\textbf{Notation :} For the current section we fix the following notations :\\ 
i) $G_{n} = \text{Gal} (\mathbb{Q}_{w}(p^{n})|\mathbb{Q}_{w})$ for all $n \in \mathbb{N}$.\\
ii) $G_{n, i} = G_{i}(\mathbb{Q}_{w}(p^{n})|\mathbb{Q}_{w})$ for all $n \in \mathbb{N}, i \in \mathbb{Z}$ with $i \geq - 1$.\\
With this notation one has $G_{n, -1} = G_{n}$. Further $G_{n, 0} = G_{n}$ since the extension $\mathbb{Q}_{w}(p^n)/\mathbb{Q}_{w}$ is totally ramified.\\~\\
\textbf{Lemma 2.3.9 :} Let $n$ be a positive integer. Assume that $P, Q \in \mathfrak{F}[p^{n}] - \mathfrak{F}[p^{n-1}]$.\\
Then there exists $\tau \in G_{n}$ such that $\tau(P) = Q$.\\~\\
\textbf{Proof :} First note that for all $\sigma \in G_{n}$, $\sigma (P) \in \mathfrak{F}[p^{n}] - \mathfrak{F}[p^{n-1}]$. Thus we have a map 
\[\phi : G_{n} \to \mathfrak{F}[p^{n}] - \mathfrak{F}[p^{n - 1}]\] given by $\phi (\sigma):= \sigma (P)$.\\
This map is injective because $P$ generates $\mathbb{Q}_{w}(p^{n})$ over $\mathbb{Q}_{w}$. \\
Now $|G_{n}| = q^{n - 1}(q - 1)$ and $|\mathfrak{F}[p^{n}] - \mathfrak{F}[p^{n - 1}]| = q^{n} - q^{n - 1} = q^{n - 1}(q - 1)$. \\
Hence the map must be onto.\\
The lemma follows from here. \hfill $\square$\\~\\
\textbf{Remark 2.3.10 :} Let $n$ be a positive integer and $P$ be a point of order $p^{n}$.\\
Let $0 \leq r \leq n$ and assume that $\sigma \in G_{n}$. Then $(\sigma(P) \, \ominus_{\mathfrak{F}} \,P) \in  \mathfrak{F}[p^{r}] $ if and only if $\sigma \in \text{Gal}\,(\mathbb{Q}_{w}(p^{n})|\mathbb{Q}_{w}(p^{n - r}))$. \\
If $(\sigma(P) \, \ominus_{\mathfrak{F}} \, P) \in \mathfrak{F}[p^{r}]$ then $\sigma([p^{r}]_{\mathfrak{F}}(P)) = [p^{r}]_{\mathfrak{F}}(P)$. Since $[p^{r}]_{\mathfrak{F}}(P)$ has order $p^{n- r}$ it generates $\mathbb{Q}_{w}(p^{n- r})$. Hence the forward implication.\\
Conversely $\sigma \in \text{Gal}(\mathbb{Q}_{w}(p^n)|\mathbb{Q}_{w}(p^{n- r}))$ implies $\sigma([p^r]_{\mathfrak{F}}(P)) = [p^r]_{\mathfrak{F}}(P)$ and hence the reverse implication.\\~\\
\textbf{Lemma 2.3.11 :} Let $n$ be a positive integer and $i$ be an integer $\geq - 1 $.\\ Let $P$ be a point of order $p^{n}$ and assume that $\sigma \in G_{n} - \{\text{Id}\}$.\\
If $(\sigma(P)\,\ominus_{\mathfrak{F}}\,P)$ has order $p^r$ for some $1 \leq r \leq n$, 
then $\sigma \in G_{n, i}$ for $0 \leq i \leq q^{n-r} - 1$ but $\sigma \notin G_{n,i}$ for $i = q^{n-r} $.\\~\\
\textbf{Proof :} Put $\Delta = \sigma(P)\,\ominus_{\mathfrak{F}}\,P$. Then $\sigma(P) = P\,\oplus_{\mathfrak{F}}\,\Delta $. Since $\sigma \neq \text{Id}$ and $P$ generates $\mathbb{Q}_{w}(p^n)$, we conclude $\Delta \neq O$, where $O$ is the identity element of the group $\mathfrak{F}[p^{n}]$. \\
Let $\mathfrak{p}$ be the unique maximal ideal in the ring of integers of $\mathbb{Q}_w(p^n)$ and let $v_{\mathfrak{p}}$ denote the valuation associated to it.\\
From the power series expansion of the formal group law we have $\sigma(P) = P + \Delta + P\Delta a$ for some integal element $a$. Thus using ultrametric triangle inequality $|\sigma(P) - P|_p = |\Delta|_p$. Hence $v_{\mathfrak{p}}(\sigma(P) - P) = v_{\mathfrak{p}}(\Delta)$.\\ But $\Delta$ has order $p^{r}$ and it generates $\mathbb{Q}_{w}(p^r)$ over $\mathbb{Q}_{w}$.
Note that $[\mathbb{Q}_{w}(p^{n}):\mathbb{Q}_{w}(p^r)] = q^{n-r}$ and the extension $\mathbb{Q}_{w}(p^{n})/\mathbb{Q}_{w}$ is totally ramified. Now use remark 2.3.5 to conclude that $v_{\mathfrak{p}}(\Delta) = q^{n- r}$.\\
From here the lemma follows just by using the definition of ramification groups and noting that $P$ is a generator for $\mathfrak{p}$.\hfill $\square$ \\~\\
\textbf{Proof of Theorem 2.1.1 (ii):} Let $k$ and $i$ be integers with $ 1 \leq k \leq n $ and $ q^{k - 1} \leq i \leq q^{k} - 1$.\\
Now from lemma 2.3.11 \[G_{n, i} = \{\sigma \in G_{n}\,|\,\text{order}(\sigma (P) \ominus_{\mathfrak{F}} P) \leq p^{n-k}\} = \text{Gal} (\mathbb{Q}_{w}(p^n)|\mathbb{Q}_{w}(p^{k}))\] where the last equality follows from remark 2.3.10 .\\
This proves theorem 2.1.1(ii). \hfill $\square$\\~\\ 
\textbf{Lemma 2.3.12 :} Let $n \in \mathbb{N}$. Then $\text{Gal}(\mathbb{Q}_{w}(p^{n})|\mathbb{Q}_{w}(p^{n-1}))$ lies in the center of $G_{n}$.\\~\\
\textbf{Proof :} Lemma 2.3.6 gives that $G_{1}$ is abelian. So it is enough to prove the lemma for $n \geq 2$.\\
Using theorem 2.1.1(ii) we conclude that $ G_{n, q^{n-1} - 1} = \text{Gal} \,(\mathbb{Q}_{w}(p^{n})|\mathbb{Q}_{w}(p^{n-1}))$ and $G_{n, q^{n-1}} = \{\text{Id}\}$.\\
Let $s \in G_{n,0}$ and $t \in G_{n,i}$ for some $i\geq 1$. Then $sts^{-1} \in G_{n,i}$.\\
$\theta_{i}$ and $\theta_{0}$ be the maps as in section 2, chapter 1.\\
From lemma 1.2.2 we know that \[\theta_{i} (sts^{-1}) = \theta_{0}(s)^{i}\theta_{i}(t).\]
Now $|G_{n,0}/G_{n,1}| = [\mathbb{Q}_{w}(p^n) : \mathbb{Q}_{w}] / [\mathbb{Q}_{w}(p^{n}) : \mathbb{Q}_{w}(p)] = (q - 1)$.\\
The map $\theta_{0}$ is a homomorphism with domain $G_{n, 0} / G_{n, 1}$ whose order is $(q - 1)$.\\
Note that $(q - 1)|(q^{n-1} - 1)$, since $n \geq 2$.\\
Put $i = q^{n-1} - 1$. Then $\theta_{0}({s})^{i} = 1$ for any $s \in G_{n, 0}$.\\
Thus \[\theta_{i} (sts^{-1}) = \theta_{i} (t).\]
We have already observed that $G_{n, i+1} = \{\text{Id}\}$. This along with the identity above and the fact that $\theta_{i}$ is an injection gives \[sts^{-1} = t.\]
True for all $s \in G_{n, 0} $ and $t \in G_{n, i}$.\\
Now $G_{n, 0} = G_{n} = \text{Gal}(\mathbb{Q}_{w}(p^n)|\mathbb{Q}_{w})$ and $G_{n, i} = \text{Gal}(\mathbb{Q}_{w}(p^n)|\mathbb{Q}_{w}(p^{n-1}))$.\\
Hence the lemma. \hfill $\square$\\~\\ 
\textbf{Notation :} By lemma 2.3.6 $G_{1} \cong \mathbb{Z}/(q - 1)\mathbb{Z}$. Let $\widetilde{\tau}$ be a fixed generator of this cyclic group.\\~\\
\textbf{Lemma 2.3.13:}  There exist a sequence $\{ \tau_{n} \}_{n \geq 1}$ such that $\tau_{n} \in G_{n}$ and the following holds :\\
a) Restriction of $\tau_{n}$ to $\mathbb{Q}_{w}(p^m)$ is $\tau_{m}$ for all $ 1 \leq m \leq n$ and $\tau_{1} = \widetilde{\tau}$.\\
b) Order of $\tau_{n}$ is $(q-1)$.\\
c) $\tau_{n}$ is in the centre of $G_{n}$.\\~\\
\textbf{Proof:} We construct a $\tau_{n}$ by induction on $n$.\\
Put $\tau_{1} = \widetilde{\tau}$. From choice of $\widetilde{\tau}$ and by lemma 2.3.6 we conclude that it satisfies all the requirements.\\
Assume that we have constructed $\tau_{n}$ for all $1 \leq n \leq m$ for some positive integer $m$. We want to construct for $\tau_{m+1}$.\\
Let $s$ be any extension of $\tau_{m}$ to $\mathbb{Q}_{w}(p^{m+1})$. Note that $[\mathbb{Q}_w(p^{m+1}):\mathbb{Q}_w(p^m)] = q$. Since $s$ restricts to $\tau_m$ whose order is $(q-1)$ one must have $(q-1) \,| \,\text{ord}(s)$ and order of $s$ is $(q - 1)p^{d}$ for some integer $d$ with $0 \leq d \leq 2$.\\
Put $t = s^{p^{d}}$. Clearly order of $t$ is $(q - 1)$ and it restricts to $\tau_{m}^{p^{d}}$.\\
Note that $\gcd (q - 1, p) = 1$. So there exists a positive integer $a$ such that $p^da \equiv 1 \mod (q - 1)$.\\
Put $\tau_{m +1} = t^{a}$. We would like to show that this choice of $\tau_{m + 1}$ has all the desired properties.\\
First observe that $\tau_{m +1}$ restricts to $\tau_{m}^{p^da} = \tau_{m}$. Further it has order $(q - 1)$ since $t$ has order $q - 1$ and $\gcd (a, q - 1) = 1$.\\
Let $\sigma \in G_{m + 1}$. We would like to show that $\sigma \tau_{m+1} \sigma^{-1} = \tau_{m+1}$.\\ 
Let $P$ be a point of order $p^{m + 1}$. It generates $\mathbb{Q}_{w}(p^{m + 1})$ over $\mathbb{Q}_{w}$. Thus if $s, t \in G_{m + 1}$, then to verify $s = t$ it is enough to check $s(P) = t(P)$.\\
By construction $\tau_{m+1}$ restricts to $\tau_{m}$. Using property c) of $\tau_{m}$ and the fact that $[p]_{\mathfrak{F}}(P) \in \mathfrak{F}[p^m]$, we have \[\sigma \tau_{m+1} \sigma^{-1} ([p]_{\mathfrak{F}}(P)) = \tau_{m+1} ([p]_{\mathfrak{F}}(P)).\]
So $\sigma \tau_{m+1}\sigma^{-1} (P) = \tau_{m+1}(P) \, \oplus_{\mathfrak{F}} \, \Delta$ for some $\Delta \in \mathfrak{F}[p]$.\\
Now by arguments as in proof of lemma 2.3.8 (here $\tau_{m+1}(P)$ is the point of order $p^{m+1}$) there is $t' \in \, \text{Gal} \, (\mathbb{Q}_{w}(p^{m+1}) | \mathbb{Q}_{w}(p^{m}))$ such that \[t'(\tau_{m+1}(P)) = \tau_{m+1}(P) \, \oplus_{\mathfrak{F}} \, \Delta.\]
Then $\sigma \tau_{m+1} \sigma^{-1}(P) = t'\tau_{m+1}(P)$. Hence $\sigma \tau_{m+1} \sigma^{-1} = t'\tau_{m+1}$. \\
Raise both sides to the power $p$. Since $t'$ lies in the center and $(t^{'})^{p} = \,\text{Id}$ we conclude that \[\sigma \tau_{m+1}^{p} \sigma^{-1} = \tau_{m+1}^{p}.\]
But $\gcd (p, q - 1) = 1$ . So there is a positive integer $b$ such that $ pb \equiv 1 $ mod$(q - 1)$.\\
Now raising to the power $b$ we conclude $\sigma \tau_{m+1} \sigma^{-1} = \tau_{m+1}$. \\
True for all such $\sigma$. So $\tau_{m+1}$ lies in the center of $G_{m+1}$.\\
Thus c) holds and $\tau_{m+1}$ is an extension of $\tau_{m}$ of desired kind.\\
Thus inductively we can construct a sequence $\{\tau_n\}_{n \in \mathbb{N}}$ having all the required properties.\hfill $\square$\\~\\
\textbf{Remark 2.3.14 :} i) Let $n$ be a positive integer and let $ 1 \leq m \leq n$. \\
Then using the lemma above one can put $\mathbb{Z} /(q - 1)\mathbb{Z}$ module structure on $\mathfrak{F}[p^{n}]$ such that $G_{n}$ acts via module automorphisms and $\mathfrak{F}[p^{m}]$ is a submodule wrt this module structure.\\
ii) For rest of this section we fix a $\tau_{n}$ for each positive integer $n$ as described in lemma.\\~\\
\textbf{Lemma 2.3.15 :} Let $n$ be a positive integer and $\tau_{n}$ be as above.\\
Let $P$ be a point of order $p^{n}$. Put $\tau_{n}(P) = Q$. Then $\{P, Q\}$ forms a basis of the $\mathbb{Z}/p^n\mathbb{Z}$ module $\mathfrak{F}[p^{n}]$ .\\~\\
\textbf{Proof:} We know $E[p^{n}]$ is a free $\mathbb{Z}/p^{n}\mathbb{Z}$ module of rank $2$.\\
By the isomorphism $(2.2.2)$, $\mathfrak{F}[p^{n}]$ is also a free $ \mathbb{Z}/p^{n}\mathbb{Z}$ module of rank $2$. Let $\{A_{1}, A_{2}\}$ be a basis for $\mathfrak{F}[p]$ over $\mathbb{Z}/p^{n}\mathbb{Z}$.\\
Say $P = [a_1]_{\mathfrak{F}}(A_1)\, \oplus_{\mathfrak{F}}\, [a_2]_{\mathfrak{F}}(A_2)$ where $a_1$ and $a_2$ are integers in the interval $[0, p^n-1]$. Since $P$ has order $p^n$ at least one of $a_1$ and $a_2$ is not divisible by $p$. Without loss of generality one can assume that $a_1$ is not divisible by $p$. Put $Q_1 = A_2$. Then it is easy to see that
$P$ and $Q_{1}$ forms a basis of the $\mathbb{Z}/p^n\mathbb{Z}$ module $\mathfrak{F}[p^{n}]$.\\
Let $a , b$ be integers between $0$ and $p^{n} - 1$ such that $Q = [a]_{\mathfrak{F}}(P) \, \oplus_{\mathfrak{F}} \, [b]_{\mathfrak{F}}(Q_{1})$.\\
To prove the lemma it is is enough to show that $b$ is an unit in $\mathbb{Z}/p^{n}\mathbb{Z}$ ie $p \nmid b$.\\
If possible assume that $p \mid b$. Then $[p^{n-1}]_{\mathfrak{F}}(Q) = [ap^{n-1}]_{\mathfrak{F}}P$. But 
\[\begin{split}
[p^{n-1}]_{\mathfrak{F}}(Q) = [p^{n-1}]_{\mathfrak{F}}(\tau_{n}(P)) & = \tau_{n}([p^{n -1}]_{\mathfrak{F}}(P))\\
                                                                    & = \tau_{1}([p^{n -1}]_{\mathfrak{F}}(P))\\
\end{split}
\] where the last equality follows from the observation $[p^{n -1}]_{\mathfrak{F}}(P) \in \mathbb{Q}_{w}(p)$ and the compatibility condition of lemma 2.3.13.\\
Thus $\tau_{1}([p^{n-1}]_{\mathfrak{F}}(P)) = [a]_{\mathfrak{F}}([p^{n-1}]_{\mathfrak{F}}(P))$. Clearly $\tau_1([p^{n-1}]_{\mathfrak{F}}(P)) \neq 0$. Hence $p \nmid a$. Let $d$ be order of the image of $a$ in the unit group of $\mathbb{Z}/p\mathbb{Z}$.\\
Then $\tau_{1}^{d}([p^{n -1}]_{\mathfrak{F}}(P)) = [p^{n - 1}]_{\mathfrak{F}}(P)$. Since $[p^{n -1}]_{\mathfrak{F}}(P)$ has order $p$ it generates $\mathbb{Q}_{w}$ and we conclude that $\tau_{1}^{d}$ is identity.\\
This is a contradiction since $d \leq p - 1 < p^{2} - 1 = \text{ord}(\tau_{1})$.\\
So one cannot have $p\,|\,b$.\\
The lemma follows from here.\hfill $\square$\\~\\
\textbf{Lemma 2.3.16 :} Let $n$ be an integer $\geq 2$.\\
Let $P$ be a point of order $p^{n}$ and assume that $\tau_{n}(P) = Q$.\\
Define $s_n, t_n \in G_{n}$ by putting $s_{n}(P) = P\,\oplus_{\mathfrak{F}}\, [p]_{\mathfrak{F}}(P)$ and $t_{n}(P) = P\,\oplus_{\mathfrak{F}}[p]_{\mathfrak{F}}(Q)$.\\
By lemma 2.3.9 such elements exist since both of $P \,\oplus_{\mathfrak{F}}\, [p]_{\mathfrak{F}}(P)$ and $P \,\oplus_{\mathfrak{F}}\,[p]_{\mathfrak{F}}(Q)$ are in $\mathfrak{F}[p^{n}] - \mathfrak{F}[p^{n-1}]$.\\
Let $<s_{n}>$ and $<t_{n}>$ denote the groups generated by $s_n, t_n$ respectively.\\
Then:\\
a) $|<s_n>|$ = $|<t_n>|$ = $p^{n-1}$,\\
b) $ <s_n> \cap <t_n> = \text{Id}$,\\
c) $s_nt_n = t_ns_n$.\\~\\
\textbf{Proof:} First note that 
\[s_{n}(Q) = s_{n}(\tau_{n}(P)) = \tau_{n}(s_{n}(P)) = \tau_{n}(P \oplus_{\mathfrak{F}} [p]_{\mathfrak{F}}(P)) = Q \,\oplus_{\mathfrak{F}} \,[p]_{\mathfrak{F}}(Q)\] since $\tau_{n}$ is in the center.\\
Now, \[s_{n}(t_{n}(P)) = s_{n}(P\,\oplus_{\mathfrak{F}}\,[p]_{\mathfrak{F}}(Q)) = P \, \oplus_{\mathfrak{F}} \, [p]_{\mathfrak{F}} ( P \, \oplus_{\mathfrak{F}} \, Q) \, \oplus_{\mathfrak{F}} [p^{2}]_{\mathfrak{F}}(Q)\] where we are using the expressions of $s_{n}(P)$ and $s_{n}(Q)$ and the fact that $[p]_{\mathfrak{F}}$ commutes with the action of Galois group.\\
Similarly using expression for $t_n(P)$ we conclude that \[t_{n}(s_{n}(P)) = t_{n}(P\,\oplus_{\mathfrak{F}}\,[p]_{\mathfrak{F}}(P)) = P \,\oplus_{\mathfrak{F}}\,[p]_{\mathfrak{F}}(P \,\oplus_{\mathfrak{F}}\,Q)\,\oplus_{\mathfrak{F}} [p^{2}]_{\mathfrak{F}}(Q).\]
Since $P$ generates the extension we conclude c).\\~\\
\textbf{Claim 2.3.17:} Let $m \in \mathbb{N}$ and $a$ be an integer such that $ 2 \leq a \leq p^{m}$. Then $p^{m+2} \mid p^{a}\binom{p^m}{a}$.\\~\\
\textbf{Proof of claim: } If $a = p^m$ then the claim follows easily since $p^{m} \geq m+2$. So assume that $ a < p^m$.\\
We shall use the notation $v_p(\cdot)$ to denote the $p$-adic valuation of a natural number.\\
Say, $v_p(a) = v < m$. Then $v_p(p^m - a) = v$. If one writes $a$ and $p^m - a$ in base $p$ both of them have a string of exactly $v$ consecutive zeroes at the end. But when we add them the sum represented in base $p$ is $1$ followed by a string of $m$ zeroes. So the number of carries is exactly $m - v$.\\
Using Kummer's theorem for $p$-adic valuation of binomial coefficients (see \cite{wiki}) one concludes that $\binom{p^m} {a}$ is divisible by $p^{m -v}$.\\
Since $a \geq 2$ and $p \geq 5$ one has $a - v_{p}(a) \geq 2$. The claim follows from here. \hfill $\square$ \\~\\
Hence for all $1 \leq m \leq n - 1$, \[(\text{Id}\,\oplus_{\mathfrak{F}}\, [p]_{\mathfrak{F}} )^{p^m} (P) \, \ominus_{\mathfrak{F}} \, ( P \, \oplus_{\mathfrak{F}} \, [p^{m+1}]_{\mathfrak{F}}P ) \in \mathfrak{F}[p^{n - (m+2)}]\] and \[( \text{Id}
 \, \oplus_{\mathfrak{F}} \,[p]_{\mathfrak{F}}\tau_{n})^{p^m} (P) \, \ominus_{\mathfrak{F}} \, ( P \, \oplus_{\mathfrak{F}} \, [p^{m+1}]_{\mathfrak{F}}Q ) \in \mathfrak{F}[p^{n - (m+2)}]\]
where we are using the notation $\mathfrak{F}[p^{- 1}] = \{ 0 \}$.\\
Let $s, t \in G_n$. We know that $s = t$ if and only if $s(P) = t(P)$. Now looking at the second terms in the corresponding expressions and using the linear independence property of $P$ and $Q$ from lemma $2.3.15$ we have that \[s_{n}^{p^m} \neq t_{n}^{p^m}\tag{2.3.8}\] for all $1\,\leq\,m\,\leq\,n - 2$ and for all $m$ in this range \[\text{Id} \notin \{s_{n}^{p^m}, t_{n}^{p^m}\}.\tag{2.3.9}\]
Further \[s_n^{p^{n-1}} = t_n^{p^{n-1}} = \text {Id}.\tag{2.3.10}\]
We already know that degree of $\mathbb{Q}_{w}(p^n)/\mathbb{Q}_{w}(p)$ is a power of $p$. By remark 2.3.10 we have $s_n, t_n \in \text{Gal} (\mathbb{Q}_{w}(p^n) | \mathbb{Q}_w(p))$. Thus order of $s_n$ and $t_n$ must be a power of $p$.\\
This along with $(2.3.8)$ and $(2.3.10)$ proves part a.\\
Now assume that order of $<s_n> \cap <t_n>$ is $p^{m}$ for some $0 \leq m \leq n - 1$. Then $s_n^{p^{n - m - 1}} = t_n^{p^{n - m -1}}$. $(2.3.9)$ forces $m = 0$ or $m = n -1$. Clearly $s_n \neq t_n$. So we cannot have $m = n -1$. Hence $m = 0$. This proves part b.\\
We have already proved part c.\\
This completes the proof of the lemma. \hfill $\square$\\~\\
\textbf{Proof of theorem 2.1.1 (i) :} By lemma 2.3.16 \[< s_n > \times < t_n > \hookrightarrow G_{n}.\]
By lemma 2.3.13 $\tau_{n}$ has order $q - 1$ which is coprime to $p$. From this we conclude $<\tau_{n}>$ intersects $< s_{n} > \times < t_{n} >$ trivially. Further $\tau_n$ commutes with $s_{n}$ and $t_{n}$. Hence \[< s_{n} > \times < t_{n} > \times < \tau_{n} > \hookrightarrow G_{n}.\]
Note that both the sets have same cardinality, namely $q^{n - 1}(q - 1)$.\\
So the injection must be a bijection.\\
Thus $G_{n} \cong < s_{n} > \times < t_n > \times <\tau_n>$. This proves the first part of the assertion.\\
The last part of assertion already follows from lemma 2.3.8.\hfill $\square$\\~\\
\textbf{Proof of theorem 2.1.1 (iii) :} Let $M$ be a positive integer co-prime to $p$.\\
Let $P$ be a point of order $p^{n}$.\\
Then $[M]_{\mathfrak{F}}(P)$ also has order $p^{n}$.\\
By lemma 2.3.8 we know there is a $\sigma \in G_{n}$ such that $\sigma (P) = [M]_{\mathfrak{F}}(P)$.\\ 
Put $Q = \tau_n(P)$. Now \[\sigma(Q) = \sigma (\tau_{n}(P)) = \tau_{n}(\sigma(P)) = \tau_{n} ([M]_{\mathfrak{F}}(P)) = [M]_{\mathfrak{F}}(Q).\]
Since $P$ and $Q$ generate the whole group (lemma 2.3.15) we conclude that $\sigma$ acts via multiplication by $[M]_{\mathfrak{F}}$.\\
Thus we have constructed an element in the Galois group which acts like multiplication by $M$. Can be done for any such $M$.\\
Hence the first part of the assertion.\\
The second part is already a consequence of lemma 2.3.8.\hfill $\square$
\section {Further results}
\textbf{Proposition 2.4.1 :} Let $M$  be a positive integer prime to $p$.\\
Let $E$ be an elliptic curve over $\mathbb{Q}_{w}$.\\
Then $\mathbb{Q}_{w}(M)/\mathbb{Q}_{w}$ is an unramifed extension.\\~\\
\textbf{Proof:} This is a standard result in algebraic number theory. For a proof see \cite{sil1} chapter VII, proposition 4.1. \hfill $\square$ \\~\\
\textbf{Proposition 2.4.2 :} Assume that, we are in the set up of section 1. \\
Let $M$ be coprime to $p$, and suppose that $n$ is a nonnegative integer. \\
Put $ N = p^{n}M$.\\
The following statements hold :\\
i) The composition $\mathbb{Q}_{w}(p^{n})\mathbb{Q}_{w}(M)$ is $\mathbb{Q}_{w}(N)$.\\
ii) The extension $\mathbb{Q}_{w}(N) / \mathbb{Q}_{w}(p^{n}) $ is unramified, and the extension $\mathbb{Q}_{w}(N) / \mathbb{Q}_{w}(M) $ is totally ramified.\\
iii) Restriction to $\mathbb{Q}_{w}(p^{n})$ induces an isomorphism of groups \[ \text{Gal} \, ( \mathbb{Q}_{w}(N) | \mathbb{Q}_{w}(M) ) \to \text{Gal} \, ( \mathbb{Q}_{w}(p^{n}) | \mathbb{Q}_{w}). \]
In particular $\mathbb{Q}_{w}(N) / \mathbb{Q}_{w}(M) $ is abelian.\\
iv) If $ n \geq 1$, then \[\text{Gal} \left(\mathbb{Q}_{w}(N)|\mathbb{Q}_{w}(\frac {N} {p})\right) \cong \text{Gal} (\mathbb{Q}_{w}(p^{n})|\mathbb{Q}_{w}(p^{ n - 1}))\]
by the restriction map.\\~\\
\textbf{Proof :} i) Follows from the fact that $E[N] \cong E[p^{n}] \oplus E[M]$.\\
ii) and iii) Note that $\mathbb{Q}_{w}(p^n)/\mathbb{Q}_{w}$ is totally ramified and $\mathbb{Q}_{w}(M)/\mathbb{Q}_{w}$ is unramified. Now put $ F = \mathbb{Q}_{w}$, $ K = \mathbb{Q}_{w}(p^{n})$ and $ L = \mathbb{Q}_{w}(M) $ in lemma 1.2.1.\\
iv) Note that by lemma 1.2.1 and by part (i) we have $\mathbb{Q}_{w}(p^{n - 1}M)/ \mathbb{Q}_{w}(p^{n - 1})$ is unramified. Clearly $\mathbb{Q}_{w}(p^n)/\mathbb{Q}_{w}(p^{n - 1})$ is totally ramified. Put $F = \mathbb{Q}_{w}(p^{n - 1})$, $ K = \mathbb{Q}_{w}(p^{n})$ and $L = \mathbb{Q}_{w}(p^{n - 1}M)$ in lemma 1.2.1 to conclude the present lemma. \hfill $\square$ \\~\\
\subsection{Lemmas on roots of unity}
\textbf{Notation :} i) $N$ will denote a positive integer. We shall write $N = p^nM$ where $n$ is an integer $\geq 0$ and $M$ is a positive integer coprime to $p$.\\
ii) For any $N \in \mathbb{N}$, $\mu_{N}$ will denote the group of roots of unities whose order divide $N$. $\mu_{p^\infty} \subseteq \overline{\mathbb{Q}}$ denote the group of roots of unity whose orders are power of $p$. We write $\mu_{\infty}$ for all the roots of unity in $\overline{\mathbb{Q}}$.\\
We shall continue to use this notation till the end of the thesis.\\~\\ 
With this notation we have the following lemma :\\~\\
\textbf{Lemma 2.4.3 : } $\mathbb{Q}_{w}(N) \cap \mu_{p^{\infty}} = \mu_{p^{n}} $.\\~\\
\textbf{Proof :} By Weil pairing $\mu_{p^n} \subseteq \mathbb{Q}_{w}(N) \cap \mu_{p^{\infty}}$.\\
We just need to show the inclusion in the other direction.\\
First we verify $\mathbb{Q}_{w}(p^{n}) \cap \mu_{p^{\infty}} \subseteq \mu_{p^{n}}$.\\ 
Let $\zeta \in \mathbb{Q}_{w}(p^{n}) \cap \mu_{p^{\infty}}$. Assume that $\zeta$ has order $p^m$ for some $m \geq 0$. We want to show that $m \leq n$.\\
From standard theory of cyclotomic extensions we know that if $\zeta_{m}$ is a primitive $p^m$- th root of unity ($m\geq 1$) then the extension $\mathbb{Q}_{p}(\zeta_m)/\mathbb{Q}_p$ is totally ramified extension with Galois group $\mathbb{Z}/(p - 1)\mathbb{Z} \times \mathbb{Z}/p^{m-1}\mathbb{Z}$. An application of lemma 1.2.1 shows that $\mathbb{Q}_{w}(\zeta_m)/\mathbb{Q}_{w}$ is a totally ramified extension with same Galois group.\\
Say, $n = 0$. Then using the observation above we must have $ m < 1$. Thus $ m = 0 = n$ and we are done in this case.\\
Now let $n \geq 1$.\\
Since $\zeta \in \mathbb{Q}_{w}(p^n)$ we have a surjective homomorphism $\text{Gal}(\mathbb{Q}_{w}(p^n)|\mathbb{Q}_{w}) \to \text{Gal}(\mathbb{Q}_{w}(\zeta)|\mathbb{Q}_{w})$. \\
The group on the left hand side is $\mathbb{Z}/( q - 1)\mathbb{Z} \times (\mathbb{Z}/ p^{n - 1}\mathbb{Z})^{2}$ and the group on the right hand side is $\mathbb{Z}/(p - 1)\mathbb{Z} \times \mathbb{Z}/p^{m - 1}\mathbb{Z}$. Comparing the $p$-part we conclude that $m \leq n$ as desired.\\
Now let $\zeta \in \mathbb{Q}_{w}(N)$ and its order is $p^m$ for some integer $m \geq 0$. If possible assume that $ m > n$.\\
Note $\mathbb{Q}_{w}(p^n) \subseteq \mathbb{Q}_{w}(p^n)(\zeta) \subseteq \mathbb{Q}_{w}(p^m)$ and the third one is a totally ramified extension of the first one. Hence the extension $\mathbb{Q}_{w}(p^n)(\zeta)/\mathbb{Q}_{w}(p^n)$ is totally ramified.\\
But $\mathbb{Q}_{w}(p^n) \subseteq \mathbb{Q}_{w}(p^n)(\zeta) \subseteq \mathbb{Q}_{w}(N)$ and the third one is an unramified extension of the first. Thus $\mathbb{Q}_{w}(p^n)(\zeta)/\mathbb{Q}_{w}(p^n)$ is unramified.\\
Hence this extension must be trivial.\\
Therefore $\zeta \in \mathbb{Q}_{w}(p^n)$. So $m \leq n$. A contradiction !\\
This contradiction proves the lemma. \hfill $\square$\\~\\ 
\textbf{Lemma 2.4.4 : } Let $ n \geq 1$.\\
If $ \psi \in \text{Gal} (\mathbb{Q}_{w}(N) | \mathbb{Q}_{w}(\frac {N} {p} ) )$ and $\alpha \in \mathbb{Q}_{w}(N) - \{ 0 \}$ such that $\psi(\alpha) / \alpha \in \mu_{\infty}$. Then,
\[ \frac {\psi(\alpha)} {\alpha} \in \mu_{Q(n)}, \, \text{where}\; Q(n) = \begin{cases}
                                                                                                                        q  &  \text{if} \,n \geq 2. \\
                                                                                                                        (q - 1)q & \text{if}\, n = 1.\\ 
                                                                                                                  \end{cases} \]
\textbf{Proof:} Put $\frac {\psi(\alpha)} {\alpha} = \beta$.\\
By hypothesis $\beta \in \mu_{\infty}$. Say its order is $N' = p^{n'} M'$ where $n' \geq 0$ and $p\nmid M'$.\\
Say $\beta^{p^{n'}} = \xi$. Then order of $\xi$ is $M'$. Since $\gcd (p, M') =1$ we conclude that the extension $\mathbb{Q}_{p}(\xi)/\mathbb{Q}_{p}$ is unramified. \\
We know that compositum of unramified extensions are unramified and subextensions of unramified extensions are unramified. Thus $\mathbb{Q}_{w}(M)(\xi)/\mathbb{Q}_{w}(M)$ is unramified. But $\mathbb{Q}_{w}(M)$ is the maximal unramified subextension of $\mathbb{Q}_{w}$ in $\mathbb{Q}_{w}(N)$. So $\xi \in \mathbb{Q}_{w}(M)$. In particular, $\psi(\xi) = \xi$.\\
Note that $\beta^{M'} \in \mu_{p^{\infty}}$. Hence by previous lemma $\beta^{M'} \in \mu_{p^n}$. Thus $\beta^{pM'} \in \mu_{p^{n - 1}}$ and it is also fixed by $\psi$.\\
$\gcd (p, M') = 1$. Therefore there are integers $a, b$ such that $ap^{n'} + bM' = 1$.\\
So, $\beta = \beta^{ap^{n'}}\beta^{bM'} = \xi^{a}\beta^{bM'}$. Therefore $\beta^{p} = \xi^{ap}(\beta^{pM'})^b$ and it is fixed by $\psi$.\\
Let $t$ be the order of $\psi \in \text{Gal}(\mathbb{Q}_{w}(N)/\mathbb{Q}_{w}(N/p))$. Then $\psi^{t}(\alpha^p) = \alpha^p$.\\
But $\psi(\beta^p) = \beta^p$. Hence $\psi(\psi(\alpha^p)) = \psi(\alpha^p)\beta^p$. Using this relation iteratively one has $\beta^{tp} = 1$.\\
From theorem 2.1.1 it follows that $ t | q - 1$ if $ n = 1$ and $t | p$ if $n \geq 2$.\\
Now the lemma follows from definition of $Q(n)$. \hfill $\square$

\chapter{Diophantine estimates}
\section{Introduction}
Let $K$ be a number field which has at least one real embedding.\\
Put $[K : \mathbb{Q}] = d$.\\
Let $E$ be an elliptic curve defined over $K$ which does not have complex multiplication.\\
Fix a nonzero prime ideal $\mathfrak{p}$ of $O_K$, the ring of integers of $K$ having the following properties:\\
i) $\mathfrak{p}$ does not ramify.\\
ii) $E$ has a super singular reduction at $\mathfrak{p}$ and the reduced curve does not have $j-$ invariant among $\{ 0 , 1728 \}$.\\
iii) $\mathfrak{p} \cap \mathbb{Z} = p\mathbb{Z}$ where $p$ is a positive prime $\geq \, 2^{d + 2}$  and
the natural Galois representation \[ \text{Gal} (\overline{\mathbb{Q}} / K) \to \text{Aut}_{\mathbb{Z}_{p}} E[p^{\infty}] \] is onto. where $E[p^{\infty}]$ is the $p$-adic Tate module and we are considering its automorphisms as $\mathbb{Z}_{p}$ module.\\~\\
By well known results we know all but finitely many prime ideals $\mathfrak{p}$ satisfy i and iii. Since $E$ does not have complex multiplication its $j-$invariant is not among $\{0, 1728\}$. Hence for all but finitely many places the $j-$invariant of the reduced curve is not also among $\{0, 1728\}$.
So by theorem 1.3.3, we know that at least one prime ideal exists satisfying all the three properties.\\
Let $f$ be degree of residue extension associated with the prime $\mathfrak{p}$. Note that, minimal Weierstrass model of $E$ at $\mathfrak{p}$ is defined over $\mathbb{Z}_{p^f}$.
\section{Local metric estimates}
Let $p,f$ be as in section 3.1 and $w = p^{2f}$, $q = p^2$. We shall think $E$ as an elliptic curve over $\mathbb{Q}_{w}$. It's minimal model at $\mathfrak{p}$ is defined over $\mathbb{Z}_w$ and its $j$-invariant is not among $0,1728$.\\
Assume that $N$ is a positive integer such that $N = p^{n}M$ where $n$ is an integer $\geq 0$ and $M$ is a positive integer with $\gcd ( M, p)  = 1$.\\
Let $\mathbb{Q}_{p}^{\text{nr}}$ denote the maximal unramified extension of $\mathbb{Q}_{p}$ inside $\overline{\mathbb{Q}_{p}}$ and let $\phi_{w}\in \text{Gal} \, (\mathbb{Q}_{p}^{\text{nr}} / \mathbb{Q}_{p} )$ denote the lift of Frobenius composed $f$ times.\\
Note that $\phi_{w} \in \text{Gal} \,(\mathbb{Q}_{p}^{\text{nr}} / \mathbb{Q}_{w})$ .\\
First we consider the case $ p \,\nmid\,N$. We have :\\~\\
\textbf{Lemma 3.2.1 :} Suppose $ p \nmid N$ and $\alpha \in \mathbb{Q}_{w}(N)$. Then $\alpha \in \mathbb{Q}_{p}^{\text{nr}}$ and
$$ |\phi_{w}(\alpha) - \alpha^{w}|_{p} \leq p^{-1} \max \{ 1 , |\phi_{w}(\alpha) |_{p} \} \max \{ 1 , |\alpha|_{p} \}^{w} \, .$$\\~\\
\textbf{Proof :} Let $L = \mathbb{Q}_{w}(N)$. Since $ p \nmid N$ we conclude that $L$ is an unramified extension of $\mathbb{Q}_{w}$. \\
Since $\alpha \in L$ the first part of the lemma follows.\\
To prove the second part at first we assume that $\alpha$ is an integer in $L$. \\
Then $(\phi_{w}(\alpha) - \alpha^{w})$ is in the maximal ideal $pO_{L}$ where $O_{L}$ is the ring of integers of $L$.\\
Therefore, $| \phi_{w}(\alpha) - \alpha |_{p} \, \leq p^{-1} $ .\\
Hence the inequality in the statement of the lemma holds trivially.\\
Now assume $\alpha \, \notin O_{L}$. Then $\alpha \, \neq 0 $ and $\alpha^{-1} \in O_{L}$ .\\
Using the computation above $|\phi_{w} (\alpha^{-1}) - \alpha^{-w}|_{p} \, \leq p^{-1}  $.\\
Notice that $\phi_{w}(\alpha^{-1}) - \alpha^{-w} = \alpha^{-w}\phi_{w}(\alpha)^{-1} ( \phi_{w}(\alpha) - \alpha^{w} ) $.\\
Hence $|\phi_{w} (\alpha) - \alpha^{w}|_{p} \, \leq p^{-1}|\phi_{w}(\alpha)|_{p}|\alpha|_{p}^{w}$ and thus the lemma holds true in this case also. \hfill $\square$ \\~\\
Now one considers the case $ p\,|\,N $. In this case we have  the following lemma  \\~\\
\textbf{Lemma 3.2.2 :} Suppose $p \, | \, N$ and $\alpha \, \in \mathbb{Q}_{w}(N)$. Then\\
$$ | \psi (\alpha)^{q} - \alpha^{q} |_{p} \, \leq p^{-1} \max \{ 1 , | \psi (\alpha) |_{p} \}^{q} \max \{ 1 , |\alpha|_{p} \}^{q}$$
for any $\psi \, \in \text{Gal}(\mathbb{Q}_{w}(N) / \mathbb{Q}_{w}(N/p) )$.\\~\\
\textbf{Proof :} Put $ K = \mathbb{Q}_{w}(p^{n})$ and $L = \mathbb{Q}(M)$.\\
Then $KL = \mathbb{Q}_{w}(N)$ by proposition 2.4.2.\\
First assume that $\alpha$ is an integer in $\mathbb{Q}_{w}(N)$.\\
Now $\psi|_{K} \in \text{Gal} ( K | \mathbb{Q}_{w}(p^{n - 1}) )$ .\\
By theorem 2.1.1, we have $\psi|_{K} \in G_{i}( K | \mathbb{Q}_{w} )$ for $ i = q^{ n - 1} - 1$.\\
Now $\psi$ fixes $L$.\\
Hence it is in $\text{Gal} ( KL | L)$.\\ 
Using lemma 1.2.1 we conclude $\psi \in G_{i}( \mathbb{Q}_{w}(N) | \mathbb{Q}_{w} )$ for $ i = q^{n - 1} - 1 $.\\
Let $\mathfrak{P}$ be the unique maximal ideal in the ring of integers of $\mathbb{Q}_{w}(N)$.\\
Thus \[ \psi(\alpha) - \alpha \in \mathfrak{P}^{q^{ n - 1}}. \] 
Now, $e (\mathbb{Q}_{w}(N)|\mathbb{Q}_{w}) = e (\mathbb{Q}_{w}(p^{n})|\mathbb{Q}_{w}) = q^{n - 1}(q - 1)$.\\
Therefore, \[(\psi(\alpha) - \alpha)^{q} \in \mathfrak{P}^{q^{n}} \subseteq \mathfrak{P}^{q^{n - 1}(q - 1)} = \mathfrak{P}^{e(\mathbb{Q}_w(N)|\mathbb{Q}_w)}.\]
Hence $| \psi(\alpha)^{q} - \alpha^{q} |_{p} \leq p^{ - 1}$.\\
Since $\alpha$ is an integer, so is $\psi(\alpha)$ and thus the inequality follows in this case.\\ 
Now assume that $ \alpha$ is not an integer in $\mathbb{Q}_{w}(N)$.\\
Then $\alpha \neq 0$ and $ \alpha^{ - 1} $ is an integer.\\
So by previous computations $ | \psi(\alpha)^{ - q} - \alpha^{- q} |_{p} \leq p^{ - 1} $.\\
But the left hand side of the identity is nothing but $ |\frac { \psi(\alpha)^q - \alpha^{q} } { \psi(\alpha)^{q} \alpha^{q}} |_{p}$.\\
So we have $ \vert \psi(\alpha)^q - \alpha^q |_{p} \leq p^{ - 1} | \psi(\alpha)|_{p}^{q} | \alpha|_{p}^{q} $ which implies the desired inequality in this case. \hfill $\square$ \\~\\
\section{A first global estimate }
Let $K$, $E$, $\mathfrak{p}$ and $p$ be as in introduction.\\
Note that $ p \geq 5$.\\
For this section we shall only need condition i and ii for $\mathfrak{p}$.\\
Let $ f $ be the local degree at $\mathfrak{p}$ .\\
Then completion of $K$ with respect to $\mathfrak{p}$ is isomorphic to $\mathbb{Q}_{p^{f}}$ (for the rest of the part we shall work with this embedding of $K$ in $\overline{\mathbb{Q}}_{p}$) and thus $E$ can be thought as an elliptic curve over $\mathbb{Q}_{w}$. Then $K(N) \subseteq \mathbb{Q}_{w}(N)$ for each positive integer $N$. So $|\cdot|_{p}$ will determine a place of $K(N)$. We call this place to be $v$.\\
Let $L$ be a finite Galois extension of $K$.\\
Let $u$ be a place of $L$.\\
An automorphism $\sigma \in \text{Gal} ( L | K )$ determines a new place $\sigma u$ of $L$ by the formula \[ | \alpha |_{\sigma u} = | \sigma^{ - 1}( \alpha ) |_u \] for all $ \alpha \in L$.\\~\\
First we handle the unramified case:\\~\\
\textbf {Lemma 3.3.1 :} Let $N$ be a positive integer.\\
Assume that $ p \nmid N$.\\
If $\alpha \in K(N) - \mu_{\infty}$ be a nonzero algebraic number then \[ h ( \alpha ) \geq \frac {\log (p^f / 2^{d})} {d(w + 1)} . \]
\textbf{Proof : } In the proof every field is considered as a subfield of $\overline{\mathbb{Q}}_{p}$.\\
Now by discussion in chapter 2, section 1 we conclude that  $\widetilde{\phi_{w}}$ acts as multiplication by $ \pm [p^{f}]$ on $\widetilde{E}$ .\\
Let $l$ be a prime with $ l \neq p$. Then $\widetilde{\phi_{w}}$ acts as multiplication by $\pm p^{f}$ on the $l$-adic Tate module $T_{l}(\widetilde{E})$.\\
But $T_{l}(\widetilde{E}) \cong T_{l}(E)$ by the isomorphism in $(2.2.2)$. \\
Note that this isomorphism commutes with the action of the absolute Galois group of $\mathbb{Q}_{w}$ and all the points of $T_{l}(E)$ are defined over $\mathbb{Q}_{w}^{ur}$.\\
Thus $\phi_{w}$ acts on the points of $T_{l}(E)$ and the isomorphism mentioned above gives that this action is nothing but multiplication by $\pm p^{f}$.\\
True for all such primes.\\
Thus $\phi_{w}$ acts on $E[N]$ as multiplication by $\pm p^{f}$.\\
Note that $\text{Aut}_{\mathbb{Z}}(E[N]) \cong \text{Gl}_{2}(\mathbb{Z}/N\mathbb{Z})$ and the action of $\phi_{w}$ is given by a scalar matrix.\\
Hence $\phi_{w}$ commutes with the action of $ G = \text{Gal} (K(N)|K)$. So it lies in the center of the Galois group.\\
Define $x = \phi_{w}(\alpha) - \alpha^{w} \in K(N)$.\\
If $x = 0$ then $h(\alpha) = h(\phi_{w}(\alpha)) = w h(\alpha)$ which contradicts the assumption on $\alpha$. So $ x \neq 0$.\\
Using product formula one has \[\sum_{u} d_{u} \log |x|_{u} = 0, \tag{3.3.1} \] where the sum is over all places of $K(N)$.\\
Let $v$ be as above. Say, $u$ is any finite place of $K(N)$ lying over $\mathfrak{p}$.\\
Then $u = \sigma^{- 1}v$ for some $\sigma \in G$.\\
The fact that $\phi_{w}$ lies in the center implies 
\begin{equation*}
\begin{split}
|x|_{u} & = | \sigma(\phi_{w}(\alpha)) - \sigma(\alpha^{w}) |_{v}\\ 
             & = | \phi_{w} (\sigma(\alpha)) - \sigma(\alpha)^{w} |_{v}\\
             & \leq p^{- 1} \max \{ 1, |\phi_{w}(\sigma(\alpha))|_{v} \} \max \{ 1, | \sigma(\alpha)|_{v}\}^{w} 
\end{split}
\end{equation*} 
where we are using lemma 3.2.1 .\\
Now 
\begin{equation*}
\begin{split}
|\phi_{w}(\sigma(\alpha))|_{v} &= | \sigma (\phi_{w}(\alpha))|_{v}\\
                                                   &= |\phi_{w}(\alpha)|_{u}
\end{split}
\end{equation*}
and $|\sigma(\alpha)|_{v}
              = |\alpha|_{u}$ .

Thus \[|x|_{u} \leq p^{- 1} \max\{ 1, |\phi_{w}(\alpha)|_{u} \} \max \{ 1, |\alpha|_{u} \}^{w}. \tag{3.3.2} \]
The estimate given above holds for any $u$ lying over $\mathfrak{p}$.\\
For other finite places we use the trivial estimate  \[ |x|_{u} \\ \leq \max \{ |\phi_{w}(\alpha)|_u , |\alpha^{w}|_u \} \\ \leq \max \{ 1, |\phi_{w}(\alpha)|_u \} \{ 1, |\alpha|_u \}^{w} .\tag{3.3.3}\]
If $u$ is an infinite place then we have the estimate \[ |x|_{u} \\ \leq 2 \max \{ |\phi_{w}(\alpha)|_u , |\alpha^w|_u \} \\ \leq 2 \max \{ 1, |\phi_{w}(\alpha)|_u \} \{ 1, |\alpha|_u \}^w .\tag{3.3.4}\]
Putting $(3.3.2)$, $(3.3.3)$ and $(3.3.4)$ in the product formula one deduces\\
\[
\begin{split}
 0 = & \sum_{ u | \mathfrak{p} } d_{u} \log |x|_{u}  + \sum_{ u | \infty} d_{u} \log |x|_{u} + \sum_{ u \nmid \infty, u \nmid \mathfrak{p} } d_{u} \log |x|_{u}\\ 
       &  \leq - \log p \sum_{u | \mathfrak{p} } d_{u} + \log 2 \sum_{ u | \infty} d_{u} + \sum_{u} (\log^{+}|\phi_{w}(\alpha)|_{u} + w \log^{+}|\alpha|_{u}) 
\end{split}
\] 
Note that $\sum_{u |\mathfrak{p}} d_u = f\frac {[K(N) : \mathbb{Q}]} {[K : \mathbb{Q}]}$ (since $\mathbb{Q}_{w}(N)$ is unramified over $\mathbb{Q}_w$ each of these places are unramified over $\mathfrak{p}$) and $\sum_{ u | \infty} d_u = [K(N) : \mathbb{Q}] $.\\~\\
So dividing both sides by $[K(N) : \mathbb{Q}]$ one obtains \[ 0 \leq - \frac {f\log p} {d} + \log 2 + h(\phi_{w}(\alpha)) + w h(\alpha) \]
and using $h(\phi_{w}(\alpha)) = h(\alpha)$ one concludes \[ h(\alpha) \geq \frac { \log ( p^f / 2^{d} ) } { d( 1 + w)} \] as desired. \hfill $\square$\\~\\
Next we consider the ramified case:\\~\\
\textbf{Lemma 3.3.2 :} Assume that $ p | N$.\\
Suppose $\psi \in \text{Gal} (\mathbb{Q}_{w}(N)|\mathbb{Q}_{w}(N/p))$, which we identify with its restriction to $K(N)$.\\
Put $ G = \text{Gal} (K(N)| K)$ and $G_{\psi} = \{ \sigma \in G | \sigma\psi\sigma^{-1} = \psi \}$.\\
Let $v$ be the place of $K(N)$ induced by $|\cdot|_p$.\\
Then\[| G_{\psi}v| \;\geq \; p^{-4} \,\frac {f[K(N) : K]} {d_{v}}\].\\
\textbf{Proof :} Let \[ H = \text{Gal} ( K(N) | K(N/p) ) .\] Clearly $H$ is a normal subgroup of $G$ and $ \psi \in H$.\\ 
Fix an isomorphism between $E[N] \cong (\mathbb{Z} / N \mathbb{Z} )^{2}$.\\
Then each automorphism of $E[N]$ (as an abelian group)  is an element of $\text{Gl}_{2}(\mathbb{Z} / N\mathbb{Z})$.\\
An automorphism which acts trivially on $E[N/p]$ is represented by an element of $ 1 + N/p \text{Mat}_{2} (\mathbb{Z} / p\mathbb{Z} )$.\\
Since the representation $ G \to \text{Aut}_{\mathbb{Z}}E[N] $ is injective we conclude that $ |H| \leq p^{4}$.\\
Consider the conjugation action of $G$ on itself.\\
$G_{\psi}$ is the stabilizer of $\psi$ under the conjugation action.\\
Its orbit is contained in $H$ since $H$ is normal in $G$.\\
So size of the orbit can be at most $p^{4}$.\\
Hence, $|G_{\psi}| \geq \frac {|G|} {p^{4}}$.\hfill $(3.3.5)$ \\~\\
Now $G$ acts transitively on the places of $K(N)$ which extend $\mathfrak{p}$.\\
The number of such places is \[\frac {[K(N) : K]} {d_{v, \mathfrak{p}}}\]  
where $d_{v, \mathfrak{p}}$ is the local degree of $v$ over $\mathfrak{p}$.\\
We have $ d_{v, \mathfrak{p}} = \frac{d_{v}} {f} $ .\\
So the orbit $G_{\psi}v$ of $v$ under action of $G_{\psi}$ has cardinality \[ |G_{\psi}v| \geq \frac {1} { [ G : G_{\psi} ]} \frac {[K(N) : K]f} {d_{v}}  = \frac { |G_{\psi}|f} {d_{v}} \geq p^{- 4}\frac {|G|f} {d_{v}} \]  using $(3.3.5)$ and $ |G| = [K(N) : K] $ . \\
This proves the lemma. \hfill $\square$ \\~\\
\textbf{Lemma 3.3.3 :} Assume that $ p | N $ and let $n$ be the positive integer with $ p^{n} | N$.\\
If $ \alpha \in K(N) $ satisfies $ \alpha^{Q(n)} \notin \mathbb{Q}_{w}(N/p) $, then there is a $\beta \in \overline{\mathbb{Q}} - \mu_{\infty} $ with $h(\beta) \leq 2p^{4}h(\alpha)$ and \[ h(\alpha) + \max \{ 0 , \frac {1} {[\mathbb{Q}(\beta) : \mathbb{Q}]} \sum_{\tau} \log | \tau(\beta) - 1|  \} \geq \frac { f\log p} {2dp^{8} } \] where $Q(n)$ is as in lemma 2.4.4, $|\cdot|$ denotes the usual absolute value on $\mathbb{C}$ and the sum is taken over all the field embeddings $\tau : \mathbb{Q}(\beta) \to \mathbb{C}$.\\~\\
\textbf {Proof :} By hypothesis we may choose $\psi \in \text{Gal} (\mathbb{Q}_{w}(N) | \mathbb{Q}_{w}(N/p)) $ with $\psi (\alpha^{Q(n)}) \neq \alpha^{Q(n)}$.\\
Clearly $\alpha \neq 0$.\\
Define $ x = \psi ( \alpha^{Q(n)}) - \alpha^{Q(n)} $. By assumption on $\psi$ we conclude that $ x \neq 0$.\\
So \[ \sum_{u} d_{u} \log |x|_{u} = 0 \] where the sum runs over the normalized absolute values of $K(N)$.\\
Let $G_{\psi}$ and $v$ be as in lemma 3.3.2.\\
Assume $\sigma \in G$ . Then the place $\sigma v$ of $K(N)$ satisfies $|\sigma (y)|_{\sigma v} = | y|_{v}$ for all $ y \in K(N)$.\\
Hence \[ 
\begin{split}
|\psi(\alpha)^{Q(n)} - \alpha^{Q(n)} |_{\sigma v} & = | \sigma\psi \sigma^{ - 1}(\alpha)^{Q(n)} - \alpha^{Q(n)} |_{\sigma v}\\
                                                                                  & = | \psi (\sigma^{- 1} (\alpha))^{Q(n)} - \sigma^{-1}(\alpha)^{Q(n)} |_{v} .
\end{split}
\]
By definition we have $ q \,|\, Q $ .
So applying lemma 3.2.2 to $\sigma^{- 1}(\alpha)^{Q/q}$ and using the computation above we conclude that 

\[\begin{split}
 | \psi (\alpha)^{Q(n)} - \alpha^{Q(n)} |_{\sigma v} & \leq p^{ - 1} \max \{ 1, |\psi(\sigma^{- 1}(\alpha))|_{v} \}^{Q(n)} \max \{ 1, |\sigma^{ - 1}(\alpha)|_{v} \}^{Q(n)} \\
                                                                                            & = p^{ - 1} \max \{1, |\sigma\psi\sigma^{ - 1}(\alpha)|_{\sigma v} \}^{Q(n)} \max \{ 1, |\alpha|_{\sigma v} \}^{Q(n)} \\
                                                                                            & = p^{ - 1} \max \{1, |\psi(\alpha)|_{\sigma v} \}^{Q(n)} \max \{ 1, |\alpha|_{\sigma v} \}^{Q(n)}.
\end{split}\]
Hence for all $u \in G_{\psi}v$, $|x|_{u} \leq p^{ - 1} \max \{ 1, |\psi(\alpha)|_{u} \}^{Q(n)} \{ 1, |\alpha|_{u} \}^{Q(n)} $ holds.\\
If $u$ is an arbitrary finite place of $K(N)$, then \[ |x|_{u} \leq \max \{ 1, |\psi(\alpha)|_{u} \}^{Q(n)} \max \{ 1, |\alpha|_{u} \}^{Q(n)} .\]
Now define $\beta = \frac {\psi(\alpha)^{Q(n)}} {\alpha^{Q(n)}} \in \overline{\mathbb{Q}} - \{1\} $.\\
Then we have the bound \[ |x|_{u} = | \beta - 1 |_{u} |\alpha|_{u}^{Q(n)} \leq | \beta - 1|_u \max \{ 1, |\alpha|_{u}\}^{Q(n)} .\]
Now by the product formula gives \[ 0 = \sum_{u} d_{u}\log |x|_u = \sum_{ u \in G_{\psi}v } d_{u}\log |x|_{u} + \sum_{ u\,|\, \infty} d_{u}\log |x|_{u} + \sum_{u \nmid \infty, u \notin G_{\psi}v} d_{u}\log |x|_{u} .\]
Putting the estimate derived above one obtains \[ 0 \leq \sum_{ u \in G_{\psi} v } d_{u}(\log p^{-1}) + \sum_{ u | \infty } d_{u}\log |\beta - 1|_{u} + Q(n)\sum_{u} d_u \log^{+}(|\psi(\alpha)|_u) + Q(n)\sum_{u}d_{u}\log^{+}(|\alpha|_{u}). \]
First note that $d_u = d_{v} \, \forall u \in G_{\psi}v$. So, $\sum_{u \in G_{\psi}v}d_{u}\log (p^{-1}) \leq - \frac {f|G|\log p} {p^{4}}$.\\
Thus \[ \frac {f|G|\log p} {p^4} \leq \sum_{u | \infty} d_{u}\log |\beta - 1|_{u} + Q(n)\sum_{u \nmid \infty, u \notin Gv}d_{u} \log^{+}( |\psi(\alpha)|_{u}) + Q(n)\sum_{u \nmid \infty, u \notin Gv}d_{u}\log^{+}(|\alpha|_{u}).\]
Dividing both sides by $[K(N) : \mathbb{Q}]$ we get \[ \frac {f\log p} {dp^{4}} \leq \frac {1} { [K(N) : \mathbb{Q} ]} \sum_{ u | \infty} d_{u} \log (|\beta - 1|_{u}) +  Q(n)h(\psi(\alpha)) + Q(n)h(\alpha) .\]
Now \[\frac {1} {[K(N) : \mathbb{Q}]} \sum_ {u|\infty} d_{u} \log(|\beta - 1|_u) = \frac {1} {[\mathbb{Q}(\beta) : \mathbb{Q}]} \sum_{\tau} \log |\tau(\beta) - 1| \] where in the right hand side the sum runs over all the embeddings of $\mathbb{Q}(\beta) $ in $\mathbb{C}$ and the absolute value is the usual absolute value and $h(\psi(\alpha)) = h(\alpha)$\\
Hence \[ 2Q(n)h(\alpha) + \max \{ 0, \frac {1} {[\mathbb{Q}(\beta) : \mathbb{Q}]} \sum_{\tau} \log (|\tau(\beta) - 1|) \} \geq \frac {f\log p} {dp^4}.\]
Note that $ Q(n) \leq p^4$.\\
So we have \[ h(\alpha) + \max \{0, \frac {1} {[\mathbb{Q}(\beta) : \mathbb{Q}]} \sum_{\tau} \log (|\tau(\beta) - 1|) \} \geq \frac {f\log p} {2dp^8}. \] 
For our choice of $\beta$ we have $h(\beta) \leq Q(n) (h(\psi(\alpha)) + h(\alpha)) \leq 2p^{4}h(\alpha)$.\\
If $\beta$ is a root of unity then so is $\frac {\psi(\alpha)} {\alpha}$. But  then by choice of $Q(n)$, $\beta = 1$ which is a contradiction.\\
Thus $\beta$ can not be a root of unity.\\
This proves the lemma. \hfill $\square$
\section{A descent argument }
\subsection{Matrices over finite fields }
Let $p \, \geq 5$ be a prime number and let $\mathbb{F}_{n}$ be the field with $n$ elements for any prime power $n$.\\
Let $M_{2}(\mathbb{F}_{p})$ denote the algebra of $2 \times 2$ matrices with entries in $\mathbb{F}_{p}$ and $\text{Gl}_{2}(\mathbb{F}_{p})$ be the group of invertible matrices in $M_{2}(\mathbb{F}_{p})$.\\
A nonsplit Cartan subgroup is a subgroup of $\text{Gl}_{2}(\mathbb{F}_{p})$ which is a cyclic subgroup of order $ q - 1$ where $q = p^{2}$.\\
Then one has the following lemma :\\~\\
\textbf{Lemma 3.4.1 :} i) Let $G$ be a nonsplit Cartan subgroup of $\text{Gl}_{2}(\mathbb{F}_{p})$. Then there is a $\mathbb{F}_{p}$ subalgebra of $M_{2}(\mathbb{F}_{p})$ which is isomorphic to $\mathbb{F}_{q}$ and its multiplicative group is $G$.\\
ii) Let $G$ be a nonsplit Cartan subgroup of $\text{Gl}_2(\mathbb{F}_{p})$. Then the set $\{hgh^{-1} \,| \, g \in G, h \in \text{Gl}_{2}(\mathbb{F}_{p})\}$  has cardinality strictly greater than $p^{3}$ and it generates $\text{Gl}_{2}(\mathbb{F}_{p})$.\\~\\
\textbf{Proof :} See \cite{hab} lemma 6.1. \hfill $\square$ \\~\\
\subsection{Descending along $p^{n}$ torsion}
Let $K$, $E$, $\mathfrak{p}$ and $p$ be as before.\\
Let $N \in \mathbb{N}$ with $ N = p^{n}M$ where $ n \geq 0$ and $ M \geq 1$ are integers and $ p \, \nmid M$.\\
Following our previous convention we consider $\text{Gal} \,(\mathbb{Q}_{w}(N) | \mathbb{Q}_{w})$ as a subgroup of $\text{Gal}\,(K(N)|K)$.\\ 
Then we have the following lemma :\\~\\
\textbf{Lemma 3.4.2: } Assume $n \geq 1$.\\
i) The subgroup of $\text{Gal}(K(N) | K)$ generated by conjugates of $\text{Gal} (\mathbb{Q}_{w}(N)|\mathbb{Q}_{w}(N/p))$ equals $\text{Gal}\,(K(N)|K(N/p))$.\\
ii) If $\alpha \in K(N)$ with $\sigma(\alpha) \in \mathbb{Q}_{w}(N/p)$ for all $\sigma \in \text{Gal} (K(N) | K)$, then $\alpha \in K(N/p)$.\\~\\ 
\textbf{Proof :} At first we put $ G = \text{Gal} (K(N) | K)$, $ H = \text{Gal} (K(N) | K(N/p))$ and $H_{\mathfrak{p}} = \text{Gal} (\mathbb{Q}_{w}(N) | \mathbb{Q}_{w}(N/p))$.\\
Let $c(H_{\mathfrak{p}})$ denote the normal closure of $H_{\mathfrak{p}}$ in $G$.\\
Since $H$ is normal in $G$ we conclude $c(H_{\mathfrak{p}}) \subseteq H$.\\
To prove $(i)$ one would like to show equality.\\~\\
Let $ \text{res} : \text{Gal} (K(N)|K) \to \text{Gal} (K(p^{n}) | K)$ be the restriction map.\\~\\
\textbf{Remark 3.4.3 :} i) res is onto.\\
ii) Let $S$ be any subgroup of of $\text{Gal}\,(K(N)|K)$. Then
\[\text{res}\, (\text{normal closure of S in Gal}\, (K(N)|K)) = \text{normal closure of res(S) in Gal}\,(K(p^{n}|K)\] since res is onto.\\
iii) $\text{res}(H) \subseteq \text{Gal}(K(p^{n}) | K(p^{n - 1}))$ .\\
iv) res is injective on $H$.\\~\\
Choose a pair of generators $(P_{\infty}, Q_{\infty})$ of $E[p^{\infty}]$ .\\
Let $(P_{m}, Q_{m})$ be the image of $(P_{\infty}, Q_{\infty})$ in $E[p^{m}]$ for each positive integer $m$.\\
If one fixes this choice of basis then $\text{Aut}_{\mathbb{Z}}E[p^{m}] \cong \text{Gl}_{2}(\mathbb{Z}/p^{m}\mathbb{Z})$.\\
Thus for each positive integer $m$ we have a representation of \[\rho_{m} : \text{Gal}\,(K(p^{m})|K) \to \text{Gl}_{2}(\mathbb{Z}/p^{m}\mathbb{Z}).\]
Clearly this map is injective since an element in the Galois group is determined by its action on $(P_{m}, Q_{m})$.\\
By our assumption iii on $\mathfrak{p}$, $\rho_{m}$ is onto for each positive integer $m$.\\
Now we embark on proof of part $(i)$.\\~\\
\textbf{Case I:} n = 1.\\
In this case the image of $H_{\mathfrak{p}}$ in $\text{Gl}_{2}(\mathbb{Z}/p\mathbb{Z})$ is a nonsplit Cartan subgroup.\\
By lemma 3.4.1 and the isomorphism mentioned above normal closure of $\text{res}(H_{\mathfrak{p}})$ in $\text{Gal} (K(p)|K)$ is nothing but whole of $\text{Gal} (K(p)|K)$.\\  
Thus by remark 3.4.3 ii and iv we conclude that $\text{res}\, (c(H_{\mathfrak{p}})) \cong \text{Gl}_{2}(\mathbb{Z}/p\mathbb{Z})$.\\
Again by remark 3.4.3 iv $|H| \, \leq |\text{Gal}\, (K(p)/K)| = |\text{Gl}_{2}(\mathbb{Z}/p\mathbb{Z})|$.\\
Thus we have the desired equality $H = c(H_{\mathfrak{p}})$.\\~\\
\textbf{Case II:} $n \geq 2$.\\
Let $\sigma \in \text{Gal} (K(p^{n})|K(p^{n-1}))$.\\
Consider $\rho_{n}(\sigma)$.\\
It is a matrix of the form 
\[\left(\text{Id} + p^{n - 1}
 \begin{bmatrix}
         a & b\\
         c & d
\end{bmatrix}\right)
\mod p^{n}\text{Mat}_{2}(\mathbb{Z})\]
where $a,b,c,d$ are integers in the interval $[0, p - 1]$.\\
So one has a map $l : \text{Gal}\, (K(p^{n})|K(p^{n - 1})) \to M_{2}(\mathbb{Z}/p\mathbb{Z})$ defined by
\[l(\sigma) =
\begin{bmatrix}
       a & b \\
       c & d \\
\end{bmatrix}
\] as above.\\~\\
\textbf{Claim 3.4.4 :} $l$ is a homomrphism of groups.\\~\\
\textbf{Proof :} Let $\sigma_{1}, \sigma_{2} \in \text{Gal}\, (K(p^{n})| K(p^{n - 1}))$.\\
Now 
\[\begin{split}
   \sigma_{1}\sigma_{2} & \equiv (1 + p^{n - 1}l(\sigma_{1}))(1 + p^{n - 1}l(\sigma_{2}))  \\                  
                                       & \equiv ( 1 + p^{n - 1}( l(\sigma_{1}) + l(\sigma_{2}) ))
\end{split}
\] where the equivalence is taken modulo $p^{n}$.\\
This proves $l(\sigma_{1}\sigma_{2}) = l(\sigma_{1}) + l(\sigma_{2})$.\\
Hence the claim. \hfill $\square$\\~\\
Further $l$ is injective since so is $\rho_{n}$.\\
$l$ is surjective since $\rho_{n}$ is surjective and each matrix in the image which is of the form 
\[1 + p^{n - 1}
\begin{bmatrix}
   a & b\\
   c & d
\end{bmatrix}\] actually arises from an element in $H$.\\
Define $L : \text{Gal}\, (K(N)| K(N/p)) \to M_{2}(\mathbb{Z}/p\mathbb{Z})$ by $L(\sigma) = l(\text{res}(\sigma))$. Here $L$ is also injective since so is res.\\
Further let $\text{res}_{1}: \text{Gal}\, (K(N)|K) \to \text{Gal}\,(K(p)|K)$ denote the usual restriction map. Define $L_{1} :\text{Gal}\,(K(N)|K) \to \text{Gl}_{2}(\mathbb{Z}/p\mathbb{Z})$ by $L_{1} = \rho_{1} \circ \text{res}_{1}$.\\
Let $\sigma \in G$ and $\psi \in H$. Then $\sigma\psi\sigma^{-1} \in H$.\\
A straight forward calculation shows $L(\sigma\psi\sigma^{-1}) = L_{1}(\sigma) L(\psi) L_{1}(\sigma)^{-1}$.\\
Now consider $H_{\mathfrak{p}}$.\\
By results in chapter 2 we know $|H_{\mathfrak{p}}| = p^{2}$.\\
Further by part iii of theorem 2.1.1 we know that $L(H_{\mathfrak{p}})$ contains all the scalar matrices in $M_{2}(\mathbb{Z}/p\mathbb{Z})$.\\
Since $|L(H_{\mathfrak{p}})| = p^{2}$ it contains at least one matrix which is not scalar.\\
Let $\theta$ be such a matrix.\\
Then by the theorem of Cayley- Hamilton we conclude that $\mathbb{F}_{p} + \mathbb{F}_{p}\theta$ contains $\theta^{2}$ where we are identifying $\mathbb{Z}/p\mathbb{Z}$ and $\mathbb{F}_{p}$.\\
Thus $\mathbb{F}_{p} + \mathbb{F}_{p}\theta$ is a subalgebra of $M_{2}(\mathbb{F}_{p})$.\\
Now we have\\~\\
\textbf{Claim 3.4.5 :} $\theta$ has no eigenvalues in $\mathbb{F}_{p}$.\\~\\
\textbf{Proof :} By theorem $2.1.1$ we know that $\text{Gal}\,(\mathbb{Q}_{w}(p^{n})|\mathbb{Q}_{w})$ is commutative.\\
Put $G_{1} = \text{Gal}\,(\mathbb{Q}_{w}(p) | \mathbb{Q}_{w})$ and think it as a subgroup of $G$.\\
Note that $L_{1}$ is injective on $G_{1}$ and thus its image is a cyclic subgroup of $\text{Gl}_{2}(\mathbb{F}_{p})$ of order $ q - 1$.\\
Note that the elements of  $L_{1}(G_{1})$ commute with the elements of $L(H_{\mathfrak{p}})$.\\
If necessary we can translate $\theta$ by a scalar matrix and assume that $\theta$ is invertible.\\
The statement above implies centralizer of $\theta$ in $\text{Gl}_{2}(\mathbb{F}_{p})$ has order dividing $ q - 1$.\\
If $\theta$ has one (and hence both) eigenvalue on the ground field then it must be conjugate to a matrix of the form 
\[
\begin{bmatrix}
         a & 0 \\
         0 & b
\end{bmatrix} 
\,\text{or} \,
\begin{bmatrix}
          a & 1 \\
          0 & a
\end{bmatrix}
\] where $a, b \in \mathbb{F}_{p} - \{0\}$. From our prior assumption on $\theta$ in the first case we have $a \neq b$.\\
If $C_{\theta}$ is the centralizer of $\theta$ and $\theta$ is conjugate to a matrix of the first form then $C_{\theta}$ is conjugate to the subgroup 
\[
\left\lbrace
\begin{bmatrix}
         r & 0 \\
         0 & s
\end{bmatrix}
\Bigg| \, r, s \, \in \mathbb{F}_{p} -\{0\}\right\rbrace.\]
In the second case $C_{\theta}$ is conjugate to a subgroup 
\[
\left\lbrace
\begin{bmatrix}
        r & \frac {as - r} {a - 1}\\
        0 & s
\end{bmatrix}
\Bigg|\, r, s \, \in \mathbb{F}_{p} - \{0\}
\right\rbrace
\,\text{or} \,
\left\lbrace
\begin{bmatrix}
       r & s\\
       0 & r
\end{bmatrix}
\Bigg| \, r \,\in \mathbb{F}_{p} -\{0\}, s \, \in \mathbb{F}_{p}
\right\rbrace
\] where the first case occurs if $a \neq 1$ and the second case occurs if $a = 1$.\\
None of these subgroups have order dividing $ q - 1$.\\
This contradiction proves the claim.\hfill $\square$\\~\\
Using the claim one concludes that the minimal polynomial of $\theta$ (one as in the proof of claim) is irreducible over $\mathbb{F}_{p}$.\\
So $\mathbb{F}_{p} + \mathbb{F}_{p}\theta$ is a field.\\
Thus $L(H_{\mathfrak{p}})$ contains a nonsplit Cartan subgroup.\\
Hence $L_{1}(G)L(H_{\mathfrak{p}})L_{1}(G)^{- 1}$ has cardinality strictly greater than $p^{3}$. But $| M_{2}(\mathbb{F}_{p})| = p^{4}$. \\
Thus the subgroup generated by it must be whole of $M_{2}(\mathbb{F}_{p})$.\\
Since $L$ is injective we conclude that $|c(H_{\mathfrak{p}})| = p^{4}$.\\
But by same reason $|H|\, \leq p^{4}$.\\
Thus we must have $c(H_{\mathfrak{p}}) = H$.\\
So in this case we are done.\\
This finishes the proof of part (i) in the lemma.\\~\\
Now put $L = K(N/p)(\alpha)$.\\
Clearly $\text{Gal}\, (K(N)|L) \subseteq \text{Gal}\,(K(N)|K(N/p))$.\\
By our assumption on $\alpha$ we have $\sigma\psi\sigma^{- 1} \in \text{Gal}\,(K(N) |L)$ for all $\sigma \in G$ and for all $\psi \, \in H_{\mathfrak{p}}$.\\
Now by part (i) of the lemma $\text{Gal}\,(K(N)|L) = \text{Gal}\,(K(N)|K(N/p))$.\\
So $\alpha\, \in K(N/p)$ as desired.\hfill $\square$\\~\\
\textbf{Lemma 3.4.6:} $K$, $E$, $\mathfrak{p}$, $p$ be as before.\\
We assume that $p^{2} \nmid N$.\\
If $\alpha \, \in K(N) - \mu_{\infty}$  is a nonzero algebraic number then there is a nonzero $\beta \in \overline{\mathbb{Q}} - \mu_{\infty}$ such that $h(\beta) \, \leq 2p^{4}h(\alpha)$ and
\[h(\alpha) + \max \big\{ 0, \frac {1} {[\mathbb{Q}(\beta) : \mathbb{Q}]} \sum_{\tau} \log |\tau(\beta) - 1| \big\} \geq \min \, \big\{\frac {f\log p} {2dp^{8}}, \frac{\log (p^f/2^{d})} {dq(q - 1)(w + 1)} \big\}, \]
where the sum runs over all the field embeddings $\tau : \mathbb{Q}(\beta) \to \mathbb{C} $ and $|\cdot|$ denotes the usual absolute value.\\~\\
\textbf{Proof :} Replacing $N$ by $pN$ if necessary, without loss of generality one can assume that $ p | N$.\\
If there is $\sigma \in \text{Gal}\, (K(N)|K)$ such that $\sigma(\alpha)^{Q(1)} \notin \mathbb{Q}_{q}(N/p)$ then one can apply lemma 3.3.3 to $\sigma(\alpha)$ to conclude the present lemma since $h(\sigma(\alpha)) = h(\alpha)$ .\\
If $\sigma(\alpha)^{Q(1)} \in \mathbb{Q}_{w}(N/p)$ for all $\sigma \in \text{Gal}\, (K(N)|K)$ then by previous lemma one concludes that $\alpha^{Q(1)} \in K(N/p)$.\\
But by our assumption on $N$ we have $\gcd (\frac {N} {p}, p) = 1$.\\
So using lemma 3.3.1 we conclude that $h(\alpha^{Q(1)}) \, \geq \frac {\log (p^f/2^{d})} {d(w + 1)}$.\\
Thus we have $h(\alpha) \, \geq \frac {\log (p^f/2^{d})} {dq(q - 1)(w + 1)}$.\\ 
Now the lemma follows if one simply chooses $\beta = \alpha$. \hfill $\square$\\~\\
\textbf{Lemma 3.4.7: } $K$, $E$, $\mathfrak{p}$, $p$ be as before.\\
Let $N$ be a positive integer. Put $v_{p}(N) = n$.\\
Then there is $\sigma \, \in \text{Gal} (\mathbb{Q}_{w}(N)| \mathbb{Q}_{w})$ satisfying the following properties:\\
i) $\sigma$  lies in the center of $\text{Gal} (K(N)|K)$.\\
ii) $\sigma$ acts on $E[p^{n}]$ as multiplication by $2$.\\
iii) $\sigma(\zeta) = \zeta^{4}$ for all $\zeta \in \mu_{p^{n}}$.\\~\\
\textbf{Proof :} Using theorem 2.1.1 and proposition 2.4.2 we conclude there is a $\sigma \, \in \text{Gal}\,(\mathbb{Q}_{w}(N) | \mathbb{Q}_{w})$ such that $\sigma$ acts on $E[p^{n}]$ as multiplication by $2$ and on $E[N/p^{n}]$ it acts like identity.\\
Hence it lies in the center of $\text{Gal}\,(K(N)|K)$ since $E[N] = E[p^{n}] \oplus E[N/p^{n}]$ .\\
By properties of Weil pairing it follows that this $\sigma$ also satisfies property (iii).
Hence the lemma. \hfill $\square$\\~\\
\textbf{Lemma 3.4.8: } There are positive constants $C_{1}, C_{2}$ depending only on $d, f, \text{and}\, p$ such that the following holds :\\
If $\alpha \in K(E_{\text{tor}}) - \mu_{\infty}$ is nonzero then there exists a nonzero $\beta$ in $\overline{\mathbb{Q}} - \mu_{\infty}$ with $h(\beta) \,\leq\, C_{1} h(\alpha)$ and 
\[ h(\alpha) + \frac {1} {5}\max \big\{ 0, \frac {1} {[\mathbb{Q}(\beta) : \mathbb{Q}]} \sum_{\tau} \log |\tau(\beta) - 1| \big\} \geq C_{2}  \] where the sum runs through all the embeddings of $\mathbb{Q}(\beta)$ in $\mathbb{C}$ and $|\cdot|$ is the usual complex absolute value.\\~\\
\textbf{Proof :} Let $\alpha$ be as in the statement.\\
Let $N$ be a positive integer such that $\alpha \in K(N) - \mu_{\infty}$.\\
Say $N= p^{n}M$ where $n$ is a nonnegative integer and $M$ is a positive integer with $\gcd \,(p, M) = 1$.\\
Let $\sigma_{0} \, \in \text{Gal}\,(\mathbb{Q}_{w}(N)|\mathbb{Q}_{w})$ be the element constructed in lemma 3.4.7.\\
Put \[ \gamma = \frac {\sigma_{0}(\alpha)} {\alpha^{4}} \in K(N).\tag{3.4.1}\]
First we we note that $\gamma \, \notin \mu_{\infty}$. Otherwise we shall have 
\[ h(\alpha) = h(\sigma_{0}(\alpha)) = h(\gamma\alpha^{4}) = h(\alpha^{4}) = 4h(\alpha) \]
using the basic height properties. Thus $h(\alpha) = 0$. But then $\alpha$ is either 0 or a root of unity contrary to our assumption.\\
Now \[h(\gamma) \,\leq h(\sigma_{0}(\alpha)) + h(\alpha^4) = 5h(\alpha). \tag{3.4.2} \]
Let $n_{1}\, \geq 0$ be the least integer such that $\sigma(\gamma) \,\in \mathbb{Q}_{w}(p^{n_{1}}M)$ for all $\sigma \in \text{Gal}\,(K(N)|K)$. It is easy to see such a least element exists and it is $\leq \, n$. By lemma 3.4.2 we conclude that $\gamma \in K(p^{n_{1}}M)$. Now we consider two cases :\\~\\
\textbf{Case I:} $n_{1} \leq 1$.\\
Then using lemma 3.4.6 we conclude that there is a nonzero $\beta \in \overline{\mathbb{Q}} - \mu_{\infty}$ with $h(\beta) \leq 2p^{4}h(\gamma)$ and
\[ h(\gamma) + \max \big\{ 0, \frac {1} {[\mathbb{Q}(\beta) : \mathbb{Q}]} \sum_{\tau} \log \,|\tau(\beta) - 1| \big\} \, \geq c_{2} \]
where $c_{2} = \min \{\frac {f \log p} {2dp^{8}}, \frac {\log (p^f/2^{d})} {dq(q - 1)( w + 1)} \}$ .\\
Using (3.4.2) we conclude the lemma in this case with the choice of constants $C_{1} = 10p^{4}$ and $C_{2} = \frac {1} {5}c_{2}$.\\~\\
\textbf{Case II:} $n_{1} \geq 2$.\\
By minimality of $n_{1}$ there is a $\sigma \in \text{Gal}\,(K(N)|K)$ such that $\sigma(\gamma) \notin \mathbb{Q}_{w}(p^{n_{1} - 1}M)$.\\
Put $\alpha_{1} = \sigma(\alpha)$ and $\gamma_{1} = \sigma(\gamma)$. Applying $\sigma$ to both sides of $(3.4.1)$  and using the fact that $\sigma_{0}$ lies in the center we get \[ \gamma_{1} = \frac {\sigma_{0}(\alpha_{1})} {\alpha_{1}^{4}} .\tag{3.4.3}\]
Clearly $\gamma_{1} \notin \mu_{\infty}$ since $\gamma \notin \mu_{\infty}$.\\
Now we want to apply lemma 3.3.3 to $\gamma_{1}$. First we want to verify the hypothesis for $\gamma_{1}$.\\
Note that since $n_{2} \geq 2$ we have $Q(n_{2}) = q$.\\
We need to show that $\gamma_{1}^{q} \notin \mathbb{Q}_{w}(p^{n_{1} - 1}M)$.\\
Assume the contrary. Then there is a $\psi \in \text{Gal}(\mathbb{Q}_{w}(p^{n_{1}}M)|\mathbb{Q}_{w}(p^{n_{1} - 1}M))$ such that $\psi(\gamma_{1}) \neq \gamma_{1}$. But $\psi(\gamma_{1}^{q}) = \gamma_{1}^{q}$.\\
Thus there is a $\xi \in \mathbb{Q}_{w}(p^{n_{1}}M)$ such that $\psi(\gamma_{1}) = \xi \gamma_{1}$ such that $\xi^{q} = 1$ and $\xi \neq 1$.\\
We identify $\psi$ with its restriction to $K(N)$ and apply it to $(3.4.3)$ to get \[ \xi\gamma_{1} = \frac {\sigma_{0}(\psi(\alpha_{1}))} {\psi(\alpha_1)^{4}} \tag{3.4.4} \]
since $\sigma_0$ lies in the center.\\
Define $\eta = \frac {\psi(\alpha_{1})} {\alpha_{1}} \neq 0$. $(3.4.3)$ and $(3.4.4)$ implies  \[ \xi = \frac {\sigma_{0}(\eta)} {\eta^{4}}. \tag{3.4.5} \]
Since $\xi$ is a root of unity we have \[h(\eta) = h(\sigma_{0}(\eta)) = h(\xi \eta^{4}) = h(\eta^4) = 4h(\eta)\] using properties of height and thus $h(\eta) = 0$.\\
So $\eta = \frac {\psi(\alpha_1)} {\alpha_1} \in \mu_{\infty}$.\\
Fix a positive integer $m$ coprime to $p$ such that $\eta^{m} \in \mu_{p^{\infty}}$ . By lemma $2.4.3$ one concludes that $\eta^{m} \in \mu_{p^n}$. Using properties of $\sigma_{0}$ we have \[\sigma_{0}(\eta) = \xi_{1} \eta^{4}\tag{3.4.6}\] where $\xi_{1} \in \mu_{\infty}$ with $\xi_{1}^{m} = 1$.\\
Using $(3.4.5)$ and $(3.4.6)$ one concludes that $\xi = \xi_{1}$. But $\xi^{q} = \xi_{1}^{m} = 1$. Since $\gcd (q, m) = 1$ this identity gives $\xi = 1$ which is a contradiction.\\
This contradiction says $\gamma_{1}^{q} \notin \mathbb{Q}_{w}(p^{n_{1} - 1}M)$.\\
So we can apply lemma 3.3.3 to $\gamma_{1}$ . Note that $h(\gamma_{1}) = h(\gamma) \leq 5h(\alpha)$.\\
So using lemma 3.3.3 we conclude the present lemma in this case with the choice of constants $C_{1} = 10p^{4}$ and $C_{2} = \frac {1} {5} c_2$.\\~\\
Thus one has the lemma and the choice of constants $C_{1} =10p^{4}$ and $C_{2} = \frac {1} {5} \min \{\frac {f \log p} {2dp^8}, \frac {\log (p^f/2^d)} {dq(q -1)(w+ 1)}\}$ works in all cases. \hfill $\square$  \\
\section{The final estimate}
The main goal of this section is to handle the sum involving the infinite places occurring in lemma 3.8.4 to get a positive lower bound for $h(\alpha)$.\\
One can use an equidistribution theorem of Bilu to do this as done in the original paper of Habegger. But this can be done in more elementary way as shown by Frey in a subsequent work. We shall use one of her results: \\~\\
\textbf{ Lemma 3.5.1: } Let $ 0 < \delta < \frac {1} {2}$ and let $\beta \in \overline{\mathbb{Q}} - \mu_{\infty}$ be such that $[\mathbb{Q}(\beta) : \mathbb{Q}] \geq 16$ and $h(\beta)^{\frac {1} {2}} \leq \frac {1} {2}$. Then 
\[ \frac {1} {[\mathbb{Q}(\beta) : \mathbb{Q}]} \sum_{\tau} \log |\tau(\beta) - 1| \leq \frac {4} {\delta^4} h(\beta)^{\frac {1} {2} - \delta} \]
where as before $\tau$ runs over all the field embeddings and $|\cdot|$ is the usual complex absolute value.\\~\\ 
\textbf{Proof :} See \cite{frey} lemma 3.5. \hfill $\square$ \\~\\
\textbf{Proof of theorem 0.1:} We want to show that there is a positive constant $C$  depending only on $d, f, p$  such that for all nonzero $\alpha \in K(E_{\text{tor}}) - \mu_{\infty}$ we have \[ h(\alpha) \geq C .\]
One wants to use lemma 3.4.8. Let $C_{1}$ and $C_{2}$ be as in the statement of that lemma.\\
Let $\alpha$ be a nonzero element of $ K(E_{\text{tor}}) - \mu_{\infty}$ .\\
Get a corresponding $\beta$ as provided by lemma 3.4.8.\\
Now we want to consider two cases:\\~\\
\textbf{Case I:} $[\mathbb{Q}(\beta):\mathbb{Q}] \leq 16$ or $h(\beta)^{1/2} \geq \frac{1}{2}$.\\
Note that $\beta$ is nonzero and it is not a root of unity.\\
In this situation one invokes the following result:\\
Let $\epsilon > 0$. Let $\beta \in \mathbb{Q} - \mu_{\infty}$ be a nonzero algebraic number of degree of $d$. Then there is a positive constant $c(\epsilon)$ depending only on $\epsilon$ such that \[ h(\beta) \geq \frac{c(\epsilon)}{d^{1 +\epsilon}}. \tag{3.5.1}\]
Such and even stronger results are well known ( See \cite{vou}).\\
Thus in any of the cases under consideration we conclude that there is a positive universal constant $c$ such that $h(\beta) \geq c$. Hence \[h(\alpha) \geq \frac{c} {C_1}.\].\\~\\
\textbf{Case II:} $[\mathbb{Q}(\beta):\mathbb{Q}] \geq 16$ and $h(\beta)^{1/2} \leq \frac{1}{2}$.\\
Then one can use lemma 3.5.1. Putting the estimate obtained from lemma 3.5.1 in the inequality of lemma 3.4.8 we obtain
\[ h(\alpha) + \frac {4} {\delta^{4}} h(\beta)^{\frac {1} {2} - \delta} \geq C_{2} .\tag{3.5.2} \]
Now  $h(\beta) \leq C_{1}h(\alpha)$. Choose $\delta = \frac {1} {4}$. Then we have 
\[ h(\alpha) + 4^5C_{1}^{1/4} h(\alpha)^{1/4} \geq C_{2}.\tag{3.5.3}\]
Consider the situation $h(\alpha) \leq 1$.\\
Then $h(\alpha) \, \leq h(\alpha)^{1/4}$ and thus using $(3.5.3)$ we have \[ h(\alpha) \, \geq \frac {C_2^{4}} {(1+ 4^{5}C_{1}^{1/4})^4}.\]
Hence in this case, the choice $C = \min \{ 1, \frac {C_{2}^{4}} {(1 + 4^{5}C_{1}^{1/4})^{4}} \}$ works.\\~\\
Thus theorem follows if one chooses $ C = \min \{1, \frac {c} {C_1}, \frac {C_2^4} {(1 + 4^{5}C_{1}^{1/4})^4}\} > 0.$\hfill $\square$ \\~\\
\textbf{Remark 3.5.2 :} Our computation does not attempt to find the optimal results. Some arguments can be made better. For example, the estimate in $(3.5.1)$ and the choice of $\delta$ can be improved.\\

\chapter{Diophantine estimates on elliptic curves}
\section{Introduction}
\textbf{Notation :} i) We shall use the notations $K$, $E$, $\mathfrak{p}$ and $p$  in the sense of chapter 3 throughout this chapter except in section 2.\\
ii) The group law on $E$ shall be denoted by same $\pm$. Meaning will be clear from context. $O$ shall denote the identity element.\\~\\
In this chapter we shall give a proof of theorem $0.2$ .\\
Before going into the main computations we review the basic facts about the N\'eron - Tate height on elliptic curve and its decomposition in terms of local height.\\
\section{N\'eron - Tate height on elliptic curves}
\textbf{Notation:} In this section the notation $E$ will be used to denote a general elliptic curve. $F$ shall denote an arbitrary number field.\\~\\
Let $E$ be an elliptic curve over $F$ presented in a given Weierstrass form .\\~\\
\textbf{Proposition 4.2.1 :} Let $P \in E(\overline{\mathbb{Q}})$. Define 
\[ h(P) =
\begin{cases}
0 & \text{if} \, P = O, \\
\frac {h(x)} {2} & \text{if} \, P = ( x, y) \neq O. 
\end{cases}
\]
Then the limit \[\lim_{ n \to \infty}\frac {h([2^{n}]P)} {4^{n}} \]  exists.\\~\\
\textbf{Proof :} See \cite{sil1} chapter VIII, proposition 9.1. \hfill $\square$ \\~\\
\textbf{Definition 4.2.2 :} The N\'eron - Tate height on $E$ is the function $\widehat{h} : E(\overline{\mathbb{Q}}) \to \mathbb{R}$ defined by \[ \widehat{h}(P) = \lim_{n \to \infty} \frac {h([2^n]P)} {4^{n}}. \]\\
The following proposition lists the properties of N\'eron - Tate height which we shall use again and again.\\~\\
\textbf{Proposition 4.2.3 :} Let $E$ be an elliptic curve defined over a number field $F$ and let $\widehat{h}$ be the N\'eron - Tate height defined on $E$.\\
Then the following holds :\\
i) Let $P, Q \in E(\overline{\mathbb{Q}})$. Then \[ \widehat{h}(P + Q) + \widehat{h}(P - Q) = 2 (\widehat{h}(P) + \widehat{h}(Q) ).\]
ii) $\widehat{h}([m]P) = m^{2}\widehat{h}(P)$ for all $ m \in \mathbb{Z}$ and for all $P \in E(\overline{\mathbb{Q}})$.\\
iii) $\widehat{h}$ is a quadratic form. In other words, $\widehat{h}$ is even and the pairing \[<,> : E(\overline{\mathbb{Q}})\times E(\overline{\mathbb{Q}}) \to \mathbb{R}\] given by $< P, Q > = \widehat{h}(P+Q) - \widehat{h}(P) - \widehat{h}(Q)$ is bilinear.\\
iv) Let $P \in E(\overline{\mathbb{Q}})$. Then $\widehat{h}(P) \geq 0$. Further $\widehat{h}(P) = 0$ if and only if $P$ is a torsion point.\\
v) Let $Q$ be a torsion point. Then $\widehat{h}(P + Q) = \widehat{h}(P)$ for all $P \in E(\overline{\mathbb{Q}})$.\\~\\
\textbf{Proof :} See \cite{sil1} chapter VIII theorem 9.3.\hfill $\square$\\~\\
\subsection{Local Height Functions}
In this subsection we briefly recall the facts about local height functions. For details one is referred to \cite{sil 2}, chapter VI.\\~\\
\textbf{Proposition 4.2.4 :} Let $K$ be a field which is complete with respect to an absolute value $|\cdot|_{v}$ and let \[v(\cdot) = - \log |\cdot|_{v} \] denote the corresponding additive absolute value. Let $E/K$ be an elliptic curve. Choose a Weierstrass equation for $E/K$,
\[ E : y^{2} + a_{1}xy + a_{3}y = x^{3} + a_{2}x^{2} + a_{4}x + a_{6}, \] and let $\Delta$ be the discriminant of this equation.\\
a) There exists a unique function \[ \lambda : E(K) - \{O\} \to \mathbb{R} \] with the following properties:\\
i) $\lambda$ is continuous on $E(K) - \{O\}$ with respect to the $v$-adic topology and is bounded on the complement of any $v$-adic neighborhood of $O$.\\
ii) The limit \[ \lim_{ P \to O} \{\lambda(P) + \frac {1} {2}v(x(P)) \} \] exists where the limit is being taken with respect to $v$-adic topology.\\
iii) For all $P \in E(K)$ with $[2]P \neq O$, \[\lambda([2]P) = 4 \lambda(P) + v((2y + a_{1}x + a_{3})P) - \frac {1} {4}v(\Delta). \]     
b)$\lambda$ is independent of the choice of Weierstrass equation for $E/K$.\\
c) Let $L/K$ be a finite extension and $\overline{v}$ is the extension of $v$ to $L$. Then \[ \lambda_{\overline{v}}(P) = \lambda_{v}(P) \hspace{1cm} \text{for all } P \in E(K) - \{O\}.\]\\
\textbf{Proof :} See \cite{sil 2} Chapter VI, Theorem 1.1. \hfill $\square$  \\~\\
\textbf{Remark 4.2.5 :} Let $\overline{K}$ denote a fixed algebraic closure of $K$ and let $\overline{v}$ denote a fixed extension of $v$ to $\overline{K}$.\\ 
Then using part c of proposition 4.2.4 one can $\lambda$ extend uniquely on $E(\overline{K}) - \{O\}$.\\~\\
\textbf{Proposition 4.2.6 :} Let $F$ be a number field. Let $M_{F}$ denote the set of normalized absolute values on $F$. Put $d_{v} = [F_{v} : \mathbb{Q}_{v}]$ for each $v \in M_{F}$. Let $E/F$ be an elliptic curve. For each $v \in M_{F}$, let $\lambda_{v}: E(F_{v}) - \{O\} \to \mathbb{R}$ be the local height function as described in proposition 4.2.4. Then \[ \widehat{h}(P) = \frac {1} {[F: \mathbb{Q}]} \sum_{v \in M_{F}} d_{v}\lambda_{v}(P) \] for all $P \in E(F) - \{O\}$.\\~\\
\textbf{Proof :} See \cite{sil 2} chapter VI, theorem 2.1. \hfill $\square$ \\~\\
\textbf{Lemma 4.2.7 :} Let $E$ and $F$ be as in statement of proposition 4.2.6. Let $v$ be a finite place of $E$. Assume that $E$ has a good reduction at $v$. Fix a minimal Weierstrass equation for $E$ at $v$. Then $\lambda_{v}(P) = \frac {1} {2} \max \{ - v(x(P)), 0\}$ where $x(P)$ denotes the x- coordinate of $P$ with respect to the Weierstrass presentation we have fixed. \\~\\
\textbf{Proof :} See \cite{sil 2} chapter VI, theorem 4.1. \hfill $\square$ 
\section{First height estimates} 
\textbf{Notation :} i) Let $l$ be a prime. We use the notation $E[l^{\infty}]$ to denote all points on $E$ (defined over $\overline{\mathbb{Q}}$) which are annihilated by a power of $l$.\\
ii) $N$ shall denote a positive integer. We shall write $N = p^n M$ where $n$ is an integer $\geq 0$ and $M$ is a positive integer such that $p \nmid M$.\\
iii) We fix an embedding of $K$ in $\overline{\mathbb{Q}}_p$ by completion at $\mathfrak{p}$. Thus $E$ can be thought as an elliptic curve over $\mathbb{Q}_{w}$ as before.\\
iv) $|\cdot|_p$ on $\mathbb{Q}_{w}(N)$ determines a place of $K(N)$ for each $N$. We shall call this place $v$.\\
v)Let $\widetilde{E}$, $\widetilde{\phi}$ be as in chapter 2. $\widetilde{\phi}$ acts on $\widetilde{E}$ as multiplication by an integer. We shall denote this integer by $a_w$.\\
vi) For an extension of local fields $L/K$, $e(L/K)$ shall denote the ramification index of the extension.\\
vii) Let $L$ be a finite extension of $K$. Let $u$ be a fixed place of $K$. For any place $\overline{u}$ of $L$ which extends $u$ we put $d_{\overline{u},u} = \frac{d_{\overline{u}}}{d_{u}}$ and define \[ \widehat{h}_{u}(P) = \frac {1} {[L : K]} \sum_{\overline{u} | u} d_{\overline{u}, u} \lambda_{\overline{u}}(P) \hspace{1 cm} \text{for all}\, P \in E(L). \] With this notation one has \[ \widehat{h}(P) = \frac {1} {[K : \mathbb{Q}]}\sum_{u \in M_{K}} d_{u} \widehat{h}_{u}(P) \hspace{1 cm} P \in E(L).\] \\~\\
The goal of this section is to prove a series of lemmas which will eventually lead to the proof of theorem 0.2.\\
These lemmas are similar to the results proved in chapter 3 and in most cases are elliptic analogues of them.\\
We assume that $E$ is presented in a Weierstrass form over $K$ such that it is minimal thought as a curve in $\mathbb{Q}_w$. By hypothesis, it has a good supersingular reduction. \\
\subsection{Lemmas on torsion points}
The next two lemmas are analogues of lemma 2.4.3 and 2.4.4 respectively :\\~\\
\textbf{Lemma 4.3.1 :} $E(\mathbb{Q}_{w}(N)) \cap E[p^{\infty}] = E[p^n].$\\~\\
\textbf{Proof :} Clearly $E[p^{n}] \subseteq E(\mathbb{Q}_{w}(N)) \cap E[p^{\infty}]$. We would like to show the reverse inclusion.\\
Let $T \in E(\mathbb{Q}_{w}(N)) \cap E[p^{\infty}]$ . Assume that $T$ has order $p^{m}$ for some integer $m \geq 0$. We would like to show $m \leq n$.\\
Consider the case $n = 0$. Then $\mathbb{Q}_{w}(N)/\mathbb{Q}_{w}$ is an unramified extension. Hence by lemma 2.3.3 we conclude that $m\geq1$ cannot hold.\\
So $m =0 = n$ and we are done in this case.\\
Now let $n \geq 1$. If already $m \leq n$, then we are done. Assume $ m > n$.\\
Using proposition 2.4.2 and theorem 2.1.1 we conclude that $e(\mathbb{Q}_{w}(N)|\mathbb{Q}_{w}) = q^{n - 1}(q - 1)$.\\
But by lemma 2.3.3 $e(\mathbb{Q}_w(T)|\mathbb{Q}_w) = q^{m - 1}(q - 1)$.\\
Since $\mathbb{Q}_w(T) \subseteq \mathbb{Q}_{w}(N)$ we conclude that $m \leq n$, a contradiction !
This contradiction proves the lemma. \hfill $\square$\\~\\
\textbf{Lemma 4.3.2 :} Assume that $n \geq 1$ and $\psi \in \text{Gal}(\mathbb{Q}_{w}(N)|\mathbb{Q}_{w}(N/p))$. Let $A \in E(\mathbb{Q}_{w}(N))$ such that $\psi(A) - A \in E_{\text{tor}}$. Then,
\[\psi(A) - A \in E[Q(n)]\] where $Q(n)$ is as in lemma 2.4.4.\\~\\
\textbf{Proof :} Put $B = \psi(A) - A$ and assume that order of $B$ is $N' = p^{n'}M'$ where $n' \geq 0$ and $p\nmid M'$.\\
Let $T= [p^{n'}](B)$ and note that its order is $M'$ which is coprime to $p$. Arguing exactly like lemma 2.4.3 (here we need to use proposition 2.4.2 instead of theory of cyclotomic extensions) we conclude that $T \in E[M]$, in particular $\psi(T) = T$.\\
The order of $[M'](B)$ is $p^{n'}$. So by lemma 4.3.2 $[M'](B) \in E[p^n]$. Then $[pM'](B) \in E[p^{n - 1}]$.\\
Therefore $\psi([pM'](B)) = [pM'](B)$.\\
Now one can imitate the proof of lemma 2.4.4 to get $[pt](B) = O$ where $t$ is order of $\psi \in \text{Gal}\,(\mathbb{Q}_{w}(N)|\mathbb{Q}_{w}(N/p))$.\\
From here the lemma follows easily. \hfill $\square$\\~\\

\subsection{Local estimates}
In this subsection we shall consider $E$ and $\widetilde{E}$ as projective curves. So without loss of generality, the points can be assumed to have primitive integral coordinates and we have a well defined reduction map.\\~\\ 
\textbf{Lemma 4.3.3 :} Say $p \nmid N$ and $A \in E(\mathbb{Q}_{w}(N))$.\\
Then $A \in E(\mathbb{Q}_{p}^{ur})$. If additionally we have $\phi_{w}(A) \neq [a_{w}]A$ then, \[\lambda_{v}(\phi_{w}(A) - [a_w](A)) \geq \frac {1} {2}\log p .\]\\~\\
\textbf{Proof :} The first part of the assertion follows from proposition 2.4.1.\\
Note that $v$ extends the place $\mathfrak{p}$. Thus it is a place of good reduction.\\
Since $\widetilde{\phi} = [a_w]$ on $\widetilde{E}$, $\phi_{w}(A) - [a_w](A) \equiv 0 \mod p$.\\
The extension $\mathbb{Q}_{w}(N)/\mathbb{Q}_{w}$ is unramified.
Now the lemma follows from lemma 4.2.7. \hfill $\square$ \\~\\
\textbf{Lemma 4.3.4 :} Let $n \geq 1$ and $\psi \in \text{Gal}\,(\mathbb{Q}_{w}(N)/\mathbb{Q}_{w}(N/p)\big)$.\\
Assume that $A \in E(\mathbb{Q}_{w}(N))$ with $\psi(A) \neq A$. Then \[ \lambda_{v}(\psi(A) - A) \geq \frac {\log p} {2(q - 1)}.\]\\
\textbf{Proof :} As in proof of lemma 3.2.2 $\psi \in G_{i}(\mathbb{Q}_{w}(N)|\mathbb{Q}_{w})$ for $ i = q^{n - 1} - 1$.\\
Let $\mathfrak{P}$ the unique maximal ideal of the ring of integers of $\mathbb{Q}_{w}(N)$.\\
So $\psi(A)$ and $A$ are same element of $E$ reduced modulo $\mathfrak{P}^{q^{n - 1}}$. Let $x$ denote the first coordinate of $\psi(A) - A$. Then $\log |x|_p  \leq - \frac {q^{n - 1}} {e} \log p$ where $e = e(\mathbb{Q}_{w}(N)|\mathbb{Q}_{w}) = q^{n - 1}(q - 1)$ using theorem 2.1.1 and proposition 2.4.2.\\
Now the lemma follows from lemma 4.2.7. \hfill $\square$\\~\\

\subsection{Lower bound on $\widehat{h}_{\mathfrak{p}}$ }
\textbf{Lemma 4.3.5 :} Assume that $n = 0$. If $A \in E(K(N)) - E_{\text{tor}}$ there is a nontorsion point $B \in E(\overline{\mathbb{Q}})$ with $\widehat{h}(B) \leq 2( w+ 1)\widehat{h}(A)$ such that \[ \widehat{h}_{\mathfrak{p}}(B) \geq \frac {1} {2}\log p.\]\\
\textbf{Proof :} Since $p \nmid N$, $K(N)/K$ is unramified at $\mathfrak{p}$ and $\phi_w$ acts on $K(N)$. Further since completion of $K$ at $\mathfrak{p}$ is $\mathbb{Q}_{p^f} \subseteq \mathbb{Q}_{w}$ we conclude that $\phi_w$ is trivial on $K$. Thus $\phi_w \in \, \text{Gal}\,(K(N)|K)$. As in proof of lemma 3.3.1 it lies in the center of the Galois group.\\
Put $B = \phi_{w}(A) - [a_w]A$. First note that $B$ is not a torsion point. Indeed, if it were a torsion point then $\widehat{h}(A) = \widehat{h}(\phi_{w}(A)) = w \widehat{h}(A)$ (recall $a_w = \pm p^f$) which implies $\widehat{h}(A) = 0$ contradicting our hypothesis $A \notin E_{\text{tor}}$.\\
Using parallelogram identity 
\[
\begin{split}
\widehat{h}(B) & \leq \widehat{h}(\phi_{w}(A) - [a_w](A)) + \widehat{h}(\phi_{w}(A) + [a_w](A)) \\
                         & = 2 (\widehat{h}(\phi_w(A)) + \widehat{h}([a_w](A))) \\
                         & = 2 ( w + 1) \widehat{h}(A). \\
\end{split}
\] 
Note that $\text{Gal}\,(K(N)|K)$ acts transitively on the places lying over $\mathfrak{p}$. Since $\phi_w$ lies in the center of the Galois group we have \[\lambda_{\sigma^{ - 1}v}(B) = \lambda_{v}(\phi_{w}(\sigma(A)) - [a_w](\sigma(A)))\]
for all $\sigma \in \text{Gal}\,(K(N)|K)$.\\
Clearly $|\text{Gal}\,(K(N)|K)v| = \sum_{u| \mathfrak{p},\mathfrak{p}\in M_{K(N)}} d_{u,\mathfrak{p}} $.\\
Using the estimate given in lemma 4.3.3 we conclude that the present lemma follows for this choice of $B$. \hfill $\square$\\~\\
\textbf{Lemma 4.3.6 :} Assume that $n \geq 1$. If $A \in E(K(N))$ satisfies $[Q(n)](A) \notin E(\mathbb{Q}_{w}(N/p))$, there is a nontorsion point $B \in E(\overline{\mathbb{Q}})$ with $\widehat{h}(B) \leq 4\widehat{h}(A)$ and
\[\widehat{h}_{\mathfrak{p}}(B) \geq \frac {\log p} {2p^6}. \]\\
\textbf{Proof :} It follows from hypothesis that there is a $\psi \in \text{Gal}\,(\mathbb{Q}_{w}(N)/\mathbb{Q}_{w}(N/p))$ such that $\psi([Q(n)](A)) \neq [Q(n)](A)$. Fix such a $\psi$. We shall use the same notation for its restriction to $K(N)$.\\
Let $B = \psi(A) - A \in E(\overline{\mathbb{Q}})$. By choice of $\psi$ it follows that $[Q(n)](B) \neq O$.\\
Using lemma 4.3.2 we conclude that $B \notin E_{\text{tor}}$.\\
An application of parallelogram law as in the previous lemma gives \[\widehat{h}(B) \leq 2(\widehat{h}(\psi(A)) + \widehat{h}(A)) = 4\widehat{h}(A).\]
Let $G_{\psi}$ denote the centralizer of $\psi$ in the group $ G = \text{Gal}(K(N)|K)$.\\
By lemma 3.3.2 we have \[|G_{\psi}v| \geq p^{- 4} \frac {|G|} {d_{v, \mathfrak{p}}}.\tag{4.3.1}\] 
For any $\sigma \in G_{\psi}$ we have $\lambda_{\sigma^{- 1}v} (B) = \lambda_v(\sigma(B)) = \lambda_v(\psi(\sigma(A)) - \sigma(A))$ and $\psi(\sigma(A)) \neq \sigma(A)$.\\
Applying lemma 4.3.4 we conclude that \[\lambda_{\sigma^{- 1}v}(B) \geq \frac {\log p} {2( q - 1)}\tag{4.3.2}\] for all $\sigma \in G_{\psi}$.\\
Note that if $u \in G_{\psi}v$ then $d_{u, \mathfrak{p}} = d_{v, \mathfrak{p}}$.\\
$\mathfrak{p}$ is a place of good reduction. Hence for each $u | \mathfrak{p}$, $u$ is a place of good reduction and $\lambda_{u}(B) \geq 0$.\\
 
Now 
\[
\begin{split}
\widehat{h}_{\mathfrak{p}}(B) & = \frac {1} {[K(N):K]} \sum_{u | \mathfrak{p}} d_{u, \mathfrak{p}}\lambda_{u}(B) \\
                                                   & \geq \frac {1} {[K(N):K]} \sum_{u \in G_{\psi}v} d_{u, \mathfrak{p}} \lambda_{u}(B)\\
                                                   & \geq \frac {\log p} {2(q - 1)[K(N):K]}\, |G_{\psi}v| d_{v, \mathfrak{p}}  \\
                                                   & \geq  \frac {\log p} {2p^4(q - 1)}\\
                                                   & \geq \frac {\log p} {2p^6}
\end{split}
\] 
where we are using the estimates in $(4.3.1)$ and $(4.3.2)$ and the fact that $|G| = [K(N):K]$.\\
Therefore this choice of $B$ does the job. \hfill $\square$\\~\\ 
\textbf{Lemma 4.3.7 :} Assume that $ n \leq 1$. If $A \in E(K(N)) - E_{\text{tor}}$, there exists a nontorsion point $B \in E(\overline{\mathbb{Q}})$ with $\widehat{h}(B) \leq 4p^4( w + 1)( q - 1)^2\widehat{h}(A)$ and \[ \widehat{h}_{\mathfrak{p}}(B) \geq \frac {\log p} {2p^6}.\]
\textbf{Proof :} Without loss of generality we can assume that $n = 1$. Now consider two cases :\\~\\
\textbf{Case I :} There is a $\sigma \in \text{Gal}\,(K(N)|K)$ such that $[Q(1)](\sigma(A)) (= \sigma([Q(1)](A)) ) \notin E(\mathbb{Q}_{w}(N/p))$.\\
Apply lemma 4.3.6 to $\sigma(A)$ to get a $B \in E(\overline{\mathbb{Q}}) - E_{\text{tor}}$ such that $\widehat{h}(B) \leq 4\widehat{h}(\sigma(A)) = 4\widehat{h}(A)$ and \[\widehat{h}_{\mathfrak{p}}(B) \geq \frac {\log p} {2p^6}. \]
This $B$ has the desired properties.\\~\\
\textbf{Case II :} For all $\sigma \in \text{Gal}\,(K(N)|K)$, $\sigma([Q(1)](A)) \in E(\mathbb{Q}_{w}(N/p))$ . Then by lemma 3.4.2 we have $\sigma([Q(1)](A)) \in E(K(N/p))$. Now $p \nmid N/p$. So we can use lemma 4.3.5 to get a $B \in E(\overline{\mathbb{Q}}) - E_{\text{tor}}$ such that $\widehat{h}(B) \leq 2(w + 1)\widehat{h}([Q(1)](A)) = 2p^4(w + 1)(q - 1)^2\widehat{h}(A)$ and \[ \widehat{h}_{\mathfrak{p}}(B) \geq \frac {1} {2} \log p .\]
Clearly this $B$ does the job. Hence we are done in this case also.\\~\\
Thus we have found a $B$ with desired properties in all cases. \hfill $\square$\\~\\
\textbf{Lemma 4.3.8 :} If $A \in E(K(E_{\text{tor}})) - E_{\text{tor}}$, there is a nontorsion point $ B \in E(\overline{\mathbb{Q}})$ with $\widehat{h}(B) \leq 40p^4(w + 1)( q - 1)^2\widehat{h}(A)$ and \[ \widehat{h}_{\mathfrak{p}}(B) \geq \frac {\log p} {2p^6}.\]\\
\textbf{Proof :} Fix a positive integer $N$ such that $A \in E(K(N))$. Let $n , M$ be as before.\\
Use lemma 3.4.7 to get a $\sigma_0 \in \text{Gal}\,(K(N)|K)$ with properties as stated there.\\
Put $C = \sigma_{0}(A) - [2](A) \in E(K(N))$. Arguing as in proof of lemma 4.3.5 and using the hypothesis $A \notin E_{\text{tor}}$, we conclude that $C \notin E_{\text{tor}}$. \\
Using parallelogram law as before \[\widehat{h}(C) \leq 2(\widehat{h}(\sigma_{0}(A)) + \widehat{h}([2](A))) \leq 10 \widehat{h}(A) \tag{4.3.3}.\] 
Fix the least integer $n' \geq 0$ such that $C \in E(K(p^{n'}M))$. Clearly $n' \leq n$. Write $N' = p^{n'}M$. Consider two cases :\\~\\
\textbf{Case I :} $n' \leq 1$.\\
Apply lemma 4.3.7 to $C$ to get a $B \in E(\overline{\mathbb{Q}}) - E_{\text{tor}}$ such that $\widehat{h}_{\mathfrak{p}}(B) \geq \frac {\log p} {2p^6}$ and it satisfies \[\widehat{h}(B) \leq 4p^{4}(w + 1)(q - 1)^2\widehat{h}(C) \leq 40p^{4}(w + 1)(q - 1)^{2}\widehat{h}(A)\] where we are using $(4.3.3)$.\\~\\
\textbf{Case II :} $n' \geq 2$.\\
Use lemma 3.4.2 to get a $\sigma \in \text{Gal}\,(K(N')|K)$ with $\sigma(C)( = C') \notin E(\mathbb{Q}_{w}(N'/p))$. Fix a $\psi \in \text{Gal}\,(\mathbb{Q}_{w}(N')/\mathbb{Q}_{w}(N'/p))$ such that $\psi(C') = C'$.\\
Set $A' = \sigma(A)$ to get \[ C' = \sigma_{0}(A') - [2](A') \in E(K(N')) \tag{4.3.4}\] since $\sigma_0$ lies in the center.\\
We would like to apply lemma 4.3.6 to $C'$. For that we need to verify $[Q(n')](C') = [q](C') \notin E(\mathbb{Q}_{w}(N'/p))$.\\
Assume the contrary. Write $T = \psi(C') - C'$.\\
Since $[q](C') \in E(\mathbb{Q}_{w}(N'/p))$ , $\psi([q](C')) = [q](C')$ and $T \in E[q]$.\\
By choice of $\psi$, $T \neq O$. Applying $\psi$ to both sides of $(4.3.4)$ and using the fact $\sigma_0$ lies in the center we have \[ C' + T =\psi(C') = \sigma_{0}(\psi(A')) - [2](\psi(A')).\tag{4.3.5} \] Put $P = \psi(A') - A' \in K(N)$. $(4.3.4)$ and $(4.3.5)$ gives $T = \sigma_{0}(P) - [2](P)$.\\
Since $T \in E_{\text{tor}}$, we have $\widehat{h}(P) = \widehat{h}(\sigma_{0}(P)) = \widehat{h}([2](P)) = 4\widehat{h}(P)$ and thus $P\in E_{\text{tor}}$.\\
Choose a positive integer $M_1$ such that it is coprime to $p$ and $[M_{1}](P) \in E[p^{\infty}]$. By lemma 4.3.1 $[M_1](P) \in E[p^n]$ .\\
By construction of $\sigma_0$, $\sigma_{0}([M_1](P)) = [2M_{1}(P)]$. Thus $[M_1](T) = O$ and $T \in E[M_1]$. But $T \in E[q]$ and $\gcd (M_1, q) = 1$. So $T = O$. A contradiction ! \\
So $[Q(n')](C') \notin \mathbb{Q}_{w}(N'/p)$. Now use lemma 4.3.6 to get a $B \in E(\overline{\mathbb{Q}}) - E_{\text{tor}}$ such that $\widehat{h}_{\mathbb{p}}(B) \geq \frac {\log p} {2p^6}$ and \[\widehat{h}(B) \leq 4\widehat{h}(C') = 4\widehat{h}(C) \leq 40\widehat{h}(A)\]
where in the last step we are using $(4.3.3)$.\\~\\
Thus in all cases we have constructed a $B$ with desired properties. \hfill $\square$ \\~\\   
\section {Further estimates on $\widehat{h}_{v}$ }
For all $P \in E(\overline{\mathbb{Q}}) - \{O\}$ we have \[ \widehat{h}(P) = \frac {1} {[K :\mathbb{Q}]}\sum_{v \in M_K} d_{v}\widehat{h}_{v}(P). \tag{4.4.1} \]
We have computed some estimates for $\widehat{h}_{\mathfrak{p}}$. The goal of the current section is to find estimates for $\widehat{h}_{v}$ when $v \neq \mathfrak{p}$.\\
Habegger does it by using two equidistribution theorems. We shall follow his path.\\~\\
\subsection{Infinite places}
\textbf{Notation :} In this subsection $|\cdot|$ is the usual complex absolute value.\\
At first we recall some generalities.\\
Let $E$ be ab elliptic curve over a number field $F$ and $v$ be an infinite place of $F$.\\
Up to complex conjugation, $v$ determines an embedding of $F$ in $\mathbb{C}$. Fix one of the conjugate embeddings, call it $\sigma_{0}$.\\
We thus obtain an elliptic curve over $\mathbb{C}$ which we shall denote by $E_v$.\\
Let $\lambda_v : E_v(\mathbb{C}) - \{O\} \to \mathbb{R}$ be the local height function. By Weierstrass uniformization there is a $\tau \in \mathbb{C}$ with $\text{Im}(\tau) > 0$ such that there is an analytic isomorphism $\mathbb{C}/\mathbb{Z} + \tau\mathbb{Z} \to E_{v}(\mathbb{C})$. Write $q = e^{2\pi i \tau}$ and note that $|q| < 1$.\\
If $A \in E_{v}(\mathbb{C}) - \{O\}$ is the image of $z \in \mathbb{C}$ and $q(z) = e^{2\pi i z}$ then,
\[
\lambda_{v}(A)\\
         = - \frac {1} {2}b_2(\frac {\text{Im}(z)} {\text{Im}(\tau)})\log |q| -\log |1 - q(z)| - \sum_{n\geq1}\log |(1 - q^nq(z))(1 - q^nq(z)^{- 1})|     
\tag{4.4.2} \]\\
where $b_2(X) = X^2 - X + \frac {1} {6}$ . (See \cite{sil 2} )\\
The group $E_v(\mathbb{C})$ endowed with complex topology is compact and let $\mu_{v, E}$ denote the unique Haar measure on it of total mass $1$.\\
Now we need an equidistribution theorem which is originally due to Szpiro, Ullmo, Zhang. We shall quote the simplified version given in $\cite{hab}$ :\\~\\
\textbf{Proposition 4.4.1 :} Let $P_1, P_2, P_3 ,\cdots \in E(\overline{\mathbb{Q}}) - E_{\text{tor}}$ be a sequence of points such that $\widehat{h}(P_k) \to 0$. Let $f :E_{v}(\mathbb{C}) \to \mathbb{R}$ is a continuous function, then \[ \lim_{k \to \infty} \frac {1} {[F(P_k) : F]}\sum_{\sigma}f(\sigma(P_k)) = \int f\mu_{E, v},\] where $\sigma$ runs over all the all the field embeddings $F(P_k) \to F$ which extend $\sigma_{0}$.\\~\\
Now are back to our set up.\\
Say, $v$ is an infinite place of $K$. Call the corresponding embedding $\sigma_0$ as before.\\
Then \[\widehat{h}_{v}(P)  = \frac {1} {[K(P) : K]} \sum_{u | v, u \in M_{K(P)}} d_{u, v} \lambda_u(P)\] for $P \in E(\overline{\mathbb{Q}}) - \{O\}$.\\
Clearly, \[\widehat{h}_{v}(P) = \frac {1} {[K(P) : K]}\sum_{\sigma} \lambda_{v}(\sigma(P)) \] where the sum runs over all $\sigma : F(P) \to \mathbb{C}$ such that it extends $\sigma_{0}$.\\~\\
Let $P_1, P_2, P_3, \cdots \in E(\overline{\mathbb{Q}}) - E_{\text{tor}}$ such that $\widehat{h}(P_k) \to 0$.\\
We would like to use proposition 4.4.1 to obtain an estimate for $\liminf \widehat{h}_{v}(P_k)$. One can not directly use proposition 4.4.1 since one can not continuously extend $\lambda_v$ to $O$.\\
So for each positive integer $m$ we define $\lambda_{v, m}(P) = \min \{ \lambda_{v}(P), m \} $ for all $P \in E(\mathbb{C}) - \{O\}$ and $\lambda_{v, m}(O) = m$.\\
It is easy to see that $\{\lambda_{v, m}\}_{m \geq 1}$ is a sequence of continuous function which increases pointwise to $\lambda_v$. By monotone convergence theorem, we have that $\lambda_v$ is measurable and $a_m \to \int \lambda_{v}\,\mu_{v, E}$ where $a_m = \int \lambda_{v, m}\mu_{v, E}$ for each $m \geq 1$.\\
Now note that \[\widehat{h}_{v}(P_{k}) \geq \frac {1} {[K(P_k) : K]} \sum_{\sigma} \lambda_{v, m}(\sigma(P))\] holds for each positive integer $k, m$.\\
At first let $k \to \infty$ and use proposition 4.4.1 to conclude that \[\liminf \widehat{h}_{v}(P_k) \geq a_m.\]
This holds for each positive integer $m$ . Letting $m \to \infty$ we have \[\liminf \widehat{h}_{v}(P_k) \geq \int \lambda_{v}\mu_{v, E}.\]
Now we have the following lemma:\\~\\
\textbf{Lemma 4.4.2 :} $\int \lambda_v\mu_{v, E} = 0.$\\~\\
\textbf{Proof :} Let $\tau$ denote the corresponding parameter in Weierstrass uniformization of $E_v$. The points on $E_v$ are parameterized by the fundamental parallelogram \[ \{ x + y\tau \; | \; 0 \leq x, y < 1 \} \subseteq \mathbb{C}\] and for this parameterization the Haar measure is nothing but the usual measure on unit square.\\
We shall use the expression in $(4.4.1)$  and integrate term by term. This is allowed since the series converges absolutely ( follows from the fact that $|q| < 1$).\\
In the following computation $z = x + y\,\tau$. So $\text{Im}(z) = y\,\text{Im}(\tau)$.\\
Now \[\int_{0}^{1}\int_{0}^{1} - \frac {1} {2}b_{2}(\frac {\text{Im} (z)} {\text{Im} (\tau)}) \log |q| dy dx = - \frac {\log |q|} {2} \int_{0}^{1} (\int_{0}^{1} \, b_2(y) dy) dx.\]
Easy computation shows that the inner integral is zero.\\
Hence the first term integrates to $0$.\\
Now let $\epsilon \in \{\pm 1\}$ and $n \geq 1$.\\
Then 
\[
\begin{split}
\int_{0}^{1}\int_{0}^{1} \log | 1 - q^{n}q(z)^{\epsilon} | dxdy & = \int_{0}^{1}( \int_{0}^{1} \log | e^{- 2\pi i \epsilon x} - e^{2\pi i \tau (n + \epsilon y)}| dx) dy \\
                                                                                                        & = \int_{0}^{1} \log^{+} (|\,\text{exp}\,(2\pi i \tau (n + \epsilon y))|) dy \\
\end{split}
\] 
where in the last step we are using the Jensen's formula. Since $\text{Im}(\tau) > 0$, $ n\geq 1$ and $ 0 \leq y < 1$ we conclude that $|\,\text{exp}\,(2\pi i \tau ( n+ \epsilon y))| < 1$. So the last integral is $0$.\\
Similar computation shows $| 1 - q(z)|$ integrates to zero.\\
Hence the lemma. \hfill $\square$ \\~\\
The discussion of this subsection can be summarized in the following lemma:\\~\\
\textbf{Lemma 4.4.3 :} Let $v$ be an infinite place of $K$ and let $P_1, P_2, P_3, \cdots \in E(\overline{\mathbb{Q}}) - E_{\text{tor}}$ such that $\widehat{h}(P_k) \to 0$. Then \[ \liminf \widehat{h}_v(P_k) \geq 0.\]

\subsection{Finite places}
We shall recall some general facts about elliptic curves over non-Archimedean local fields. For details see \cite{sil 2}, chapter VI, section 4.\\
Let $E$ be an elliptic curve over a number field $F$. Assume that $v$ is a finite place of $F$. Consider $E$ as a curve in $F_v$.\\
If $v$ is a place of good reduction then for all $A \in E(\overline{F_v}) - \{O\}$, $\lambda_{v}(A) \geq 0$.\\
If $v$ is a place of split multiplicative reduction , then by Tate uniformization there is $q_v \in F_{v}^{\times}$ with $|q_v|_v < 1$ such that there is a surjective homomorphism $\phi : \overline{F_v}^{\times} \to E(\overline{F_v})$ with kernel $q_v^{\mathbb{Z}}$, the cyclic group generated by $q_v$. \\
It is easy to see that for any $A \in E(\overline{F_v}) - \{O\}$ there is an unique $q_{0}(A) \in \overline{F_v}^{\times}$ such that $ |q_v|_v < |q_{0}(A)|_v \leq 1$. Now the local height is given by
\[\lambda_v(A) = - \frac {1} {2}b_{2}(\frac {\log |q_0(A)|_v} {\log |q_v|_v})\log |q_v|_v - \log|1 - q_{0}(A)|_v \tag{4.4.3}\] 
Further using ultrametric triangle inequality $ |1 - q_0(A)|_v \leq 1$. Thus for all $A \in E(\overline{F_v}) - \{O\}$ we have \[ \lambda_v(A) \geq - \frac {1} {2}b_2(\frac {\log |q_0(A)|_v} {\log|q_v|_v}) .\tag{4.4.4} \]
It is known that $\phi$ commutes with the action of the absolute galois group of $F_v$. So $\lambda_v$ is invariant under action of $\text{Gal}(\overline{F}_v/F_v)$.\\
Note that we have a well defined map $l_v : E(\overline{F_v}) \to \mathbb{R}/\mathbb{Z}$ given by \[l_v(A) = \frac {\log |q(A)|_v} {\log |q_v|_v} + \mathbb{Z}\] where $q(A) \in \overline{F_v}^{\times}$ such that $\phi(q(A)) = A$.\\
We identify the topological group $\mathbb{R}/\mathbb{Z}$ with unit circle and equip it with the Haar measure $\mu_{\mathbb{R}/\mathbb{Z}}$ of lotal mass 1.\\
Now we have the following theorem which is due to Chambert-Loir. We quote the simplified version in \cite{hab} :\\~\\
\textbf{Proposition 4.4.4 :} Let $P_1, P_2, P_3, \cdots \in E(\overline{\mathbb{Q}}) - E_{\text{tor}}$ be a sequence of points such that $\widehat{h}(P_k) \to 0$.\\
If $f : \mathbb{R}/\mathbb{Z} \to \mathbb{R}$ is a continuous function, then 
\[ \lim_{k \to \infty} \frac {1} {[F(P_k) : F]} \sum_{\sigma} f(l_{v}(\sigma (P_k))) = \int f\mu_{\mathbb{R}/\mathbb{Z}} \]
where $\sigma$ runs over all field embeddings $F(P_k) \to \overline{F_v}$ which are identity on $F$.\\~\\
Now we are back to our set up.\\
By standard theory of elliptic curves we know that there is a finite, galois extension $L/K$ such that thought as an elliptic curve over $L$, $E$ has either good or split multiplicative reduction at finite places.\\
Let $\mathfrak{p}_1,\cdots,\mathfrak{p}_s$ be the prime ideals of ring of integers of $K$ which appear in the denominator of $j$-invariant of $E$. (since $\mathfrak{p}$ is a place of good reduction we conclude that $\mathfrak{p} \neq \mathfrak{p}_i$ for any $i$ in the range)\\
The reduction type at some finite place $u \in M_L$ is given by the following rule :\\
if $u \nmid \mathfrak{p}_i$ for any $i$ then $u$ is a place of good reduction,\\
if $u \, | \mathfrak{p}_i$ for some $i$ then $u$ is a place of split multiplicative reduction.\\
Let $v$ be a finite place of $K$ such that $v \notin \{\mathfrak{p}, \mathfrak{p}_1,\cdots,\mathfrak{p}_s \}$.\\
Then for any $u \in M_L$ with $u | v$, $u$ is a place of good reduction and hence $\lambda_u$ is nonnegative. Thus $\widehat{h}_v(P) \geq 0$ for all $P \in E(\overline{\mathbb{Q}}) - E_{\text{tor}}$.\\
Now let $v = \mathfrak{p}_i$ for some $i$.\\
Fix a place $u \in M_L$ such that $u |v$. Fix an algebraic closure $\overline{L_u}$ of $L_u$. Since $\lambda_u$ is invariant under the action of absolute galois group of $L_u$ we conclude from the definition of $\widehat{h}_{v}$ that \[\widehat{h}_{v}(A) = \frac {1} {[L:K]} \sum_{\tau}\frac {1} {[L(A) : L]} \sum_{\sigma} \lambda_{u}(\sigma(A)) \tag{4.4.5}\]
where the first sum runs over all the field embeddings $\tau : L \to \overline{L_u}$ such that $\tau|_K$ is the map $K \to K_v$ and the second sum runs over field embeddings $\sigma: L(A) \to \overline{L_u}$ and $A \in E(\overline{\mathbb{Q}}) - \{O\}$.\\
Fix a $\tau$ as above. So we have an elliptic curve $E$ over $\overline{L_u}$. \\
Let $q_u$ be the uniformizer for this curve.\\
Note that $b_2(0) = b_2(1)$. Thus it defines a continuous function on $\mathbb{R}/\mathbb{Z}$. Call it $B_2$.\\
Further it follows from definitions that \[b_2(\frac {\log |q_0(A)|_u} {\log |q_u|u}) = B_2(l_u(A))\] for all $A \in E(\overline{\mathbb{Q}}) - \{O\}$.\\  
Now let $P_1, P_2, P_3 , \cdots \in E(\overline{\mathbb{Q}}) - E_{\text{tor}}$ be a sequence of points such that $\widehat{h}(P_k) \to 0$. \\
Choose a $\tau_{1} \in \text{Gal}\,(\overline{\mathbb{Q}} | K)$ such that $\tau \circ \tau_{1} : L \to L_{u}$ is the identity mapping. Then the extensions of this map are given precisely by $\sigma \circ \tau_1$ where $\sigma$ is an extension of $\tau$.\\
Consider the sequence of points $\{\tau_1^{- 1}(P_k) \}_{k \geq 1}$. Clearly they are elements of $E(\overline{\mathbb{Q}}) - E_{\text{tor}}$. Further $\widehat{h}(\tau^{-1}(P_k)) \to 0$.\\
Using proposition 4.4.4 we have 
\[
\begin{split}
\lim_{k \to \infty} \frac {1} {[L(P_k) : L]}\sum_{\sigma} B_{2}(l_u(\sigma(P_k))) & = \lim_{k \to \infty} \frac {1} {[L(\tau_{1}^{- 1}P_k) :L]}\sum_{\sigma_{1}=\sigma \circ \tau_{1}} B_{2}(l_u(\sigma_{1}(\tau_1^{- 1}P_k)))\\
                                                                                                                                & = \int B_{2}\mu_{\mathbb{R}/\mathbb{Z}}
\end{split}
\] where we are using that $[L(P_k): L] = [L(\tau_1^{- 1}(P_k)) : L]$ for each $k$.\\
One can evaluate the integral easily as before, and see that it evaluates to $0$.\\
Use $(4.4.4)$ to conclude $\liminf \frac {1} {[L(P_k): L]}\sum_{\sigma}\lambda_{v}(\sigma(P_k)) \geq 0$.\\
True for each $\tau$ and there are only finitely many of them.\\
Using $(4.4.5)$ one has $\liminf \widehat{h}_v(P_k) \geq 0$.\\
The discussion of this subsection can be summarized as follows :\\~\\
\textbf{Lemma 4.4.5 :} Let $v \neq \mathfrak{p}$ be a finite place of $K$. Let $P_1, P_2, \cdots \in E(\overline{\mathbb{Q}}) - E_{\text{tor}}$ be such that $\widehat{h}(P_k) \to 0$. Then $\liminf \widehat{h}_v(P_k) \geq 0$.  
\section {Proof of theorem 0.2}
Assume that theorem 0.2 does not hold.\\
Then there is a sequence of distinct points $A_1, A_2, A_3 \cdots \in E(\overline{\mathbb{Q}}) - E_{\text{tor}}$ such that $\widehat{h}(A_k) \to 0$.\\
Use lemma 4.3.8 to $A_k$ to get a $B_k \in E(\overline{\mathbb{Q}}) - E_{\text{tor}}$ such that \[\widehat{h}(B_k) \leq 40p^4(w + 1)(q - 1)^{2}\widehat{h}(A_k)\] and \[\widehat{h}_{\mathfrak{p}}(B_k) \geq \frac {\log p} {2p^6}\] for each $k \geq 1$. Clearly $\widehat{h}(B_k) \to 0$.\\
By lemma 4.4.3 and 4.4.5 we have $\liminf \widehat{h}_{v}(B_k) \geq 0$ for each $v \neq \mathfrak{p}$.\\
Use $(4.4.1)$ to conclude that $\liminf \widehat{h}(B_k) \geq \frac {d_{\mathfrak{p}}\log p} {2[K : \mathbb{Q}]p^6} > 0$.\\
A contradiction.\\
This contradiction proves the theorem. \hfill $\square$

\begin{appendices}
\chapter{Division points of formal groups}
\section{Introduction}
Let $p$ be a prime number and let $K$ be a finite extension of $\mathbb{Q}_p$. Put $O_K$ to be the ring of integers of $K$, let $\mathfrak{p}_K$  
denote the unique maximal ideal of $O_K$ and $v_{\mathfrak{p}_K}(\cdot)$ be the valuation associated to it. Fix an algebraic closure $\overline{\mathbb{Q}}_p$ and $|.|_p$ be an fixed extension of the absolute value. Let $\overline{O}$ be the ring of integers of $\overline{\mathbb{Q}}_p$ and $\overline{\mathfrak{p}}$ be the unique maximal ideal of $\overline{O}$. Clearly $\overline{\mathfrak{p}} \cap K = \mathfrak{p}_K$.\\
Fix a generator $\pi$ of $\mathfrak{p}_K$ and let $K_{R} = \mathbb{Q}_p(\pi)$. Use $R$ to denote the ring of integers of $K_R$, $\mathfrak{p}_R$ to denote the maximal ideal and $f_R$ to denote the degree of the residue extension.\\ 
Let $\mathfrak{F}$ be a (one dimensional, commutative) formal group law defined over $O_K$ which satisfies following additional conditions :\\
i) $\mathfrak{F}$ has a formal $R$ module structure ,\\
ii) if \[[\pi](X) = \pi X + a_2X^2 +a_3X^3 +a_4X^4 +\cdots \] then $\min \, \{ i \geq 2 | |a_i|_{p} = 1 \} = p^h$ for some positive integer $h$. 
The integer $h$ will be called the height of group law. This condition is satisfied unless all the $a_i$-s are in maximal ideal (see \cite{haz2}, 18.3.1).\\
$\mathfrak{F}$ defines a $R$ module structure on $\mathfrak{p}_K$ which naturally extends to a $R$ structure on $\overline{\mathfrak{p}}$. We shall denote the corresponding addition by $\oplus_{\mathfrak{F}}$ to distinguish it from usual addition.\\
For each $n \geq 1$, let $\mathfrak{F}[\pi^n]$ denote the $\pi^n$-torsion submodule of $\overline{\mathfrak{p}}$ and let $K(\pi^n)$ be the subfield of $\overline{\mathbb{Q}}_p$ generated by $\mathfrak{F}[\pi^n]$ over $K$. We shall adopt the convention $\mathfrak{F}[\pi^0] = \{0\}$.\\
Note that from the condition $(i)$ on $\mathfrak{F}$ implies \[[\pi^n](X) = \pi^n X + \text{higher degree terms} \tag{1.1}\] for each $n \geq 1$.\\
First we shall prove :\\~\\
\textbf{Proposition 1.1 :} With the set up described above one has $f_R\,|\,h$.\\~\\
The main goal of this appendix is to prove :\\~\\
\textbf{Theorem 1.2 :} Let $n \geq 1$. Put $q = p^h$ and $h_{r} = \frac{h}{f_R}$. Then the following statements are true :\\
i) $\mathfrak{F}[\pi^n] \cong (R/\pi^nR)^{h_{r}}$ as $R$ modules.\\
ii) If $z \in \mathfrak{F}[\pi^n] - \mathfrak{F}[\pi^{n-1}]$, then $K(z)/K$ is a totally ramified extension of degree $q^{n-1}(q - 1)$.\\~\\
\textbf{Theorem 1.3 :} Let $n, q, h_r$ be as in statement of theorem 1.2.\\
Let $f$ be the degree of the extension of residue fields associated to the extension $K/\mathbb{Q}_p$. We shall assume that $h\,|\,f$.\\
Put $K^{h_r}_{R}$ to be the unique unramified extension of degree $h_r$ of $K_R$ in $\overline{\mathbb{Q}}_p$. Use $R^{h_r}$ to denote the ring of integers and $U_{i,K^{h_r}_R}$ ($i \geq 0$) to denote the $i$-th unit group of this field. Then :\\
i) $K(z) = K(\pi^n)$ for any $z \in \mathfrak{F}[\pi^n] - \mathfrak{F}[\pi^{n-1}]$.\\
ii) $K(\pi^n)/K$ is a Galois extension with \[\text{Gal}\,(K(\pi^n)|K) \cong U_{0,K_{R}^{h_r}}/U_{n,K_{R}^{h_r}}.\]
iii) Let $k$ and $i$ be integers such that $1 \leq k \leq n$ and $q^{k-1} \leq i \leq q^k - 1$. Then \[G_{i}(K(\pi^n)|K) = \text{Gal}(K(\pi^n)|K(\pi^{k})).\]\\
We shall give a proof of the theorem modulo the following assumption :\\~\\
\textbf{Assumption :} $K(z) = K(\pi^{n+1})$ for any $z \in \mathfrak{F}[\pi^{n+1}] - \mathfrak{F}[\pi^{n}]$ ie an analogue of lemma 2.3.7 in this case. We shall discuss more about it in Appendix-B. \\~\\
\textbf{Remark 1.4 :} i) If $K/\mathbb{Q}_p$ is an unramified extension then one can take $\pi = p$. In this case the condition $(i)$ on $\mathfrak{F}$ is verified trivially and $R = \mathbb{Z}_p$.\\
ii) Note that our hypothesis on $\mathfrak{F}$ is quite weaker than the Lubin-Tate hypothesis (see \cite{neu}, chapter 3, section 6) though we prove similar Galois theoritic properties that one expects from the Lubin-Tate theory (\cite{neu}, chapter 3, section 7 and 8).\\
iii) As a consequence of our computations we shall be able to prove that the $\pi$-adic Tate module $T(\pi) = \varprojlim \mathfrak{F}[\pi^n]$ has a $R^{h_r}$ module structure and it will turn out to be a free $R^{h_r}$ module of rank $1$. We shall also see that the absolute Galois group $\text{Gal}(\overline{\mathbb{Q}}_p/K)$ acts on $T(\pi)$ via $R^{h_r}$ module morphisms.\\
iv) We shall use the notations introduced in this section throughout this appendix.\\~\\

\section{Structure of the torsion subgroups}
We start by noting a standard result in commutative algebra :\\~\\
\textbf{Proposition 2.1 (Preparation Theorem) :} Let $(A, \mathfrak{m})$ be a local ring which is complete with respect to $\mathfrak{m}$-adic topology. Now \[f(X) = a_0 + a_1X + a_2X^2 +\cdots \in A[[X]]\] be such that there exists an $i \geq 0$ satisfying $a_i \notin \mathfrak{m}$. Put \[s(f) = \min_{i \geq 0} \{ i | a_i \notin \mathfrak{m} \}.\] Then there is a unique ordered pair $(u(X), F(X))$ satisfying the following conditions:\\ $u(X)$ is a unit in $A[[X]]$, $F(X)$ is a monic polynomial of degree $s(f)$ such that all the coefficients (except the leading one) of $F(X)$ come from $\mathfrak{m}$ and  $f(X) = u(X)F(X)$.\\~\\
\textbf{Proof :} See \cite{Bourbaki}, chapter VII, section 3. \hfill $\square$\\~\\
In our context $A = O_K$, $\mathfrak{m} = \mathfrak{p}_K = \pi O_K$.\\
Let $f_n(X) = [\pi^n](X) \in O_K[X]$ for all $n \in \mathbb{N}$.\\
By condition $(ii)$ on $\mathfrak{F}$ we have $s(f_n) = p^{nh}$.\\
Write $f_n(X) = u_n(X)F_n(X)$ as in proposition 2.1.\\
Note that if $u(X) \in O_K[[X]]$ is an unit then the constant term of $u(X)$ (say $c_0$) is a unit in $O_K$. An easy application of ultrametric triangle inequality gives $|u(z)|_p = |c_0|_p \neq 0$ for all $z \in \overline{\mathfrak{p}}$.\\
Clearly $z\in \mathfrak{F}[\pi^n]$ implies $f_n(z) = 0$. Using the observation noted above we have $F_n(z) = 0$. Since degree of $F_n(X)$ is exactly $p^{nh}$ this implies \[|\mathfrak{F}[\pi^n]| \leq p^{nh} \tag{2.1} \]    
holds for all $n \in \mathbb{N}$.\\~\\
Now we shall prove a sequence of lemmas which will give us information about the torsion points :\\~\\
\textbf{Lemma 2.2 :} Let $n \in \mathbb{N}$. $z \in \overline{\mathbb{Q}}_p$ be such that $F_n(z) = 0$. Then $z \in \overline{\mathfrak{p}}$.\\~\\
\textbf{Proof :} For simplicity write $s(f_n) = s_n$. Then \[ F_n(X) = X^{s_n} + b_1X^{s_n - 1} + \cdots + b_{s_n}\tag{2.2} \] where $ b_i \in \mathfrak{p}_K$ for all $1 \leq i \leq s_n$.\\
By hypothesis $F_n(z) = z^{s_n} + \cdots + b_{s_n} = 0$. Now again using ultrametric triangle inequality we see that one must have $|z|_p < 1$. From here the lemma follows easily. \hfill $\square$\\~\\
\textbf{Lemma 2.3 :} Let $n \in \mathbb{N}$. Then $z \in \mathfrak{F}[\pi^n]$ if and only if $F_n(z) = 0$.\\~\\
\textbf{Proof :} The discussion before lemma 2.2 proves if $z \in \mathfrak{F}[\pi^n]$ then $F_n(z) = 0$.\\
Conversely, if $F_n(z) = 0$ then $[\pi^n](z) = f_n(z) = 0$. By lemma 2.2 we have $z \in \overline{\mathfrak{p}}$. Hence $z \in \mathfrak{F}[\pi^n]$.\\
This proves the lemma. \hfill $\square$\\~\\
\textbf{Lemma 2.4 :} Let $n \in \mathbb{N}$. Put $F_0(X) = X$. Then the following holds :\\
i) $F_{n-1}(X) | F_n(X)$ in $O_K[X]$.\\
ii) Let $g_{n}(X) = \frac{F_n(X)}{F_{n-1}(X)} \in O_{K}[X]$. Then $g_{n}(X)$ is an irreducible Eisenstien polynomial of degree $q^{n - 1}(q - 1)$.\\
iii) $\gcd\,(g_{n}(X), F_{n-1}(X)) = 1$ and $F_n(X) = Xg_{1}(X) \cdots g_{n}(X)$ is a prime factorization of $F_n(X)$ in distinct prime factors.\\
iv) $F_n(X)$ has exactly $p^{nh}$ distinct roots in $\overline{\mathbb{Q}}_p$.\\~\\
\textbf{Proof :} Let $n \geq 1$ and $s_n, b_1,\cdots, b_{s_n}$ be as in proof of lemma 2.2.\\  
Note that in $[\pi^n](X)$, the constant term is $0$ and the coefficient of $X$ is $\pi^n$. But $u_n(X)$ is an unit in $O_K[[X]]$. Hence $b_{s_n} = 0$ and $v_{\mathfrak{p}_K}(b_{s_n -1}) = n$.\\ 
Now we shall prove the lemma induction on $n$.\\
\emph{Base case :} By definition $F_{0}(X) = X$. So first we need to show that $X | F_1(X)$. But this easily follows from the fact that the constant term of $[\pi](X)$ is zero.\\   
Put $g_{1}(X) = \frac{F_{1}(X)}{X} \in O_{K}[X]$. We know that it is a monic polynomial of degree $(q - 1)$. Write \[g_1(X) = X^{q-1} + a_1X^{q-2} +a_2X^{q-3} +\cdots +a_{q-1}.\] By discussion in the beginning of the proof $v_{\mathfrak{p}_K}(a_{q-1}) = 1$. Now $F_{1}(X) = Xg_{1}(X)$. Reducing both sides modulo $\mathfrak{p}_K$ we see that the left hand side is a monomial. Hence $g_{1}(X)\,\text{mod}\,\mathfrak{p}_{K}$ must be a monomial. Thus $a_{i} \in \mathfrak{p}_{K}$ for all $1 \leq i \leq q -1$. This with $v_{\mathfrak{p}_{K}}(a_{q-1}) = 1$ proves $(ii)$.\\
$(iii)$ follows directly from $(ii)$ and the observation that $a_{q-1} \neq 0$.\\
$(iv)$ follows from the prime factorization given in $(iii)$.\\
Thus we have proved the base case.\\
Assume that the lemma is true for some $k \in \mathbb{N}$. We would like to prove it for $n = k + 1$.\\
\emph{Induction step :} Clearly $z \in \mathfrak{F}[\pi^k]$ implies $z \in \mathfrak{F}[\pi^{k+1}]$. Thus each root of $F_{k}$ is also a root of $F_{k+1}$. By induction hypothesis (part $(iv)$) each root of $F_{k}$ has multiplicity $1$. Hence $F_{k}(X)\,|\,F_{k+1}(X)$ in $K[X]$. Further since $F_{k}(X)$ is monic and $F_{k}(X),\,F_{k+1}(X)$ have coefficients in $O_K$, by division algorithm the quotient $\frac{F_{k+1}(X)}{F_{k}(X)}$ also has coefficients in $O_K$. This proves $(i)$.\\
Let $g_{k+1}(X) = \frac{F_{k+1}(X)}{F_{k}(X)}.$ Clearly it has degree $q^{k+1} - q^{k} = q^k(q - 1)$ and it is monic. Write \[g_{k+1}(X) = X^{q^k(q-1)} + a_1X^{q^k(q-1)-1} + \cdots + a_{q^k(q-1)}\] where $a_i \in O_K$ for all $1 \leq i \leq q^k(q-1)$.\\
Now $F_{k+1}(X) = F_{k}(X)g_{k+1}(X)$. The left hand side is a monomial modulo $\mathfrak{p}_K$. So $g_{k+1}(X)$ is also a monomial modulo $\mathfrak{p}_K$. Thus $a_i \in \mathfrak{p}_K$ for all $1 \leq i \leq q^k(q-1)$.\\
Let $b$, $b'$ be coefficients of $X$ in $F_{k+1}(X),\,F_{k}(X)$ respectively. Then \[b = b'a_{q^k(q-1)}.\]
By the discussion in the beginning $v_{\mathfrak{p}_K}(b) = k+1$, $v_{\mathfrak{p}_K} = k$. Thus from the relation above we conclude that $v_{\mathfrak{p}_K}(a_{q^k(q-1)}) = 1$.\\
Hence $g_{k+1}(X)$ is an Eisenstein polynomial. Hence $(ii)$.\\
Note that $\text{deg}_{X}(F_k(X)) = q^k \leq q^k(q-1) = \text{deg}_{X}(g_{k+1}(X))$. But  $g_{k+1}(X)$ is irreducible while $F_{k}(X)$ is reducible. So $\gcd\,(F_k(X), g_{k+1}(X)) = 1$. The second part follows from induction hypothesis on $F_k(X)$, part $(ii)$ and the fact that $\gcd\,(F_k(X), g_{k+1}(X)) = 1.$\\
$(iv)$ follows from the prime factorization in $(iii)$ (since we are working in a field of characteristic $0$ each polynomial is separable).\\
Thus we have proved the lemma for $n = k+1$.\\
Now the lemma follows from principle of mathematical induction. \hfill $\square$ \\~\\
\textbf{Corollary 2.5 :} Let $n \geq 1$. Then :\\
i) $|\mathfrak{F}[\pi^n]| = p^{nh}$,\\
ii) $z \in \mathfrak{F}[\pi^n] - \mathfrak{F}[\pi^{n-1}]$ if and only if $g_{n}(z) = 0$.\\~\\
\textbf{Proof :} $(i)$ follows from lemma 2.3 and lemma 2.4 part $(iv)$.\\
$(ii)$ follows from lemma 2.3 and lemma 2.4 part $(iii)$. \hfill $\square$\\~\\
\textbf{Proof of theorem 1.2(ii) :} Let $n \in \mathbb{N}$ and assume that $z \in \mathfrak{F}[\pi^n] - \mathfrak{F}[\pi^{n-1}]$. By corollary 2.5, $z$ is a root of $g_n(X)$. But by lemma 2.4 $(ii)$ we have $g_n(X)$ is an Eisenstein polynomial of degree $q^{n-1}(q-1)$. Now from standard theory of Eisenstein polynomials it follows that $K(z)/K$ is a totally ramified extension of degree $q^{n-1}(q-1)$. \hfill $\square$\\~\\
\textbf{Proof of proposition 1.1 :} $\mathfrak{F}[\pi]$ is a finite $R$ module which is annahilited by $\pi$. Hence it is a finite dimentional $R/\pi R$ vector space. Say the dimension is $d$. Then $|\mathfrak{F}[\pi]| = |R/\pi R|^d$. The left hand side is $p^h$ and the right hand side is $p^{df_R}$. Hence $h = df_R$. Thus $f_R | h$ as desired. \hfill $\square$ \\~\\
\textbf{Corollary 2.6 :} $\mathfrak{F}[\pi] \cong (R/\pi R)^{h_r}$ as $R$ module.\\~\\
\textbf{Proof :} From the proof of proposition 1.1 it follows that $\mathfrak{F}[\pi] \cong (R/\pi R)^{d} $ as $R/\pi R$ vector spaces where $d = \frac{h}{f_R} = h_r$. Since $R/\pi R$ vector space structure on both sides descend from $R$ module structure, we conclude the result. \hfill $\square$\\~\\
\textbf{Proof of theorem 1.2(i) :} Let $n \in \mathbb{N}$.\\
Clearly $\mathfrak{F}[\pi^n]$ is annhillited by $\pi^n$. By structure theorem of finitely generated modules over PID there are integers $k, n_1, n_2, \cdots n_k$ with $k \in \mathbb{N}$ and $1\leq n_1 \leq n_2 \cdots \leq n_k \leq n$ such that \[\mathfrak{F}[\pi^n] \cong R/\pi^{n_1}R \oplus R/\pi^{n_2}R \oplus \cdots \oplus R/\pi^{n_k}R \tag{2.3}\] as $R$ modules.\\
If $M$ is a $R$ module let $M_{\pi}$ be the $\pi$-torsion submodule of $M$. It is easy to see that $\mathfrak{F}[\pi^n]_{\pi} = \mathfrak{F}[\pi]$. Now from $(2.3)$ we have \[\mathfrak{F}[\pi] \cong (R/\pi^n_1)_{\pi} \oplus \cdots \oplus (R/\pi^n_k)_{\pi}\tag{2.4}\] as $R$ module.\\
Now $(R/\pi^{n_i}R)_{\pi} = \pi^{n_i - 1}R/\pi^{n_i}R \cong R/\pi R$ as $R$ module, for all $1 \leq i \leq n$. Using $(2.4)$ we have $\mathfrak{F}[\pi] \cong (R/\pi R)^k$ as $R$ module. So $|\mathfrak{F}[\pi]| = p^{f_Rk}$. Comparing with the corollary 2.6 we conclude that $k = h_r$.\\
Now left hand side of $(2.3)$ has cardinality $p^{nh}$ while the cardinality of right hand side is $p^{f_R(n_1 + \cdots + n_k)}$. So $nh_r = n_1 + \cdots + n_k$. But $k = h_r$ and $n_i \leq n$ for all $1 \leq i \leq h_r$. Thus $n_i = n$ for all $1 \leq i \leq h_r$.\\
From here the statement follows by the isomorphism in $(2.3)$. \hfill $\square$\\~\\
Before concluding the section we note a corollary :\\~\\
\textbf{Corollary 2.7 :} Let $n \in \mathbb{N}$ and assume that $z \in \mathfrak{F}[\pi^n] - \mathfrak{F}[\pi^{n-1}]$. Let $\mathfrak{p}_{z}$ be the unique maximal ideal in the ring of integers of $K(z)$. Then $z$ is a generator of $\mathfrak{p}_z$ and $|z|_{p} = |\pi|_p^{\frac{1}{q^{n-1}(q-1)}}$.\\~\\
\textbf{Proof :} By lemma 2.4 and corollary 2.5 $z$ is a root of an irreducible Eisenstein polynomial $g_n(X) \in O_{K}[X]$. This proves the first part.\\
The second part follows from theorem 1.2(i) and the assertion in first part. $\square$\\~\\ 
\textbf{Remark 2.8 :} Results in this section give a new proof of lemma 2.3.4.\\

\section{Computation of Galois group}
Let $n \in \mathbb{N}$. We shall continue to use the notations $f_n, F_n, g_n, u_n$ as introduced in previous section. The notation $\mu_n$ will be used to denote the group of $n$-th roots of unities in $\overline{\mathbb{Q}}_p$.\\
In this section we shall work with the additional assumption $h | f$ as in statement of theorem 1.3. At first we gather information about the Galois group and ramification groups by simple results in algebraic number theory and use them to compute the ramification groups and the Galois group.  
We begin by noting a simple fact :\\~\\
\textbf{Lemma 3.1 :} Let $n \in \mathbb{N}$. Then $K(\pi^n)/K$ is a Galois extension.\\~\\
\textbf{Proof :} By lemma 2.3 $K(\pi^n)$ is the splitting field of $F_{n}(X)$ over $K$. Hence the conclusion. \hfill $\square$\\~\\        
\textbf{Lemma 3.2 :} i) Let $z_0 \in \mathfrak{F}[\pi] - \{0\}$. Then $K(z_0) = K(\pi)$.\\
ii) $\text{Gal}(K(\pi)/K) \cong \mathbb{Z}/(q-1)\mathbb{Z}.$\\~\\
\textbf{Proof :} By lemma 2.4(ii) and corollary 2.5(ii) we know that the elements of $\mathfrak{F}[\pi] - \{0\}$ are exactly the conjugates of $z_0$ over $K$. But they generate $K(\pi)$ over $K$. Thus $K(\pi)$ is nothing but the Galois closure of $K(z_0)$ over $K$. So to prove the first part of the lemma it is enough to show that $K(z_0)/K$ is Galois.\\
By theorem 1.2(ii) we have $K(z_0)/K$ is a totally ramified extension of degree $(q-1)$. Since $p\,|\,q$ this extension is tamely ramified. Let $\mathfrak{p}_{K(z_0)}$ denote the unique maximal ideal of the ring of integers of $K(z_0)$. Applying a standard result in algebraic number theory (proposition-12 in chapter-2 of Lang's book, \cite{lang2}) we conclude that there is a generator $\Pi$ of $\mathfrak{p}_{K(z_0)}$ and a generator $\pi'$ of $\mathfrak{p}_K$ such that $\Pi^{q-1} = \pi'$. Note that $X^{q-1} - \pi'$ is an irreducible polynomial over $K$. Thus $K(\Pi)$ has degree $q-1$ over $K$ and $K(\Pi) = K(z_0)$. \\
The residue degree corresponding to the extension $K/\mathbb{Q}_p$ is $f$. Hence it contains all the $(p^f - 1)$-th roots of unity. Thus it contains all the $(q-1)$-th roots of unity (recall that $q = p^h$ and $h\,|\,f$). This along with previous observation implies that $K(z_0)/K$ is a Kummer extension with enough roots of unities in the base field. \\
Now from standard theory of Kummer extensions we conclude that $K(z_0)/K$ is Galois and $\text{Gal}(K(z_0)/K) \cong \mathbb{Z}/(q-1)\mathbb{Z}$.\\
This proves the lemma. \hfill $\square$\\~\\
\textbf{Corollary 3.3 :} i) Let $\Pi$ be as in the proof of lemma 3.2. Let $\sigma \in \text{Gal}(K(\pi)|K)$. By the proof of the lemma we know that there is an unique $\zeta \in \mu_{q-1}$ (depending on $\sigma$) such that $\sigma(\Pi) = \zeta \Pi$. Let $z \in \mathfrak{F}[\pi] - \{0\}$. Then $|\frac{\sigma(z)}{z} - \zeta|_p < 1$.\\
ii) If $\zeta$ is the root of unity associated to $\sigma$ then $\zeta^i$ is the root of unity associated to $\sigma^i$ and $|\frac{\sigma^i(z)}{z} - \zeta^i|_p < 1$ for all $i \geq 1$.\\
iii) If $\sigma$ is a generator $\text{Gal}(K(\pi)|K)$ then the corresponding $\zeta$ is a generator of $\mu_{q-1}$.\\~\\
\textbf{Proof :} $(i)$ Write $z = \sum_{i \geq 1} a_i\Pi^i$ where $a_i \in O_K$ for each $i \geq 1$ and $a_1 \in U_K$. This is possible since by corollary 2.7 and lemma 3.2 $z$ has order $1$ in $K(\pi)$.\\
Now $\sigma(z) = \sum_{i \geq 1} \zeta a_i \Pi^i$. Thus $\frac{\sigma(z)}{z} = \zeta (1+a)$ where $a \in O_{K(\pi)}$, the ring of integers of $K(\pi)$ with $|a|_p < 1$.\\ From here the proof of corollary is clear. \\
$(ii)$ Note that if $\sigma(\Pi) = \zeta\Pi$ then $\sigma^i(\Pi) = \zeta^i\Pi$ for all $i \geq 1$. Now the result follows from first part.\\
$(iii)$ This follows from the observation that if $\sigma^i \neq \sigma^j$ for some $i,j \geq 1$ then $\zeta^i \neq \zeta^j$. \hfill $\square$ \\~\\ 
\textbf{Corollary 3.4 :} Let $\sigma \in \text{Gal}(K(\pi)|K)$ and $z \in \mathfrak{F}[\pi] - \{0\}$. Then $\sigma(z) \in \mathfrak{F}[\pi] - \{0\}$.
Further, the correspondence $\text{Gal}(K(\pi)|K) \to \mathfrak{F}[\pi] - \{0\}$ given by $\sigma \to \sigma(z)$ is bijective.\\~\\
\textbf{Proof :} First part of the assertion follows from the fact that the group law is defined over $K$.\\
The correspondence is injective because $z$ generates $K(\pi)$ over $K$. Note that the cardinality of both sides is same. Hence an injective map actually induces an bijective correspondence. \hfill $\square$\\~\\
\textbf{Lemma 3.5 :} Let $z_1, z_2 \in \mathfrak{F}[\pi^n] - \mathfrak{F}[\pi^{n-1}]$. Then there is an unique $\zeta \in \mu_{q-1}$ such that $|\frac{z_1}{z_2} - \zeta|_p < 1$. Further the same $\zeta$ satisfies $|\frac{[\pi^{n-1}](z_1)}{[\pi^{n-1}](z_2)} - \zeta|_p < 1$.\\~\\
\textbf{Proof :} Note that if $\zeta_1, \zeta_2 \in \mu_{q-1}$ then 
\[|\zeta_1 - \zeta_2|_p = 
\begin{cases}
1 \,\text{if} \,\zeta_1 \neq \zeta_2 \\
0 \,\text{if} \,\zeta_1 = \zeta_2 .\\
\end{cases}
\]
From here using ultrametric triangle inequality uniqueness part follows.\\
We prove rest of the assertion by induction on $n$.\\
By corollary 3.4 there is a $\sigma \in \text{Gal}(K(\pi)|K)$ such that $\sigma(z_2) = z_1$. Let $\zeta \in \mu_{q-1}$ be element as in statement of corollary 3.3. Clearly this $\zeta$ satisfies the desired properties.\\
Assume that the statement holds for some positive integer $k$. We would like show that the statement is true for $k+1$.\\
Let $z_1, z_2 \in \mathfrak{F}[\pi^{k+1}] - \mathfrak{F}[\pi^k]$. Hence $|\frac{z_1}{z_2}|_p = 1$.\\
Let $[\pi](X) = \pi X + a_2X^2 + a_3X^3 + \cdots$. Note that $\pi | a_i$ for each $2 \leq i \leq q-1$ and $|a_q|_p = 1$. \\
Note that by corollary 2.7 $|z_1|_p = |\pi|_p^{\frac{1}{q^k(q-1)}}$.\\
Thus \[|z_1^q|_p = |\pi|_p^{\frac{1}{q^{k-1}(q-1)}} > |\pi|_p.\] Hence looking at the power series expansion we conclude that $\frac{[\pi](z_1)}{z_1^q} = a_q + A_1$ where $A_1 \in K(z_1)$ with $|A_1|_p < 1$.\\
Similarly $\frac{[\pi](z_2)}{z_2^q} = a_q + A_2$ where $A_2 \in K(z_2)$ with $|A_2|_p < 1$.\\
Using ultrametric triangle inequality $|\frac{[\pi](z_1)}{z_1^d} - \frac{[\pi](z_2)}{z_2^d}|_p < 1$. Note that \[\Big|\frac{[\pi](z_1)}{z_1^q} - \frac{[\pi](z_2)}{z_2^q}\Big|_p = \Big|\frac{[\pi](z_2)}{z_1^q}\Big|_p \Big|\frac{[\pi](z_1)}{[\pi](z_2)} - \frac{z_1^q}{z_2^q}\Big|_p = \Big|\frac{[\pi](z_1)}{[\pi](z_2)} - \frac{z_1^q}{z_2^q}\Big|_p \] where the last equality follows from corollary 2.7 and the fact that $[\pi](z_2) \in \mathfrak{F}[\pi^k] - \mathfrak{F}[\pi^{k-1}]$. Thus $|\frac{[\pi](z_1)}{[\pi](z_2)} - \frac{z_1^q}{z_2^q}|_p < 1$.\\
Now let $L = K(z_1,z_2)$. Assume that the degree of the residue extension of $L$ over $\mathbb{F}_p$ is $f_1$. Clearly $h | f_1$ and $\mu_{p^{f_1} - 1} \subseteq L$. Further $\mu_{p^{f_1} - 1} \cup \{0\}$ forms a complete set of representatives for residue extension corresponding to $L$. Since $|\frac{z_1}{z_2}|_p = 1$ there is a $\zeta_1 \in \mu_{p^{f_1} -1}$ such that $|\frac{z_1}{z_2} - \zeta_1|_p < 1$. Hence $|\frac{z_1^q}{z_2^q} - \zeta_1^q|_p < 1$. Note that \[\Big|\frac{[\pi](z_1)}{[\pi](z_2)} - \frac{z_1^q}{z_2^q}\Big|_p = \Big|\frac{[\pi](z_1)}{[\pi](z_2)} - \zeta + \zeta - \frac{z_1^q}{z_2^q}\Big|_p\] where $\zeta$ is the unique element in $\mu_{q-1}$ which satisfies $|\frac{[\pi^k](z_1)}{[\pi^k](z_2)} - \zeta|_p < 1$. By induction hypothesis $|\frac{[\pi](z_1)}{[\pi](z_2)} - \zeta|_p < 1$.\\~\\
\textbf{Claim 3.6 :} With notations as above $\zeta_1^q = \zeta$.\\
\textbf{Proof of claim :} Say, $\zeta_1^q \neq \zeta$. Then $|\zeta_1^q - \zeta|_p = 1$ and \[ |\frac{z_1^q}{z_2^q} - \zeta|_p = |\frac{z_1^q}{z_2^q} - \zeta_1^q + \zeta_1^q - \zeta|_p = 1\] since $|\frac{z_1^q}{z_2^q} - \zeta_1^q|_p < 1$. But then \[|\frac{[\pi](z_1)}{[\pi](z_2)} - \frac{z_1^q}{z_2^q}|_p = |\frac{[\pi](z_1)}{[\pi](z_2)} - \zeta + \zeta -\frac{z_1^q}{z_2^q}|_p = 1\].\\  
Contradiction !\\
This contradiction proves the claim. \hfill $\square$\\~\\
Hence $\zeta_1^{q(q-1)} = 1$. But order of $\zeta_1$ is coprime to $p$. Hence $\zeta_1^{q-1} = 1$. So $\zeta_1^q = \zeta_1$. Hence $\zeta_1 = \zeta$. \\
Thus there is a $\zeta \in \mu_{q-1}$ satisfying $|\frac{z_1}{z_2} - \zeta|_p < 1$ and this $\zeta$ satisfies $|\frac{[\pi^k](z_1)}{[\pi^k](z_2)} - \zeta|_p < 1$.
So the statement is true for $k+1$.\\
Hence we are done by the principle of mathematical induction. \hfill $\square$\\~\\
\textbf{Remark 3.7 :} This lemma gives a kind of analytic interpretation of the $R^{h_r}$ module structure that we want.\\~\\
\textbf{Lemma 3.8 :} Assume that $z_1, z_2 \in \mathfrak{F}[\pi^n] - \mathfrak{F}[\pi^{n-1}]$. Then there is a $\tau \in G_n$ such that $\tau(z_1) = z_2$.\\~\\
\textbf{Proof :} Follows from the fact that $z_1, z_2$ are roots of same irreducible polynomial over $K$ (corollary 2.5 $(ii)$). \hfill $\square$\\~\\
From now on we shall work with the assumption as stated in introduction and prove a sequence of lemmas which are analogues of lemmas in chapter 2. Proofs are more or less similar and we shall indicate if there is any significant difference.\\~\\
\textbf{Lemma 3.9 :} Let $G = \text{Gal}(K(\pi^{n+1})|K(\pi^n))$. Then $G \cong (\mathfrak{F}[\pi], \oplus_{\mathfrak{F}})$ as an abelian group.\\~\\
\textbf{Proof :} Similar to the proof of lemma 2.3.8. \hfill $\square$\\~\\
\textbf{Notation :} As before \[ G_{n} = \text{Gal}(K(\pi^n)|K)\] and \[G_{n,i} = G_{i}(K(\pi^n)|K).\] Since $K(\pi^n)|K$ is totally ramified we conclude that $G_{n,0} = G_n$.\\~\\
\textbf{Remark 3.10 :} Assume that $z \in \mathfrak{F}[\pi^n] - \mathfrak{F}[\pi^{n-1}]$ and let $0 \leq r \leq n$ and $\sigma \in G_n$.
Then $(\sigma(z) \ominus_{\mathfrak{F}} z) \in \mathfrak{F}[\pi^r]$ if and only if $\sigma \in \text{Gal}(K(\pi^n)|K(\pi^{n-r}))$. Proof is similar to the proof of remark 2.3.10.\\~\\
\textbf{Lemma 3.11 :} Let $z \in \mathfrak{F}[\pi^n] - \mathfrak{F}[\pi^{n-1}]$ and $\sigma \in G_n - \{\text{Id}\}$. Assume that $(\sigma(z) \ominus_{\mathfrak{F}} z) \in \mathfrak{F}[\pi^r] - \mathfrak{F}[\pi^{r-1}]$ for some $1 \leq r \leq n$ (note that $\sigma(z) \neq z$). Then $\sigma \in G_{n,i}$ for $0 \leq i \leq q^{n-r} - 1$ but $\sigma \notin G_{n,i}$ for $i = q^{n-r}$.\\~\\
\textbf{Proof :} Similar to the proof of lemma 2.3.11. \hfill $\square$\\~\\
\textbf{Proof of theorem 1.3(iii) :} Let $k$ and $i$ be integers with $1 \leq k \leq n$ and $q^{k-1} \leq i \leq q^{k}-1$ and assume that $z \in \mathfrak{F}[\pi^n] - \mathfrak{F}[\pi^{n-1}]$. Now by lemma 3.11 \[G_{n,i} = \{ \sigma \in G_n | (\sigma(z) \ominus_{\mathfrak{F}} z) \in \mathfrak{F}[\pi^{n-k}]\} = \text{Gal}(K(\pi^n)|K(\pi^k))\] where the last equality follows from remark 3.10.\\
This proves theorem 1.3(iii). \hfill $\square$\\~\\
\textbf{Lemma 3.12 :} $\text{Gal}(K(\pi^n)|K(\pi^{n-1}))$ lies in the center of $G_n$.\\~\\
\textbf{Proof :} Similar to the proof of lemma 2.3.12. \hfill $\square$\\~\\
We have already seen that $G_1 \cong \mathbb{Z}/(q-1)\mathbb{Z}$. Let $\widetilde{\tau}$ be a generator of this cyclic group.\\~\\
\textbf{Lemma 3.13 :} There exists a sequence $\{\tau_n\}_{n \geq 1}$ such that $\tau_n \in G_n$ and the following holds :\\
a) Restriction of $\tau_n$ to $K(\pi^n)$ is $\tau_m$ for all $1 \leq m \leq n$ and $\tau_1 = \widetilde{\tau}$.\\
b) Order of $\tau_n$ is $(q-1)$.\\
c) $\tau_n$ is in the center of $G_n$.\\~\\
\textbf{Proof :} Similar to the proof of lemma 2.3.13. \hfill $\square$\\~\\
\textbf{Remark 3.14 :} Put $\bigcup_{n \geq 1} K(\pi^n) = K(\pi^{\infty})$. This is an infinite Galois extension of $K$. The previous construction defines an element $\tau_{\infty} \in \text{Gal}(K(\pi^{\infty})|K)$ such that image of $\tau_{\infty}$ in $G_n$ is $\tau_n$. It is easy to check that $\tau_{\infty}^{q-1} = \text{Id}$ and $\tau_{\infty}$ lies in the centre of $\text{Gal}(K(\pi^{\infty})|K)$ (since this is true for restriction of $\tau_{\infty}$ to $K(\pi^n)$ for all $n \in \mathbb{N}$ one can verify by letting the automorphisms act on elements). Further, since order of $\widetilde{\tau}$ is exactly $(q-1)$ we conclude that order of $\tau_{\infty}$ is exactly $(q-1)$. We shall write $G_{\infty} = \text{Gal}(K(\pi^{\infty})|K)$.\\~\\
\textbf{Lemma 3.15 :} Let $z \in \mathfrak{F}[\pi^n] - \mathfrak{F}[\pi^{n-1}]$. Then $\{z,\tau_{\infty}(z),\cdots, \tau_{\infty}^{h_r - 1}(z)\}$ forms a $R/\pi^nR$ module base for $\mathfrak{F}[\pi^n]$.\\~\\
\textbf{Proof :} Since $\mathfrak{F}[\pi]$ is already a free $R/\pi^nR$ module of rank $h_r$ it is enough to show that the set $\{z, \tau_{\infty}(z), \tau_{\infty}^2(z), \cdots, \tau_{\infty}^{h_r-1}(z) \}$ is linearly independent over $R/\pi^nR$.\\
For convenience write $z_i = \tau_{\infty}^{i-1}(z)$ for each $1 \leq i \leq h_r$.\\ 
First consider the case $n = 1$ ie $z \in \mathfrak{F}[\pi] - \{0\}$. \\
Let $a_1,\cdots,a_{h_r} \in R$ be such that $[a_1](z_1) \oplus_{\mathfrak{F}} [a_2](z_2) \oplus_{\mathfrak{F}} \cdots \oplus_{\mathfrak{F}}[a_{h_r}](z_{h_r}) = 0.$\\
Note that $\mu_{p^{f_R}-1} \subseteq R$ and $\mu_{p^{f_R}-1} \cup \{0\}$ forms a set of representatives for the field of residue corresponding to $R$. Write \[a_i = \sum_{j \geq 0} c_{i,j}\pi^j \] where $1 \leq i \leq h_r$ and $c_{i,j} \in \mu_{p^{f_R} - 1} \cup \{0\}$.\\
Now by corollary 3.3 there is a primitive $(q-1)$-th root of unity $\zeta$ such that we have $z_i = \zeta^{i-1}z_1 + a_i$ for some $a_i$ satisfying $|a_i|_p < |z_1|_p$.\\ 
Note that $F(X,Y) = X + Y + \text{higher degree terms}$. Further, if $a \in R$ then $[a](X) = aX + \text{higher order terms}$.\\
Thus \[[a_1](z_1)\oplus_{\mathfrak{F}}\cdots\oplus_{\mathfrak{F}}[a_{h_r}](z_{h_r}) = (\sum_{1 \leq i \leq h_r}c_{i,0}\zeta^{i-1})z_1 + a' \] for some $a'$ satisfying $|a'|_p < |z_1|_p$ ($|z_1|_p = |z_j|_p < 1,\;\forall 1 \leq j \leq h_r$).  \\
Since the left hand side is zero one must have $|\sum_{1 \leq i \leq h_r} c_{i,0}\zeta^{i-1}|_p < 1$.\\
Clearly, $R^{h_r} = R[\zeta]$ (see \cite{Serre}, chapter 3, section 6). Let $k_{R^{h_r}}$, $k_R$ be the field of residues corresponding to $R^{h_r}$ and $R$ respectively. Note that under reduction modulo $\pi$ map the set $\mu_{p^{f_R}-1} \cup \{0\}$ maps bijectively to $k_R$ and the image of the set $\{1, \zeta, \cdots, \zeta^{h_r - 1} \}$ forms a $k_R$ basis for $k_{R^{h_r}}$. Since $|\sum_{1 \leq i \leq h_r} c_{i,0}\zeta^{i-1} |_p < 1$ the element $\sum_{1 \leq i \leq h_r} c_{i,0}\zeta^{i-1}$ (of $R^{h_r}$) maps to $0$ under reduction modulo $\pi$. Using the linear independence of the image of $\{1, \zeta, \cdots, \zeta^{h_r - 1} \}$ we conclude that image of $c_{i,0}$ is $0$ under reduction modulo $\pi$ map, for all $1 \leq i \leq h_r$. Hence $c_{i,0} = 0$ for all $1 \leq i \leq h_r$. Hence $ \pi | a_{i}$ in $R$  for all $ 1 \leq i \leq h_r$.\\
We have shown that if $[a_1](z_1) \oplus_{\mathfrak{F}} \cdots \oplus_{\mathfrak{F}} [a_r](z_r) = 0$ then $a_i \in \pi R$ for all $1 \leq i \leq h_r$. This shows that $\{z_1, \cdots, z_{h_r}\} $ is linearly independent over $R/\pi R$.\\
Thus we are done in the case $n = 1$.\\
Now we shall consider the general case ie $z \in \mathfrak{F}[\pi^n] - \mathfrak{F}[\pi^{n-1}]$ for some $n \in \mathbb{N}$.\\
First we prove the following claim :\\~\\
\textbf{Claim 3.16 :} Let $a_1,\cdots,a_{h_r} \in R$ be such that at least one of $a_i$-s is an unit. Then $[a_1](z_1) \oplus_{\mathfrak{F}} \cdots \oplus_{\mathfrak{F}}[a_{h_r}](z_{h_r}) \in \mathfrak{F}[\pi^n] - \mathfrak{F}[\pi^{n-1}]$.\\~\\
\textbf{Proof of claim :} Put $z = [\pi^{n-1}]([a_1](z_1) \oplus_{\mathfrak{F}} \cdots \oplus_{\mathfrak{F}}[a_{h_r}](z_{h_r})).$ \\
So $z = [a_1]([\pi^{n-1}](z_1)) \oplus_{\mathfrak{F}} \cdots \oplus_{\mathfrak{F}}[a_{h_r}]([\pi^{n-1}](z_{h_r}))$. Note that $[\pi^{n-1}](z_i) = \tau_{\infty}^{i-1}([\pi^{n-1}](z_1))$ for all $1 \leq i \leq h_r$. From our proof of the case $n = 1$ and the hypothesis of the claim it follows that $z \neq 0$. Clearly $[\pi](z) = 0$.\\
Hence $ z \in \mathfrak{F}[\pi] - \{0\}$. The claim follows from here. \hfill $\square$\\~\\
Now $a_1,\cdots,a_{h_r} \in R$ be such that $[a_1](z_1) \oplus_{\mathfrak{F}} \cdots \oplus_{\mathfrak{F}} [a_{h_r}](z_{h_r}) = 0.$  Write $a_i = \pi^{n_i}u_i$, for all $1 \leq i \leq h_r$ where $u_i$ is an unit in $R$. Put $m = \min \{n_1, \cdots, n_{h_r} \}$. Now $[a_1](z_1) \oplus_{\mathfrak{F}} \cdots \oplus_{\mathfrak{F}} [a_{h_r}](z_{h_r}) = [\pi^m] ([a_1\pi^{-m}](z_1) \oplus_{\mathfrak{F}} \cdots \oplus_{\mathfrak{F}}  [a_{h_r}\pi^{-m}](z_{h_r})).$ Now by claim 3.16 the expression within bracket $\in \mathfrak{F}[\pi^n] - \mathfrak{F}[\pi^{n-1}].$ Since $[a_1](z_1) \oplus_{\mathfrak{F}} \cdots \oplus_{\mathfrak{F}} [a_{h_r}](z_{h_r}) = 0$ we conclude that $n \leq m$.\\This shows that if $[a_1](z_1) \oplus_{\mathfrak{F}} \cdots \oplus_{\mathfrak{F}}[a_{h_r}](z_{h_r}) = 0$ then $a_i \in \pi^nR$ for all $1 \leq i \leq h_r$. This finishes the proof of lemma. \hfill $\square$ \\~\\
Let $T(\pi)$ be the $\pi$-adic Tate module. It has a natural $R$ module structure. Let $z \in T(\pi)$ be such that its image $z_n$ in $\mathfrak{F}[\pi^n]$ satisfies $[\pi^{n-1}](z_n) \neq 0$. It is clear that such an element exists.\\~\\
\textbf{Corollary 3.17 :} $T(\pi)$ is a free $R$ module of rank $h_r$ and $\{z, \tau_{\infty}(z), \cdots, \tau_{\infty}^{h_r - 1}(z) \}$ forms a $R$ base of $T(\pi)$.\\~\\
\textbf{Proof :} From lemma 3.15 $\{z_n, \tau_{\infty}(z_n), \cdots, \tau_{\infty}^{h_r - 1}(z_n) \}$ generates $\mathfrak{F}[\pi^n]$ as $R$ module, for each $n \geq 1$. Hence the set under consideration generates $T(\pi)$ as $R$ module.\\
Now let $a_1, \cdots, a_{h_r} \in R$ be such that $[a_1](z) \oplus_{\mathfrak{F}} \cdots \oplus_{\mathfrak{F}} [a_{h_r}](\tau_{\infty}^{h_r - 1}(z)) = 0$. Using lemma 3.15 we conclude that $a_i \in \pi^{n}R$ for all $1 \leq i \leq h_r$ and for all $n \geq 1$. This proves that $a_i = 0$ for all $1 \leq i \leq h_r$. \hfill $\square$ \\~\\ 
Consider $\tau_{\infty} : T(\pi) \to T(\pi)$. Since the group law is defined over $K$ we conclude that $\tau_{\infty}$ is a $R$ module morphism. But $\tau_{\infty}$ has order $(q-1)$ in $G_{\infty}$. Hence $\tau_{\infty}^{q-1} = \text{Id}$ in $\text{End}_R(T(\pi))$. \\
Let $\zeta$ be a primitive $(q-1)$-th root of unity. We have $R^{h_r} = R[\zeta]$. The argument above shows that we have an well defined $R$-algebra homomorphism $R[\zeta] \to \text{End}_R(T(\pi))$ by sending $\zeta \to \tau_{\infty}$ and extending the map $R$-linearly.\\
This puts a $R^{h_r} = R[\zeta]$ module structure on $T(\pi)$ which commutes with $R$ module structure.\\
Now we have the following lemma :\\~\\
\textbf{Lemma 3.18 :} $T(\pi)$ is a $R[\zeta]$ module of rank $1$.\\~\\
\textbf{Proof :} Let $z$ be as before. We shall show that $T(\pi)$ is a free $R[\zeta]$ module with $\{z\}$ as a base.\\
Note that $\{z, \cdots, \tau_{\infty}^{h_r - 1}(z)\} \subseteq R[\zeta](\{z\})$. By corollary 3.17 we conclude that $T(\pi) \subseteq R[\zeta](\{z\})$. Hence $T(\pi)$ is generated by $z$ as $R[\zeta]$ module.\\
Note that $R[\zeta] = R \oplus R\zeta \oplus \cdots \oplus R\zeta^{h_r-1}$ as $R$ module. Now let $a \in R[\zeta]$ be such that $[a](z) = 0$. Write $a = a_1 + a_2\zeta + \cdots + a_{h_r}\zeta^{h_r - 1}$ with $a_1, a_2, \cdots a_{h_r} \in R$. From definition of $R[\zeta]$ module structure we have $[a](z) = [a_1](z) \oplus_{\mathfrak{F}} [a_2](\tau_{\infty}(z)) \oplus_{\mathfrak{F}} \cdots \oplus_{\mathfrak{F}} [a_{h_r}](\tau_{\infty}^{h_r - 1}(z))$. Using corollary 3.17 we conclude that $a_i = 0$ for all $1 \leq i \leq h_r$. This shows that $a = 0$.\\
This proves that $z$ freely generates $T(\pi)$ as $R[\zeta]$ module. \hfill $\square$ \\~\\
\textbf{Lemma 3.19 :} $\text{Gal}(\overline{\mathbb{Q}}_p|K)$ acts on $T(\pi)$ via $R[\zeta]$ module morphisms.\\~\\
\textbf{Proof :} Follows from the fact that $\mathfrak{F}$ is defined over $K$ and $\tau_{\infty}$ lies in the center of $\text{Gal}(K(\pi^{\infty})|K)$. \hfill $\square$ \\~\\ 
Note that $R^{h_r}$ module structure on $T(\pi)$ gives a $R^{h_r}/\pi^nR^{h_r}$ module structure on $\mathfrak{F}[\pi^n]$ for each $n$.\\
Let $z$ be as before and $z_n$ be its image in $\mathfrak{F}[\pi^n]$. We have :\\~\\
\textbf{Lemma 3.20 :} $z_n$ freely generates $\mathfrak{F}[\pi^n]$ as $R^{h_r}/\pi^nR^{h_r}$ module.\\~\\
\textbf{Proof :} Similar to the proof of lemma 3.18. \hfill $\square$\\~\\
Let $\sigma \in \text{Gal}(K(\pi^n)|K)$. Let $a \in R^{h_r}/\pi^nR^{h_r}$ be such that $\sigma(z_n) = [a](z_n)$. Since $\sigma$ is invertible $a$ must be an unit.\\
Note that the unit group of $R^{h_r}/\pi^nR^{h_r}$ is $U_{0,K^{h_r}_R}/U_{n,K^{h_r}_R}$. Hence we have a map $\phi : \text{Gal}(K(\pi^n)|K) \to U_{0,K^{h_r}_R}/U_{n,K^{h_r}_R}$. Since the action of Galois group commutes with $R^{h_r}$ module structure it is easy to see that $\phi$ is a group homomorphism.\\
$\phi$ is injective since $z_n$ generates $K(\pi^n)$ over $K$. Since both the sides have same cardinality it is also onto.\\
This proves theorem 1.3 (ii). \hfill $\square$  

\chapter[Errata and concluding remarks]{Errata and concluding remarks\footnote{Added later}}
Statement of lemma 2.3.7 (see chapter-2 and appendix A) does not hold in general. In fact, in the set-up considered in chapter-2 this assumption turns out to be equivalent to theorem 2.1.1.\\
For more details see \cite[Sec 2]{article1}. The question of Galois group associated to formal group laws is further studied in \cite{article2}, \cite{article3}. \\~\\
Set up be as in chapter-2.\\
Let $\mathfrak{F}$ denote the formal group law associated to $E$. Clearly $\mathfrak{F}$ is defined over $\mathbb{Z}_w$ and in terminology of $\cite{article1}$ it is an unramified $\mathbb{Z}_p$ module of height $2$ (here $w = p^{2f}$). All $p$-torsion points of $E$ come from $p$-torsion points $\mathfrak{F}$. Since $\mathfrak{F}$ is an unramified group law, $\text{End}_{\mathbb{Z}_w}(\mathfrak{F})$ is an integrally closed, complete subring of $\mathbb{Z}_w$ (\cite[Thm 3.1.2]{article1}). Since $\mathfrak{F}$ has $\mathbb{Z}_p$ height 2,  $\text{End}_{\mathbb{Z}_w}(\mathfrak{F})$ equals the absolute endomorphism ring  $\text{End}(\mathfrak{F})$ (\cite[Sec 2]{article2}) and there are only two possibilities for endomorphism rings :\\     
I. $\text{End}_{\mathbb{Z}_w}(\mathfrak{F}) = \mathbb{Z}_p$,\\
II. $\text{End}_{\mathbb{Z}_w}(\mathfrak{F}) = \mathbb{Z}_q$.\\
If case II holds then lemma 2.3.7 is true (\cite[Thm 2.3]{article1}) and we have already proved theorem 0.1 and theorem 0.2 in this case.\\
In following discussion we shall assume case I holds.\\
$T_p(\mathfrak{F})$ be the $p$-adic Tate module associated to $\mathfrak{F}$. It is a free $\mathbb{Z}_p$ module of rank $2$. Consider the usual representation \[\rho : G_{w} \to \text{Gl}(V_p(\mathfrak{F}))\] where $G_w = \text{Gal}(\overline{\mathbb{Q}}_w|\mathbb{Q}_w)$ and $V_p(\mathfrak{F}) = T_p(\mathfrak{F}) \otimes_{\mathbb{Z}_p} \mathbb{Q}_p$. \\ Let $H = \text{Im}(\rho).$ Clearly there is a topological isomorphism of groups \[ H \cong \text{Gal}\big(\mathbb{Q}_w(\mathfrak{F}[p^\infty])|\mathbb{Q}_w\big).\] 
This representation was studied in \cite{article1}, \cite{article2}, \cite{article3} and ramification filtration on $\text{Gal}(\mathbb{Q}_w(\mathfrak{F}[p^\infty])|\mathbb{Q}_w)$ can be given by estimates in Sen's theory of $p$-adic Lie filtration (see \cite[A.4]{galoisrep} and \cite[Sec 2]{article3}). But these estimates are not sufficient to prove results of chapter-3.\\~\\ 
Note that if the elliptic curve $E$ in chapter 2 is defined over $\mathbb{Z}_q$ one can directly appeal to Lubin-Tate theory and twisting of elliptic curves as in Habegger's paper (\cite[Sec 3]{hab}) and prove lemma 2.3.7 using Lubin-Tate theory. In global setup of chapter 3 the minimal Weierstrass model of $E$ at $\mathfrak{p}$ is defined over $\mathbb{Z}_q$ if degree of residue extension of $\mathfrak{p}$ is $1$ or $2$ ie. $f = 1,2$. \\
Thus theorem 0.1 and theorem 0.2 holds if there exists a prime $\mathfrak{p}$ such that it satisfies three conditions in section 3.1 (chapter-3) and at least one of the following conditions is satisfied: \\
I. $\text{End}_{\mathbb{Z}_w}(\mathfrak{F}) = \mathbb{Z}_q$,\\
II. $f = 1,2$. \\
Proving existence of such a prime is associated with questions related to distribution of supersingular primes and we avoid making any precise remarks about this.

\end{appendices}


\begin{thebibliography}{21}

\bibitem{hab}
P.Habegger;
\textit{Small Height and Infinite Nonabelian Extensions;}
Duke Mathematical Journal, Vol. 162, No. 11, 2013.

\bibitem{sil1}
J.H.Silverman;
\textit{Arithmetic of Elliptic Curves;}
Grad. Texts in Math, Springer.

\bibitem{sil 2}
J.H.Silverman;
\textit{Advanced Topics in the Arithmetic of Elliptic Curves;}
Grad. Texts in Math, Springer.

\bibitem{sil 3}
J.H.Silverman;
\textit{A lower bound for the cannonical height on elliptic curves over abelian extensions;}
J. Number Theory, 2004.

\bibitem{az}
F.Amoroso and U.Zannier;
\textit{A relative Dobrowolski lower bound over abelian extensions;}
Ann. Scuola Norm. Sup. Pisa Cl. Sci.(2000)

\bibitem{sch}
A. Schinzel;
\textit{On the product of conjugates outside the unit circle of an algebraic number;}
Acta Arith. (1973)

\bibitem{heights}
E.Bombeiri and W. Gubler;
\textit{Heights in Diophantine Geometry;}
Cambridge University Press

\bibitem{Serre}
Jean-Pierre Serre;
\textit{Local Fields;}
Springer

\bibitem{Elkies}
N.Elkies;
\textit{Supersingular primes for elliptic curves over real number fields;}
Composito Math. (1989)

\bibitem{oshikawa}
K Oshikawa;
\textit{On formal groups over complete discrete valuation rings with application to a theorem of Lutz;}
Journal f\"ur die reine und angewandte Mathematik 334 (1982): 79-90

\bibitem{Huse}
D.Husem\"oller;
\textit{Elliptic curves, 2nd edition;}
Grad. Texts in Math, Springer.

\bibitem{lang}
S.Lang;
\textit{Algebra;}
Grad. Texts in Math., Springer.

\bibitem{lang2}
S. Lang;
\textit{Algebraic Number Theory}
Grad. Texts in Math., Springer.

\bibitem{wiki}
Wikipedia
\textit{Kummer's theorem,}
https://en.wikipedia.org/

\bibitem{frey}
Linda Frey;
\textit{Explicit Small Heights in Infinite Non-Abelian Extensions;}
https://arxiv.org/pdf/1712.04214.pdf

\bibitem{vou}
Paul M Voutier;
\textit{An effective lower bound for the height of algebraic numbers;}
https://arxiv.org/pdf/1211.3110.pdf

\bibitem{neu}
J\"urgen Neukirch;
\textit{Class Field Theory;}
Springer.

\bibitem{Bourbaki}
Nicholas Bourbaki;
\textit{Commutative Algebra;}
Addision-Wesley.


\bibitem{haz2}
Michiel Hazewinkel;
\textit{Formal Groups and Applications;}
Academic Press, 1978.

\bibitem{article1}
Soumyadip Sahu;
\textit{Field Generated by Division Points of Certain Formal Group Laws;}
2018, arxiv:1809.00112.

\bibitem{article2}
Soumyadip Sahu;
\textit{Field Generated by Division Points of Certain Formal Group Laws II;}
2018, arxiv:1810.00330.

\bibitem{article3}
Soumyadip Sahu;
\textit{Field Generated by Division Points of Certain Formal Group Laws III;}
2018; arXiv:1811.04855.

\bibitem{galoisrep}
Jean-Marc Fontaine and Yi Ouyang;
\textit{Theory of p-adic Galois Representations;}
Springer. 

\end{thebibliography}
\end{document}